\documentclass[12pt]{article}

\usepackage{amsfonts}

\def\a{\mathfrak{a}}

\def\sii{\underlineÊ{\si}}

\def\s{\mathfrak{s}}
\def\t{\mathfrak{t}}

\def\si{\sigma}
\def\ep{\varepsilon}

\def\A{{\cal A}}
\def\CC{{\cal C}}
\def\O{\cal O}

\def\LL{\Lambda}
\def\l{\lambda}
\def\aa{\alpha}

\def\cu{C^{\infty}(U_0\backslash G, \psi)}
\def\CCU{ {\cal C} (U_0\backslash G, \psi)}
\def\ccu{C_c^{\infty}(U_0\backslash G, \psi)}
\def\cmu{C^{\infty}(U_0\cap M \backslash M, \psi)}

\def\CCMU{ { \cal C} (U_0\cap M \backslash M, \psi)}
\def\lu{L^2(U_0\backslash G, \psi)}
\def\AA{{\cal A} ^w(U_0\backslash G, \psi)}
\def\AAM{{\cal A} ^w(U_0\cap M\backslash M, \psi)}

\def\fhi{C^{\infty}({\cal O} ,Wh \otimes  i_P^G)}

\def\sgpss{sous-groupe parabolique semi-standard}

\def\DD{\Delta}

\def\sii{{\underline \sigma}}

\def\D{\mathbb{D}}

\def\N{\mathbb{N}}

\def\Z{\mathbb{Z}}
\def\R{\mathbb{R}}
\def\C{\mathbb{C}}
\def\D{\mathbb{D}}

\def\fd {\hspace{0.35cm} \raise
-0.5mm\hbox{$\blacksquare$}\\}

\def\qed{{\null\hfill\ \raise3pt\hbox{\framebox[0.1in]{}}}\break\null}
\newtheorem{theo}{Th\'eor\`eme}
\newtheorem{prop}{Proposition}
\newtheorem{lem}{Lemme}
\newtheorem{cor}{Corollaire}
\newtheorem{rem}{Remarque}
\newtheorem{defi}{D\' efinition}

\def\ste{\par\smallskip\noindent}

\def\dem{ {\em D\'emonstration: \ste }}

\def\beq{\begin{equation}}
\def\eeq{\end{equation}}

\def\A {\cal A}
\def\D{\cal D}

\def\O {\mathcal O}

\newenvironment{res}
               {\begin{equation}\begin{minipage}{0.85\textwidth}}
               {\end{minipage}\end{equation}}
\def\ber{\begin{res}}
\def\eer{\end{res}}

\catcode`\Ž=\active\def Ž{\'e}
\catcode`\=\active\def {\`e}
\catcode`\ˆ=\active\def ˆ{\`a}
\catcode`\=\active\def {\`u}
\catcode`\=\active\def {\^{e}}
\catcode`\"=\active\def "{\^{i}}
\catcode`\™=\active\def ™{\^{o}}
\catcode`\ž=\active\def ž{\^{u}}
\catcode`\=\active\def {\c{c}}
 \def\opto#1{{\smash{\mathop{\longrightarrow}\limits^#1}}}

\title{Formule de Plancherel pour les fonctions de  Whittaker sur un groupe rŽductif $p$-adique }

\author{ Patrick Delorme}
\date{ }

\begin{document}

\maketitle
\section{Introduction}
Soit $F$ un corps local non archimŽdien et soit $G$ le groupe des points sur $F$ d'un groupe rŽductif connexe dŽfini sur $F$. Soit  $(P_0, P_0^-)$ un couple de sous-groupes paraboliques minimaux opposŽs de $G$. On note $M_0$ leur sous-groupe de LŽvi commun.  Soit $U_0$ le radical unipotent de $P_0$ et $A_0$ le plus grand  tore dŽployŽ de $M_0$. Soit $K$ un bon sous-groupe compact maximal  de $G$ relativement ˆ $A_0$.  Soit $\psi$  un caractre unitaire lisse de $U_0$ non dŽgŽnŽrŽ, i.e. tel que pour toute    racine $P_0$-simple de $A_0$, $\aa$,  sa restriction au  sous-groupe radiciel $({U_0})_\alpha$ soit non triviale.  
 \\On note $\cu$ l'espace des fonctions de Whittaker lisses sur $G$, i.e. des fonctions, $f$, sur $G$,  invariantes ˆ droite par un sous-groupe compact ouvert de $G$ et telles que:
  $$f(u_0g)= \psi(u_0) f(g), g\in G, u_0 \in U_0.$$
  On dŽfinit un sous-espace vectoriel, $\CCU$, de $\cu$, qui est l'analogue pour notre situation de l'espace de Schwartz ${\cal C}(G)$ d'Harish-Chandra (cf. [W]) et qui est muni d'une topologie naturelle. On dŽfinit une transformation de Fourier des ŽlŽments de $\CCU$, dite transformation de Fourier-Whittaker, ˆ l'aide des reprŽsentations lisses irrŽductibles de carrŽ intŽgrable de sous-groupes de LŽvi de $G$. A noter que celle-ci est diffŽrente de celle introduite dans [D4], qui utilise des reprŽsentations cuspidales .\\ Le principal rŽsultat est la caractŽrisation de l'image de cette transformation et une formule d'inversion. L'espace $\CCU$ est un sous espace de $\lu$. On Žtudie Žgalement  comment se transforme le produit scalaire $L^2$ par la transformation de Fourier-Whittaker. 
 Nous avons appris  de Ramakrishnan, par l'intermŽdiaire de Laurent Clozel, que ce travail ( la formule de Plancherel pour les fonctions de Whittaker) Žtait du ˆ Harish-Chandra, malheureusement non publiŽ.
 Noter que pour les fonctions de Whittaker sur un groupe rŽductif rŽel, la formule de Plancherel a ŽtŽ Žtablie (Harish-Chandra, non publiŽ,  Wallach [Wall], Chapitre 15).  
 Notre travail  est une  suite naturelle  ˆ [D4]. Il est largement inspirŽ par la rŽdaction de Waldspurger de la formule de Plancherel d'Harish-Chandra pour $L^2(G)$ ( cf. [W]). \\
 DŽtaillons le contenu de cet article.\\  
 Si $(\pi, V)$ est une reprŽsentation lisse de $G$, on note $Wh(\pi, U_0)$ ou $Wh(\pi)$ l'espace des formes linŽaires, $\xi$,  sur $V$ telles que:
$$ \langle \xi, \pi(u_0)v \rangle  = \psi(u_0)  \langle \xi, v\rangle  , v\in V, u_0\in U_0.$$
Si $\pi$ est de longueur finie,  $Wh(\pi) $ est de dimension finie (cf.[BuHen], ThŽorme 4.2) et [D3] pour une autre dŽmonstration).   Notez que si le groupe $G$ n'est pas quasi-dŽployŽ, cette dimension n'est pas toujours infŽrieure ou Žgale ˆ  1. 
\\ Si $v\in V$, on note $c_{\xi, v}$ le coefficient gŽnŽralisŽ dŽfini par:
$$ c_{\xi, v}(g):=  \langle \xi, \pi(g) v\rangle  , \>g \in G.$$
C'est un ŽlŽment de $\cu$. \\
 Dans la suite la phrase ``Soit $P= MU$ un sous-groupe parabolique semi-standard de $G$'' voudra dire que $P$ contient $M_0$, que $M$ est son sous-groupe de LŽvi contenant $M_0$ et que $U$ est son radical unipotent. On note 
$X(M)$ le groupe des caractres non ramifiŽs unitaires  de $M$. C'est un tore compact. On note $\delta_P$ la fonction module de $P$. Soit    
$(\si, E)$ une  reprŽsentation lisse de  $M$. On note
$(i^G_P  \si, i^G_P E)$ l'induite parabolique de $(\si, E)$. On suppose $\si$ unitaire irrŽductible et on note $\O$
 l'ensemble des classes d'Žquivalences des reprŽsentations 
 $\si_\chi:= \si\otimes \chi$,  $ \chi \in X(M)$, qui est un tore compact contenu dans un tore complexe $\O_\C$. On appelle ${\cal O}$ l'orbite inertielle unitaire de $\si$, qui est munie d'une mesure non nulle et $X(M)$-invariante, convenablement normalisŽe.  On utilise dans la suite la notion de fonction sur $\O$, qui dŽpendent des objets concrets $\si_\chi$, avec des rgles de transformation pour tenir compte des Žquivalences unitaires de reprŽsentations. On note que 
$Wh(\si _\chi ):= Wh(\si_ \chi, U_0\cap M)$ est indŽpendant de $\chi \in X(M) $, car $\chi$  est trivial sur  $ U_0\cap M$. Par restriction des fonctions ˆ $K$, les reprŽsentations $i^G_P  \si_\chi$ admettent une rŽalisation dans un espace indŽpendant de $\chi$, la rŽalisation compacte.
 \\{\em Il existe une bijection naturelle, $\eta \mapsto \xi(P, \si, \eta) $, entre $Wh(\si)$ et $Wh(i^G_P \si)$. (cf. Rodier [R], Casselman-Shalika [CS], Shahidi [Sh], Proposition 3.1)}\\  
On appelle fonctionnelle de Jacquet les ŽlŽments de $Wh(i^G_P \si)$. 
On dŽfinit les fonctionnelles de Jacquet  pour un sous-groupe parabolique semi-standard de $G$,  $P=MU$, par 
 transport de structure. \ste Plus prŽcisŽment,  soit $K$ un bon sous-groupe compact maximal  de $G$ relativement ˆ $A_0$. On note $\overline{ W}^G$ (resp. $\overline{ W}^M$) le groupe de Weyl de $G$ (resp. $M$) relativement ˆ $M_0$.  On fixe un ensemble $W^G$ de reprŽsentants dans $K$  de  $\overline{ W}^G$. Pour un bon choix de $w\in W^G$  
 tel que $Q:= wPw^{-1}$ est anti-standard, on dŽfinit 
 $$ \xi(P, \si, \eta): = \xi (Q, w\si, \eta) \circ \lambda(w),  \eta \in  Wh(P, \si):= Wh(w\si),$$
o $\lambda (w) $ est la translation ˆ gauche par $w$, qui entrelace $i^G_P\si$ et $i^G_Qw\si$.
Soit $P$,  $Q$  des sous-groupes paraboliques semi-standard de $G$,   possŽdant  le m\^eme sous-groupe de LŽvi,  $M$, contenant $M_0$.\ste On suppose $\si$ de longueur finie. On introduit les intŽgrales d'entrelacement, $A(Q, P, \si)$, qui, lorsqu'elles sont dŽfinies, entrelacent  $i^G_P \si$ et  $i^G_Q \si$. \\  
Les intŽgrales d'entrelacement transforment les fonctionnelles de Jacquet en des fonctionnelles de Jacquet, ce qui permet d'introduire les matrices $B$: \\
{\em Il existe une unique fonction rationnelle sur $X(M)$ ˆ valeurs dans $End (Wh (Q, \si), Wh (P, \si))$, $\chi \mapsto B(P, Q, \si_\chi)$,  telle que:}
$$\xi(Q, \si_\chi, \eta) \circ A(Q, P, \si_\chi)= \xi (P,\si_\chi,  B(P, Q, \si_\chi) \eta). $$ 
On dŽfinit les intŽgrales de Jacquet pour $\phi \in Wh(\si) \otimes i^G_PE$ par: $$ E^G_P(\phi) = c_{\xi, v}  \in  \cu,  $$ si $\phi= \eta \otimes v$, avec $\eta \in Wh(\si)$, $v\in i^G_P E$ et 
 $\xi= \xi(P, \si, \eta)$.\\
Soit $ (\pi, V)$ {une reprŽsentation lisse irrŽductible et de carrŽ intŽgrable de   $G$. On note $A_G$ le plus grand tore dŽployŽ du centre de $G$.  Alors pour tout $\xi \in Wh( \pi)$ et $v \in V$, $\vert c_{\xi, v}\vert$ est de carrŽ intŽgrable sur $ A_G U_0 \backslash G $. Gr\^ace au Lemme de Schur, on voit  qu'il  existe un unique produit scalaire sur $Wh(\pi)$ tel que:
$$\int _{A_G U_0\backslash G}  c_{\xi, v}(g) \overline{ c_{\xi', v'}(g)    }dg= (\xi, \xi') (v, v'), \>\xi, \xi' \in Wh(\pi), v,v' \in V.$$
On dŽfinit maintenant la transformŽe de Fourier-Whittaker de  $f\in \CCU$, ${\cal F}f$,   comme suit. Pour tout sous-groupe parabolique anti-standard de $G$, i.e. contenant $P_0^-$,  $P=MU$, et toute reprŽsentation lisse irrŽductible de carrŽ intŽgrable  de $M$,  $(\si, E)$, le ThŽorme de reprŽsentation de Riesz montre qu'il existe un  unique ŽlŽment de $ Wh(\si)\otimes i^G_PE$,   ${\cal F}f(P, \si)$,  tel que:  $$(({\cal F}f)(P, \si), \eta \otimes  v) = \int_{U_0\backslash G } f(g) \overline {E^G_P(\si,  \eta ,  v) (g)}dg, \eta \in Wh(P, \si), v\in i^G_PE,$$
l'intŽgrale Žtant convergente d'aprs les propriŽtŽs des intŽgrales de Jacquet et d'aprs l'hypothse $f \in \CCU$. \\
 Alors $ {\cal F}f $  vŽrifie les propriŽtŽs suivantes:\ste
 1) a) L'application $\chi \mapsto ({\cal F}f)(P, \si_\chi) $,  est une fonction $C^{\infty}$. En particulier, dans la rŽalisation compacte, cette application est ˆ valeurs dans un espace vectoriel de dimension finie.\\   
  b) Si $(\si, E)$ et $(\si_1, E_1) $ sont unitairement Žquivalentes, $({\cal F}f)(P, \si)$ et $({\cal F}f)(P, \si_1)$ vŽrifient une relation de  transport de structure.\\
  c) On peut dŽfinir pour $g\in G$, $\rho_\bullet (g) ({\cal F}f)$, en posant $(\rho_\bullet (g)( {\cal F}f ))(P, \si) = (Id \otimes i^G_P \si(g)) ( {\cal F}f )(P, \si)$. 
  Alors, si $f$ est invariante ˆ droite par   un sous-groupe compact ouvert de $G$, $H$, on a $\rho_{\bullet}(h) {\cal F}f =  {\cal F}f $ pour tout $h\in H$. \\
  On rŽsume les propriŽtŽs a), b), c) en disant que pour toute orbite inertielle unitaire, $\O$,   d'une  reprŽsentation lisse irrŽductible et de carrŽ intŽgrable de $M$,  $(\si, E)$, ${\cal F} (P, .)$ dŽfinit un ŽlŽment de 
  $C^{\infty}(\O_u, Wh(P, )\otimes   i^G_P)$, cet espace Žtant muni d'une topologie naturelle. Alors $\rho_\bullet $ dŽfinit une reprŽsentation lisse de $G$ sur cet espace. 
 On note $\Theta$ l'ensemble des couples $(P,\O)$ o $P=MU$ est un sous-groupe parabolique anti-standard de $G$ et $\O$ est l'orbite inertielle unitaire d'une reprŽsentation lisse  irrŽductible et de carrŽ intŽgrable de $M$. 
 Soit $(P, \O)\in \Theta$ et $P^-$ le sous-groupe parabolique opposŽ ˆ $P$ relativement ˆ $M$.
 \\ {\em Il existe une fonction $\mu$, rationnelle, non identiquement nulle, et  $C^{\infty} $ sur $ \O$, telle que l'on ait l'ŽgalitŽ de fonctions rationnelles sur $\O$:
 $$A(P,P^-, \si) A(P^-,P, \si) = \mu(\si)^{-1}Id.$$
 Soit   $\phi \in C^{\infty}(\O_u, Wh(P,) \otimes  i^G_P)  $.  On dŽfinit le paquet d'ondes $f_\phi$ par: 
$$f_\phi(g):= \int_{\O}\mu(\si) E(P, \si, \phi(\si))(g)d\si, g \in G,$$
On dŽmontre que $f_\phi\in \CCU$ en utilisant un rŽsultat semblable de [W] et le fait suivant, utilisŽ de manire rŽcurrente:
\\{\em  Si $(\pi,V)$ est une reprŽsentation lisse de $G$, $v\in V, \xi\in Wh(\pi).$ Alors la restriction de $c_{\xi,v}$ ˆ $A_0^-$ est Žgale  ˆ  la restriction ˆ $A_0^-$ d'un coefficent lisse, que l'on peut prŽciser.\\
 On dŽfinit: $${\cal S}= \oplus_{(P, \O) \in \Theta} \fhi.$$
  Soit $P= MU$ et $P'= M'U'$ deux sous-groupes paraboliques semi-standard de $G$.  On note note ${\overline W}(M'\vert G \vert M) $ l'ensemble quotient, par l'action ˆ gauche de  ${\overline  W}^{M'}$, de
   $\{ s\in  {\overline  W}^G \vert s. M \subset M'\}$.  On note $W(M'\vert G \vert M)$ un ensemble de reprŽsentants dans $W^G$ de ce quotient, dont nous ne dŽtaillerons pas le choix dans cette introduction. 
 Soit $(P, \O)\in \Theta$ et  $P'=M'U'$ un sous-groupe parabolique anti-standard de $G$,  $s\in W(M'\vert G\vert M)$ tels que $M$ et $M'$ soient conjuguŽs.  
On dŽfinit  une fonction rationnelle sur $\O$ par: 
 $$C^0(s, P', P,  \si ) : = B(s.P, P', s\si)^{-1} \otimes  A(P', s.P, s\si) \lambda (s),$$
 On introduit l'assertion suivante, dite Assertion A:\\
{\em Si $(P=MU, \O)\in \Theta$ et $Q$ est un sous-groupe parabolique de sous-groupe de LŽvi $M$, on a l'ŽgalitŽ de fonctions rationnelles sur $\O$:
$$B(Q, P, \si)^*= B(P,Q, \si), $$
o $^*$ dŽsigne l'adjoint.}\\
Cette assertion implique facilement que $C^0(s, P', P,  \si )$  est unitaire et s'Žtend en une fonction holomorphe au voisinage de $\O$.  La preuve de l'Assertion A nŽcessite quelques dŽtours comme son analogue dans [D4].\\ Mentionnons que la preuve de l'analogue de cette assertion pour les espaces symŽtriques rŽductifs rŽels  a ŽtŽ l'une des clefs de la preuve de formule de Plancherel pour ces espaces (cf. e.g., [D2]).\\
 On note ${\cal S}^{inv}$ l'ensemble des ŽlŽments  $(\phi_{P, \O})_{(P, \O)\in \Theta} $ de ${\cal S}$ tels que  que pour tout $(P, \O), P', s, \si $ comme ci dessus:
 $$\phi_{P', s\O}(s\si)= C^0(s, P', P,  \si )\phi_{P, \O}(\si).$$
 On dŽfinit une application linŽaire ${\cal W}$ de ${\cal S}$ dans $\CCU$ qui associe, pour $(P, \O)\in \Theta $,  ˆ $\phi \in \fhi$ le paquet d'ondes $f_\phi$ mulitipliŽ par une constante $c(P)^{-1}$.
 \\ Le rŽsultat principal s'Žnonce alors:
\\{\em  La transformation de Fourier-Whittaker, ${\cal F}$, est une bijection entre $\CCU$ et ${\cal S}^{inv}$. L'inverse de cette bijection est la restriction de  ${\cal W}$  ˆ ${\cal S}^{inv}$. Pour une topologie naturelle sur ${\cal S}$, ${\cal F}$ et ${\cal W}$ sont continues.
 \\En particulier: 
$$f= \sum_{(P=MU,\O)\in \Theta} c(P)^{-1} \int_\O \mu(\si)E^G_P ( ({\cal F}f)(P, \si) )d\si  , f \in \CCU .$$}
 On Žtudie aussi la transformation du produit scalaire $L^2$ sur $\CCU$ par ${\cal F}$. \\ 
Donnons une idŽe de notre  preuve.\\ L'injectivitŽ de ${\cal F}$ rŽsulte d'un rŽsultat de Joseph Bernstein (cf. [B1]). \\
  Puis on commence par introduire constant faible pour les fonctions dites  tempŽrŽes de $\cu$.
  Le terme constant faible se dŽduit naturellement du terme constant introduit dans [D4].
  Les intŽgrales de Jacquet tempŽrŽes et on calcule leur terme constant faible en procŽdant comme dans [D4].
 On  a besoin d'estimŽes qui se rŽduisent souvent ˆ des rŽsultats de [W] gr\^ace au lien entre les coefficients gŽnŽralisŽs et les coefficients lisses mentionŽ ci-dessus. 
On poursuit  comme dans [W], section VI, en introduisant la transformŽe unipotente de $f \in \CCU$ relative ˆ $P=MU$, sous-groupe parabolique anti-standard de $G$, $f^P$ dŽfini par:
$$f^P(m):= \delta_P^{1/2}(m) \int _Uf(mu)du, m\in M,$$
l'intŽgrale Žtant convergente car $f\in \CCU$.
Nous n'avons pas ŽtŽ capable de montrer directement que $f^P\in \CCMU$. Cela rŽsulte nŽammoins du rŽsultat principal. 
\\ Soit $(P, \O)\in \Theta$. On introduit la notion d'ŽlŽment $\phi \in \fhi$ trs rŽgulire. Nous ne donnerons pas ici la dŽfinition prŽcise. Disons seulement  que si $\phi$ est trs rŽgulire, elle est polynomiale et, notamment,  les fonctions sur $\O$, $C^0(s,P',P, \si) \phi(\si)$ sont holomorphes au voisinage de $\O$ dans $\O_\C$. Si l'Assertion A est vraie, tout $\phi\in \fhi$ polynomiale est trs rŽgulire. 
\\ Pour $ \phi$ trs regulire (ou pour $\phi\in \fhi$ si l'assertion A est vraie),  $f^P_\phi$ est un ŽlŽment de $\CCMU$ que l'on calcule. La preuve est inspirŽe d'un rŽsultat analogue de [W]. Pour $\phi$ trs rŽgulire, la preuve utilise des dŽplacements de contour d'intŽgrales qui rendent le calcul lŽgŽrement plus simple, ˆ nos yeux, que dans [W]. Cela reste une partie difficile de ce travail. Ensuite on raisonne par continuitŽ et densitŽ.
On en dŽduit, sous les m\^emes hypothses, le calcul de la transformation de Fourier-Whittaker de $f_\phi$. 
On en dŽduit aussi une formule pour le produit scalaire $L^2$ de deux paquets d'ondes $f_\phi$ et $f_{\phi'}$ :$$(f_\phi, f_{\phi'})_G:= \int_{U_0\backslash G}f_{\phi}(g)\overline{ f_{\phi'}(g)}dg.$$
En Žchangeant le r™le de $\phi$ et $\phi'$, on obtient une autre expression de ce produit scalaire. En faisant varier $\phi$ et $\phi'$ parmi les fonctions trs rŽgulires, dans l'ŽgalitŽ de ces deux  expressions de $(f_\phi, f_{\phi'})_G$, on en dŽduit l'assertion A . Les consŽquences  de cette Assertion sur les fonctions $C^0$ sont alors acquises.\\
On montre alors facilement que l'image de ${\cal F}$ est dans ${\cal S}^{inv}$.
Il est alors aisŽ  de montrer que  la restriction de ${\cal FW}$ ˆ ${\cal S}^{inv} $ est Žgale ˆ l'identitŽ, gr\^ace au calcul mentionnŽ ci-dessus du calcul de la transformŽe de Fourier-Whittaker des paquets d'ondes. Joint ˆ l'injectivitŽ de ${\cal F}$, ceci achve essentiellement la preuve du rŽsultat principal.

\section{Notations, Rappels}
\setcounter{equation}{0}
 \subsection{} 

 On reprend essentiellement les notations de [D4].
Soit $F$ un corps local non archimŽdien. 
On considre divers
groupes algŽbriques dŽfinis sur
$F$ et on utilisera des abus de terminologie du type suivant: ``soit $A$ un tore
dŽployŽ" signifiera ``soit
$A$ le groupe des points sur $F$ d'un tore dŽfini et dŽployŽ
sur
$F$". Avec ces conventions, soit $G$ un groupe  algŽbrique linŽaire rŽductif et
connexe.
On fixe un tore dŽployŽ maximal, $A_0$, de $G$ et on note $M_0$ son centralisateur dans $G$.  On fixe $P_0$ un sous-groupe parabolique minimal de $G$ qui admet $M_0$ comme sous-groupe de LŽvi.  On notera $U_0$ le radical unipotent de $P_0$.\ste
Si $P$ est un sous-groupe parabolique de $G$, on dit que $P$ est semi-standard (resp. standard )  si $M_0\subset P$ (resp. $P_0\subset P$). Si $P$ est semi-standard, il  possde un unique sous-groupe de LŽvi, $M$,  contenant $M_0$.  On dit que $M$ est un sous-groupe de LŽvi semi-standard.  
 \ste L'expression  ``$P=MU$ est sous-groupe parabolique semi-standard de $G$'' signifiera que $P$ est un tel sous-groupe, que $M$ est son sous-groupe de LŽvi semi-standard et que $U$ est son radical unipotent.  On notera  $P^-=MU^-$ le sous-groupe parabolique opposŽ ˆ $P$ de sous-groupe de LŽvi $M$. Un sous-groupe 
parabolique semi-standard de $G$, $P$,  sera dit anti-standard si $P^-$ est
standard. \ste 
 Si $H$ est un groupe algŽbrique, on note $Rat (H)$ 
le groupe des  caractres algŽbriques de $H$ dŽfinis sur $F$.\\ On fait le choix d'une unifomisante de $F$. Ce choix 
permet de relever le quotient de tout tore dŽployŽ $A$ par son sous-groupe compact maximal,$A^0$,  dans un rŽseau 
 de $A$, $\Lambda (A)$ : cÕest l'image
du rŽseau des sous-groupes ˆ  un paramtre de $A$ par l'Žvaluation en l'uniformisante. \\ Si $V$
est un espace vectoriel, on note $V'$ son dual et,  s'il est rŽel,  on note  $V_{\C} $
son complexifiŽ. \\ On note $A_{G}$ le plus grand tore dŽployŽ dans le centre
de $G$.  On note
$\a_{G}= Hom_{\Z} (Rat (G),\R)$. La restriction des caractres algŽbriques 
de $G$ ˆ $A_{G}$  induit un isomorphisme: 
\beq Rat(G)\otimes_{\Z}  \C \simeq Rat(A_{G})\otimes_{\Z} \C. \eeq  
On dispose de l'application canonique: \beq H_{G}: G\rightarrow
\a_{G}, \eeq
 dŽfinie par: 
\beq  \label{HG} e^{ \langle H_{G}(x), \chi\rangle   }= \vert \chi (x)\vert_{F}, \> x\in G, \chi \in
Rat (G ) .\eeq o
$\vert. \vert_F$ est la valuation normalisŽe de $F$. Le noyau de 
$H_{G}$, qui  est notŽ   $G^{1}$,  est l'intersection des noyaux des
caractres de $G$ de la forme $\vert \chi \vert_{F}$, $\chi \in Rat (G)$.
On notera
$X(G)_\C= Hom (G/G^{1}, \C^{*})$. C'est le groupe des caractres non ramifiŽs de $G$. \ste   On a des notations similaires pour des
sous-groupes de LŽvi semi-standard. Si $P$ est un sous-groupe
parabolique semi-standard de $G$,  on notera
$\a_{P}= \a_{M_{P}}$, $H_{P}=H_{M_{P}}$. On note $\a_{0}= \a_{M_{0}}$,
$H_{0}=H_{M_{0}}$.\\
   On note $\a_{G, F}$, resp. ${\tilde
\a}_{G, F}$,  l'image de $G$, resp. $A_{G}$, par $H_{G}$. Alors
$G/G^{1}$ est un rŽseau isomorphe ˆ
$\a_{G,F}$.
Soit $M$ un
sous-groupe de LŽvi semi-standard.  Les inclusions $A_{G}\subset
A_{M}\subset M\subset G$ dŽterminent   un morphisme  surjectif
 $\a_{M, F}\rightarrow  \a_{G, F}$, resp.  un morphisme injectif
 ${\tilde
\a}_{G, F}   \rightarrow {\tilde
\a}_{M, F}$, qui se prolonge de manire unique  en une  application linŽaire
surjective entre $\a_{M}$ et $\a_{G}$, resp. injective entre $\a_{G}$ et
$\a_{M}$.
La deuxime  application permet d'identifier $\a_{G}$ ˆ un sous-espace de
$\a_{M}$ et le noyau de la premire, $\a^{G}_{M}$, vŽrifie ;
\beq \a_{M}= \a^{G}_{M}\oplus \a_{G}.\eeq
Soit $P=MU$ un \sgpss. On note $\Sigma(A_M)$ (resp. $\Sigma(P)$)
l'ensemble des racines de
$A_M$ dans l'algbre de Lie de $G$ (resp. $P$) qui s'identifie ˆ un
sous-ensemble de $\a_M '$. On note
$\Delta(P)$ l'ensemble des racines simples de $\Sigma(P)$. Si $\aa\in
\Sigma(A_M)$, on note $U_\aa$ le sous-groupe radiciel de  $U$ correspondant ˆ $\aa$.   On peut associer ˆ tout $\aa \in \Sigma(A_M) $ une coracine $\check{\aa}\in \a_M$ (cf. [A], section 3).\ste  
On note ${}^-{\a_P^G}'$ ( resp. ${{}^-{\overline \a}_P^G}'$) l'ensemble des $\lambda \in \a_M'$ de la forme 
$$\lambda= \sum _{\aaÊ\in \Delta(P)}x_\aa \aa$$
avec $x_\aa< 0$ ( resp. $x_\aa\leq 0$). 
\ber \label{lchi}   Soit  $\l \mapsto \chi_\l $ l'application $(\a_{G}')_{\C}\rightarrow X(G)_\C\rightarrow 1$  qui, en utilisant la dŽfinition   de $\a_G$,
associe ˆ
$\chi\otimes s$,  le caractre $g\mapsto \vert \chi (g)\vert_F
^{s}$. \eer  Le noyau est un rŽseau et cela dŽfinit sur $X(G)_\C$ une structure de
variŽtŽ algŽbrique complexe pour laquelle $X(G)_\C\simeq \C^{*d}$, o
$d=dim_{\R}\a_{G}$. Pour $\chi \in X(G)_\C$, soit $\lambda \in
\a_{G,\C}'$ un ŽlŽment se projetant sur $\chi $ par l'application
(\ref{lchi}). La partie rŽelle $Re\>\lambda \in \a_{G}'$ est indŽpendante du
choix de $\lambda$. Nous la noterons $Re \> \chi$. Si $\chi \in Hom (G,
\C^{*})$, le caractre $\vert \chi \vert$ appartient ˆ $X(G)$. On pose
$Re\> \chi=  Re\> \vert \chi\vert $. De mme, si $\chi \in Hom(A_{G},
\C^{*})$, le caractre $\vert \chi \vert$ se prolonge de faon unique en un 
ŽlŽment de $X(G)$  ˆ valeurs dans $\R^{*+}$, que l'on note
encore  
$\vert
\chi \vert$  et on pose $Re \> \chi=  Re\> \vert \chi
\vert$.\ste Soit 
$
\> X(G):=
\{
\chi
\in X(G)\vert Re
\>  \chi =0\}$ l'ensemble des ŽlŽments unitaires de $X(G)$. \\
Les notations ainsi dŽfinies s'appliquent ˆ tous les sous-groupes de LŽvi semi-standard de $G$.
\\ On choisit $K$ un sous-groupe compact maximal de $G$, dont on suppose qu'il est le fixateur d'un point spŽcial de l'appartement associŽ ˆ $A_0$ dans l'immeuble de $G$. Pour le  rŽsultat suivant, cf.  [C], Prop. 1.4.4:
\ber \label{iwa} Il existe une suite dŽcroissante de sous-groupes compacts ouverts de $G$, $K_n$, $n\in \N$,  telle que pour tout  $n\in \N^*$, $H= H_n$
est normal  dans  $K= H_0$ et pour tout sous-groupe parabolique standard de $G$, $P$,  on a: \ste 
1) $ H= H_{U^-} H_M H_U$ o $H_{U^-} = H \cap U^-$, $H_M=
H\cap M$,  $H_U= H \cap U. $\ste
   2) Pour tout $a\in A_M^-:=\{a \in A_M\vert \vert \alpha(a) \vert_F \leq 1, \aa \in \Sigma(P) \}$, $aH_Ua^{-1}\subset H_U$, $a^{-1}H_{U^-}a \subset 
K_{U^-}.$ \ste
3) Le groupe $H_M$ vŽrifie 1) et 2) relativement aux sous-groupes paraboliques de $M$ contenant $P_0\cap M$.
 \ste 
4) La suite  $H_n$ forme une base de voisinages de l'identitŽ de  l'ŽlŽment neutre de  $G$. \eer 
Soit $H$ un sous-groupe compact ouvert de $G$. 
 On dit que  $H$  possde une factorisation  d' Iwahori par rapport ˆ  $(P,P^-)$ si 1) et 2) sont satisfaits.
  \subsection{Choix des mesures}
 On munit $G$ (resp. $K$) d'une mesure de Haar, $dg$, (resp.  $dk$) telle que le volume de $K$ soit Žgal ˆ 1. Si $H$ est un sous-groupe compact ouvert de $G$, on note $vol(H)$ son volume pour la mesure de Haar sur $G$. On note $e_H$ le multiple de la fonction indicatrice de $H$ dont l'intŽgrale sur $G$ est Žgale ˆ 1.\\  Pour tout sous-groupe fermŽ, $H$,  de $G$, on note $dh$ une mesure de Haar ˆ gauche sur $H$, dont le choix sera Žventuellement spŽcifiŽ et 
on note
$\delta_H$ la fonction module de
$H$. \ste 
Pour tout espace totalement discontinu $Z$, on note $C_c^{\infty}(Z)$ (resp.  $C(Z)$, resp. $C_c(Z)$)
l'espace des fonctions localement constantes, ˆ support compact  (resp. continues, resp. continues ˆ support compact) sur  $Z$ ˆ valeurs dans  
 $\C$. \ste
Soit $P=MU$ un sous-groupe parabolique semi-standard
de $G$. Alors $G=PK$ et le groupe $M\cap K$ vŽrifie relativement ˆ $M$ les mmes propriŽtŽs que $K$ relativement ˆ $G$. Pour $g\in G$, on choisit $u_P(g)\in U$, $m_P(g) \in
M$, $k_P(g) \in K$,  de sorte que $g=u_P(g)m_P(g)k_P(g)$. 
On note $du^-$ la mesure de Haar sur $U^-$ telle que:
\beq \label{fintu} \int_{U^-} \delta_P(m_P(u^-)) du^-=1\eeq
et de mme pour $U$.
Il esiste une
unique mesure de Haar sur
$M$,
$dm$, telle que pour tout $f
\in C_c(G)$ on a (cf. [W], I (1.2)): 
 \ber \label{iwas} $$\int_{G}f(g)dg= \int_{U\times M\times K} f(umk) \delta_P(m)^{-1} dkdmdu. $$
  $$\int_{G}f(g)dg= \int_{U\times M\times U^-} f(umu^-)
\delta_P(m)^{-1}dudmdu^-.$$\eer
  Soit $f$ une fonction  continue sur $K$ et invariante ˆ gauche par  $K\cap P$. Alors (cf. par exemple [K], ch. V, section 6, consŽquence 7, pour la version rŽelle): 
\ber \label{intu-}$$ \int_Kf(k)  dk= \int_{U^-}
f(k_P(u^-))\delta_P(m_P(u^-) ) du^-, $$
\eer
o l'intŽgrale du membre de droite est absolument convergente.
En particulier si $f$ est une fonction continue sur $G$   telle que
$f(umg)= 
\delta_P(m)  f(g)$ pour $u\in U, m\in M, g\in G $, on a:
\ber \label{intu--}$$ \int_Kf(k)  dk=  \int_{U^-}
f(u^-)du^-. $$ \eer
o l'intŽgrale du membre de droite est absolument convergente. \\
On fixe  la  mesure de Haar de masse totale 1 sur  $A_G\cap K$ et sur $X(A_G)$,  qui s'identifie au dual unitaire de   $A_G/A_G \cap K$. On fixe sur $A_G/A_G \cap K$ la mesure de Haar duale de celle sur  $X(A_G)$. Des mesures ainsi choisies,  on dŽduit une mesure sur $A_G$, notŽe $da_G$. L'homomorphisme de restriction dŽtermine un morphisme surjectif de $X(G)$ sur  $X(A_G)$. On fixe sur $X(G)$ une mesure de Haar telle que ce morphisme prŽserve localement les mesures de Haar choisies.\\
On choisit un ensemble de reprŽsentants dans $K$, $W^G$,  du quotient du  normalisateur de $M_0$
dans $G$ par son centralisateur, $\overline{W}^G$, qui existe parce que $K$ est le fixateur d'un point spŽcial de l'appartement associŽ ˆ $A_0$. On le choisit de telle sorte qu'il contienne l'ŽlŽment neutre de $G$,  qu'on notera $1_G$ ou seulement $1$ s'il n'y a pas de confusion. On munit $\a_0$ d'une norme euclidienne invariante par l'action naturelle de $\overline{W}^G$ sur $\a_0$. 
\ber \label{xpoint}Si $x\in G$ et $Y$ est une partie  de $G$,  on notera $x.Y:=\{xyx^{-1}\vert y \in Y\}$. 
De mme si $H$ est un sous-groupe  de $G$, et $(\si, E)$ est une reprŽsentation de $H$, on notera $x\si$ la reprŽsentation de $x.H$ dans $xE:=E$ dŽfinie par $x\si(xh x^{-1})= \si(h), h\in H$. \ste
On rappelle que le dual d'un espace vectoriel $V$ est notŽ $V'$. Si $T$ est une application linŽaire entre deux espaces vectoriels complexes, on note $T^t$ sa transposŽe. Si $(\si,E) $ est une reprŽsentation de $H$, on note $\si'$ la reprŽsentation de $H$ dans $E'$ dŽfinie par $\si' (h)= \si(h^{-1})^t, \>h \in H$. \ste 
On notera aussi, pour toute application dŽfinie sur $H$, $f$, on dŽfinit: $$(\l(h)f)(h')= f(h^{-1} h'), \> (\rho(h)f)(h')= f(h'h), \> h, h'\in H.$$\eer
 \subsection{ReprŽsentations induites} 
 Les reprŽsentations lisses de $G$ et de ses sous-groupes fermŽs seront toujours ˆ coefficients complexes. \\ Si $(\pi,V)$ est une reprŽsentation lisse de $G$, on note $\pi'$ (resp. $\check{\pi}$) la reprŽsentation contragrŽdiente  de $G$ dans le dual algŽbrique, $V'$,  ( resp. lisse, $\check{V}$) de $V$.\\
  Soit $P=MU$ un sous-groupe parabolique semi-standard de $G$. 
  Si $(\si, E)$  est une reprŽsentation lisse de $M$, on l'Žtend en une reprŽsentation de $P$ triviale sur $U$. Soit $\chi \in X(M)_\C$. On note $E_\chi$ l'espace de la reprŽsentation $\si\otimes \chi$, qu'on notera  $\si_\chi$, et on note
$i^G_{P}E_\chi$ l'espace des fonctions $v$ de $G$ dans $E$, invariantes ˆ droite par au moins
un sous-groupe compact ouvert de $G$ et telles que $v(mug)= \delta_P(m)^{1/2}\si_\chi (m) f(g)$ pour tout
$m\in M$, $u\in U$, $g\in G$. On note $i^G_{P}\si_\chi$ la
reprŽsentation de $G$  dand 
$i^G_{P}E_\chi$ par translations ˆ droite.  On note $i^K_{P\cap K}E$ l'espace des fonctions
$v$ de $K$ dans $E$, invariantes ˆ droite par un sous-groupe compact ouvert de $K$ et telles
que $f(pk)= \si(p)f(k)$ pour tout $k\in K$ et $p\in P\cap K$. La restriction des fonctions ˆ
$K$ dŽtermine un isomorphisme de $i^G_{P}E_\chi$ sur $i^K_{P\cap K}E$. On notera ${\overline
i}^G_{P}\si _ \chi$ la reprŽsentation de $G$ dans $i^K_{P\cap K}E$ dŽduite de
$i^G_{P}\si _\chi$ par transport de structure. Cette reprŽsentation sera appelŽe  la rŽalisation compacte de 
$ i^G_P\si_\chi$ dans cet espace indŽpendant de
$\chi$. 
Si $v \in i^K_{P\cap K}E$, on note
$v_\chi$ l'ŽlŽment de  $i^G_{P}E_\chi$ dont la restriction ˆ $K$ est Žgale ˆ $v$.   
\ste Si $\si$ est unitaire,  on munit $i^K_{P\cap K}E$ du produit scalaire dŽfini par: 
  \beq (v, v')= \int_{K}(v(k),v'(k))dk, \> v, v' \in i^K_{P\cap K}E .\eeq
 Alors, muni de ce produit scalaire, la reprŽsentation $\overline{ i}^G_P \si_\chi$ est unitaire pour $\chi$ unitaire et par consŽquent, par transport de structure, $i^G_P\si _\chi$ Žgalement.  On a alors, gr\^ace ˆ  (\ref{intu--}):
\beq \label{prodscalu-}(v_\chi, v'_\chi) =\int_{U^-} (v_\chi(u^{-}), v'_\chi (u^-)) du^-. \eeq
\section{Fonctionnelles de Whittaker tempŽrŽes, terme constant faible }
\setcounter{equation}{0}
\subsection{Fonctions et fonctionnelles de Whittaker}
Soit $U_0$ le radical unipotent de $P_0$.  Soit $\psi$  un caractre unitaire lisse de $U_0$ non dŽgŽnŽrŽ, i.e. tel que pour toute    racine $P_0$-simple de $A_0$, $\aa$,  sa restriction au  sous-groupe radiciel $({U_0})_\alpha$ soit non triviale.  
 \\On note $\cu$ l'espace des fonctions de Whittaker lisses sur $G$, i.e. des fonctions, $f$, sur $G$,  invariantes ˆ droite par un sous-groupe compact ouvert de $G$ et telles que:
  $$f(u_0g)= \psi(u_0) f(g), g\in G, u_0 \in U_0.$$ 
  On note $\ccu$ l'espace des ŽlŽments de $\cu$ qui sont ˆ support compact modulo $U_0$.
  \\ 
 Si $(\pi, V)$ est une reprŽsentation lisse de $G$, on note $Wh(\pi, U_0)$ ou $Wh(\pi)$ l'espace des formes linŽaires, $\xi$,  sur $V$ telles que:
$$ \langle \xi, \pi(u_0)v \rangle  = \psi(u_0)  \langle \xi, v\rangle  , v\in V, u_0\in U_0.$$
Si $\pi$ est de longueur finie,  $Wh(\pi) $ est de dimension finie (cf.[BuHen], ThŽorme 4.2) et [D3] pour une autre dŽmonstration).   Notez que si le groupe n'est pas quasi-dŽployŽ, cette dimension n'est pas toujours infŽrieure ou Žgale ˆ  1, comme expliquŽ dans l'introduction de [BuHen]. 
\\ Si $v\in V$, on note $c_{\xi, v}$ le coefficient gŽnŽralisŽ dŽfini par:
\beq \label{cxivl} c_{\xi, v}(g):=  \langle \xi, \pi(g) v\rangle  , \>g \in G.\eeq
C'est un ŽlŽment de $\cu$. \\
\subsection{Rappel sur le terme constant des fonctions de Whittaker}
On note, pour $\varepsilon >0$,  $$A_0^-(\varepsilon):=\{ a\in A_0\vert \vert \alpha(a)\vert_F \leq \varepsilon, \aa \in \DD(P_0)\} \> \> et \>\> A_0^-= A_0^- (1).$$
$$M_0^-(\varepsilon):=\{ m\in M_0\vert \langle  \alpha, H_0 (m) \rangle  \leq log\> \varepsilon, \aa \in \DD(P_0)\} \> \> et \>\> M_0^-= M_0^- (1).$$ Comme dans [D3], Lemme 3.1,  on voit que:  \ber \label{H'} Pour tout sous-groupe compact ouvert $H$ de $G$, il existe un
sous-groupe compact ouvert de $G$,  $H'$, contenu dans $H$, tel que pour toute reprŽsentation lisse de $G$,
$(\pi,V)$, tout $\xi
\in Wh(\pi)$ et tout
$v\in V^H$ on ait:
$$ \langle \xi, \pi(a)v\rangle  =  \langle  e_{H'}\xi, \pi(a) v\rangle   , a \in A^-_0$$
o $e_{H'}\xi$ est l'ŽlŽment de ${\check V}$ dŽfinit par $ \langle e_{H'}\xi, v\rangle  : =
 \langle \xi, \pi(e_{H'}) v\rangle  , \> v\in V$, ce qui se rŽŽcrit:
 $$c_{\xi,v}(a)= c_{e_{H '}\xi, v}(a), a \in A_0^-. $$\eer 
 On rappelle (cf. [D3],  Lemme 5.4):
\ber \label{fa-} Soit $H$ un sous-groupe compact ouvert de $G$. Il existe $C>0$ 
tel que, pour tout $f\in \cu$, invariante ˆ droite par $H$,  $f(a)=0$ si $\vert a^\aa \vert_{F} >C$ pour au moins un ŽlŽment $\aa$ de  $\Delta(P_0)$,  i.e. pour  $a$ ŽlŽment du complŽmentaire de $A_0^-(C).$
\eer 
Ceci implique \ber\label{trans} Soit $H$ un sous-groupe compact ouvert de $G$. Il existe $a_0 \in A_0$ tel que pour tout $f \in \cu$ invariante ˆ droite par $H$, la restriction de $\rho(a_0)f$ ˆ $A_0$ est ˆ support dans $A_0^-$. \eer 
 Soit $I$ un sous-ensemble fini de $K$  tel que $K \subset IH$ et soit $F_0$ un sous-ensemble fini de $M_0$ tel que $G= U_0A_0 F_0 K$. Appliquant (\ref{fa-}) aux translatŽes ˆ droite de $f$ par les ŽlŽments de $F_0I$,  on en dŽduit:
 \ber \label{suppf} Il existe $C>0$ tel que  pour tout  $f \in \cu$ invariante ˆ droite par $H$, $f $ est ˆ support dans $U_0 M_0^-(C) K$ \eer
Rappelons la caractŽrisation du terme constant des ŽlŽments de $\cu$ et
des fonctionnelles de Whittaker.
Soit $P=MU$ un sous-groupe parabolique standard de $G$.
Avec les notations ci-dessus, on note $(\pi_P, V_P)$ le produit tensoriel entre, d'une part 
 la reprŽsentation de $M$ dans   le  quotient de  $V$ par le 
$M$-sous-module engendrŽ  par les  $\pi(u)v-v$, $u\in U, v\in V$, et d'autre part 
la reprŽsentation  de  $M$ sur  $\C$ donnŽe par  $\delta_P^{-1/2}$.   On 
note, pour $v \in V$, $j_P(v) $ ou $v_P$ sa projection naturelle dans  
$V_P$.
\ste Soit $\Theta_P$ l'ensemble des ŽlŽments de  $\DD(P_0)$  qui  sont racines de $A_0$ dans l'algbre de Lie de $M$. 
On note, pour $\varepsilon >0$:
$$A_0^-(P, < \varepsilon):=\{ a \in A_0^-\vert \vert \alpha(a)\vert_F < \varepsilon, \aa \in \DD(P_0)\setminus \Theta_P\}.$$
D'aprs [D3] ThŽorme 3.4, Remarque  3.5 et  Proposition 3.6, on dispose d'une unique   application linŽaire  $Wh(\pi)\mapsto Wh(\pi_P)$, $\xi \to
j_{P^-}(\xi)$, notŽe aussi $\xi \mapsto \xi_P$ pour plus de commoditŽ,  qui vŽrifie:  
\ber  \label{jpxi}Pour tout sous-groupe compact ouvert $H$ de $G$, il existe
$\varepsilon_H>0$, indŽpendant de $P$,  avec les propriŽtŽs suivantes:
\ste Pour toute reprŽsentation  lisse
$(\pi, V)$ et 
$\xi
\in Wh(\pi)$, on a:   
$${\delta_P^{1/2}}(a)\langle \xi_P, \pi_P(a) v_P\rangle  _P= \langle \xi, \pi(a)
v \rangle, \> a \in A_0^{-}(P,<\varepsilon_H), \> v\in V^H.$$\eer 
Le terme constant  le long de $P$ d'un ŽlŽment $f$ de $\cu$ a
ŽtŽ dŽfini dans [D3],  DŽfinition 3. C'est un ŽlŽment $f_P$  de $\cmu$ tel que:
 \ber \label{ffp} $$\delta_P^{1/2}(a)  f_P(a)= f(a), \> a \in
A^-_0(P, <\varepsilon_H). $$\eer
\ber \label{covct} L'application $f\mapsto f_P$ est un morphisme de $P$-modules entre $\cu$ et $\cmu$,
 o $P$ agit par reprŽsentation rŽgulire
droite sur le premier espace et $M$ (resp. $U$) agit par reprŽsentation rŽgulire
droite tensorisŽe par $\delta_P^{1/2}$(resp. trivialement)  sur le second. 
\eer \label{ctcxiv}
 Soit $ (\pi, V)$ une reprŽsentation lisse de $G$ et $\xi\in Wh(\pi)$. On  a (cf. [D3], Proposition 3.13) :
\beq \label{tctcoef} (c_{\xi,v})_P= c_{ \xi_P, v_P}, \> v\in V. \eeq 
 \subsection{Fonctions et fonctionnelles de Whittaker tempŽrŽes ou  de carrŽ intŽgrable }
 Suivant [W], section I.1, on  fixe un plongement algŽbrique $\tau$ de $G$ dans $GL_n(F)$ tel que $\tau(K)\subset GL_n({\cal O})$ o $\cal O$ est l'aneau des entiers de $F$. Si $(a_{i,j})$ (resp. $b_{i,j}$ ) sont les  coefficients de la matrice $\tau(g)$ (resp. $ \tau(g)^{-1}$) on note $\Vert  g \Vert = sup_{i,j} sup ( \vert a_{i,j}\vert_F , \vert b_{i,j}\vert_F )$. 
 On pose \beq \label{sii} \sii (g)= log \Vert  g \Vert ,\eeq
 o l'on a soulignŽ $\si$ pour Žviter la confusion avec les reprŽsentations. \\
 On introduit la fonction $\Xi$ d'Harish-Chandra (cf. [W] section II.1)   et  l'espace de Schwartz-Harish-Chandra $ \CC(G) $ (cf. [W], section III.6).\\
 Pour $g \in G$, on note $\delta_0(g)$ ( resp. $m_0(g)$, etc..) au lieu de $\delta_{P_0}(g)$ (resp. $m_{P_0}(g)$, etc..).\\
 On note $\CCU$ l'espace des ŽlŽments $f$  de $\cu$  tels que pour tout $d\in \N$:
 \beq \label{nud} \nu_d(f) := sup_{g \in G} \delta^{1/2} _{0} (m_{0}(g))^{-1} (1+ \vert H_{0} (m_{0} (g))\vert )^d \vert f(g)\vert<\infty .Ê\eeq 
  Pour tout sous-groupe compact, $H$, de $G$, on  note $\CCU^H$, l'espace des $f\in \CCU$ invariantes ˆ droite par $H$, que l'on munit de la topologie dŽfinie par les $\nu_d$. Montrons:
 \ber \label{dense} $C_c^{\infty} (U_0 \backslash G, \psi)^H$ est dense dans $\CCU^H$ .\eer En effet soit $f \in \CCU$. 
D'aps (\ref{suppf}), il existe $C>0$ tel que $f$  a son support contenu dans $U_0 M_0^-(C) K$. Notons, pour $n \in \N^* $,  $ u_n $ l'indicatrice  de l'ensemble compact $U_n:= \{ m\in M_0\vert  -log n  \leq \langle  \alpha, H_0 (m) \rangle  \leq log\> C, \aa \in \DD(P_0)\}$.  Il existe une constante $C'$ tellle pour tout $n \in \N^*$: $$ \Vert H_0(m)\Vert > C' log \>n, m\in M_0^-(C) \setminus U_n. $$ On en dŽduit :$$ \nu_d(f-u_n f)\leq (1+ C' log \>n) \nu_{d+1} (f).$$ Donc $(u_nf) $ est une suite d'ŽlŽments de $\ccu$ qui converge vers $f$ dans $\CCU^H$.\\ On munit $ \CCU$ de la topologie limite inductive des $\CCU^H$.  
 \ber \label{deftemp} On dit que $f\in \cu$ est tempŽrŽe si il existe $d\in \N$ telle que: 
  $$sup_{g \in G} \delta_0(m_{0}(g))^{-1} (1+ \Vert H_0(m_{0} (g))\Vert )^{-d}  \vert f(g)\vert < \infty.$$ \eer
On note $C^w( U_0\backslash G, \psi) $ l'espace vectoriel des fonctions tempŽrŽes. Muni de l'action rŽgulire droite c'est un $G$-module lisse.  \\
D'aprs la premire formule de (\ref{iwas}), il existe un entier $d \in \N$ tel que \beq \label{intcv} \int_{U_0\backslash G} \delta_0(m_0(g)) (1+ \Vert H_0(m_0(g) ) \Vert )^{-d}  dg <\infty .\eeq 
 Donc: \ber  \label{intff} Si $f_0 \in C^w( U_0\backslash G, \psi) $ et $f \in \CCU$, $ f \overline{ f_0} $ est ŽlŽment de $L^1(U_0\backslash G)$ et l'application $f\mapsto \int _{U_0\backslash G} f (g)\overline{ f_0(g)}dg $ est une forme linŽaire continue sur $\CCU$.  \eer 
Soit $(\pi, V) $ une reprŽsentation admissible de $G$ et $\xi \in Wh(\pi)$. On dit que $\xi$ est tempŽrŽe si l'espace $\A(\xi)$ formŽ des 
$c_{\xi,v} , v \in V$ est contenu dans $ C^w(U_0\backslash G, \psi) $. 
On note $ {\cal A}^w (U_0\backslash G, \psi)$ la rŽunion des $\A(\xi)$, lorsque $\pi$ dŽcrit les reprŽsentations lisses admissibles de $G$ et $\xi$ dŽcrit les fonctionnelles de Whittaker tempŽrŽes  de $ \pi$. 
Montrons  que: \ber Un ŽlŽment $f$ de  $C^w ( U_0\backslash G, \psi)$  appartient ˆ  ${\cal A}^w( U_0\backslash G, \psi)$ si et seulement si $f$ est $ZB(G)$-finie, o $ZB(G) $  le centre de Bernstein de $G$. \eer En effet, comme tout vecteur d'un 
module admissible est $ZB(G)$-fini, la partie seulement si est claire. Supposons maintenant que $f\in C^w(U_0\backslash G, \psi)$ est $ZB(G)$-finie.  Alors le $G$-module, $(\pi, V)$, engendrŽ par $f$ sous la reprŽsentation rŽgulire droite est un module de longueur finie (cf. [D4], Lemme 2) , donc admissible. La mesure de Dirac en l'ŽlŽment neutre fournit un ŽlŽment $\xi$ de $Wh(\pi)$ et $c_{\xi,f}=f$. Donc $f \in{\cal A} ^w (U_0\backslash G, \psi)$ Ceci achve de prouver la partie si de notre affirmation.   \\
 Soit $(\pi, V)$ une reprŽsentation admissible de $G$. Pour $\chi\in Hom(A_G, \C^*)$ on pose: $$V_\chi
 = \{v \in V\vert Il \> existe \> d\in \N\> tel \> que \> (\pi(a)-Ê\chi(a))^d v= 0, a \in A_G\}.$$
 Si $\xi \in Wh( \pi)$, on note $\xi_\chi$ la restriction de $\xi$ ˆ $V_\chi$. On appelle exposant de $\pi$ (resp. $\xi$)  un caractre $\chi$ tel que $V_\chi$ (resp. $\xi_\chi$) soit non nul. On note ${\cal E} xp (\pi)$ (resp. ${\cal E} xp (\xi)$) l'ensemble des exposants de $\pi$ (resp. $\xi$).\\
 Une  fonctions mesurable, $f$,   sur $G$ telle que $f(ug)= \chi(u)f(g)$ pour  $u \in U_0$ et  $g \in G$, et telle que : $$\Vert fÊ\Vert^2:=  \int_{U_0Ê\backslash G} \vert f(g)\vert ^2 dg$$ sera dite fonction de Whittaker de carrŽ intŽgrable. L'espace des classes  modulo l'Žquivalence presque partout de fonctions de Whittaker de carrŽ intŽgrable dŽfinit un espace de Hilbert, $\lu$, 
 sur lequel $G$ agit unitairement  et continument par la reprŽsentation rŽgulire droite $\rho$.\\
 Soit $(\pi, V)$ une reprŽsentation lisse admissible de $G$  admettant un caractre central unitaire et $\xi\in Wh(\pi)$. On dit que $\xi$ (resp. $\pi$) est de carrŽ intŽgrable si pour tout $v\in V$, $c_{\xi, v} $ est de carrŽ intŽgrable modulo $A_GU_0 \backslash G$ ( resp. pour tout $v\in V$ et $\check{v}$ ŽlŽment du dual lisse $\check{V}$ de $V$, le coefficient lisse $c_{ \check{v},v}$ est de carrŽ intŽgrable sur  $A_G \backslash G$). 
 \subsection{Critre pour les fonctionnelles de Whittaker tempŽrŽes ou de carrŽ intŽgrable}
 \begin{prop}\label{crit}
Soit $(\pi, V)$ une reprŽsentation lisse admissible de $G$ et $\xi\in Wh(\pi)$. Les conditions suivantes sont Žquivalentes:
\\(i) $\xi$ est de carrŽ intŽgrable (resp. tempŽrŽe) 
\\(ii) pour tout sous-groupe parabolique standard $P=MU$ de $G$ et tout $\chi \in {\cal E}xp (\xi_P)$, on a $\chi \in ^- {\a_P^G}'$   
(resp.  $\chi \in ^-{{\overline \a} _P^G}'$).  
\\ (iii) pour tout sous-groupe parabolique standard maximal  $P=MU$ de $G$ et tout $\chi \in {\cal E}xp (\xi_P)$, on a $\chi \in ^- {\a_P^G}'$ (resp.  $\chi \in ^- {{\overline \a} _P^G}'$)
\end{prop}
\dem 
Le cas de carrŽ intŽgrable a tŽtŽ traitŽ dans [D4], Proposition 13.
Pour la tempŽrance,  on procde comme dans la dŽmonstration du critre analoque pour les groupes (cf. [W], Proposition III.2.2). \qed
\begin{prop} \label{squaresquare}
Soit $(\pi, V)$ une reprŽsentation lisse admissible de carrŽ intŽgrable  (resp. tempŽrŽe)  de $G$ et $\xi\in Wh(\pi)$, alors $\xi$ est de carrŽ intŽgrable (resp. tempŽrŽe).
\end{prop}
\dem En effet pour tout 
 sous-groupe paraboliqsue standard, les exposants de $\xi_P$ sont des exposants de $\pi_P$. Alors le corollaire rŽsulte de la Proposition prŽcŽdente jointe au ThŽorme 4.4.6 de [C] (resp. ˆ la Proposition III.3.2 de [W]). \qed
Le rŽsultat suivant  suivant  a ŽtŽ conjecturŽ par Lapid et Mao et  prouvŽ dans [D4], Thorme 8. IndŽpendamment Matringe ([Ma]) en a donnŽ une preuve pour certains groupes.
\ber  \label{lapid} Soit $(\pi, V)$ une reprŽsentation lisse  irrŽductible de $G$ admettant un caractre central unitaire. \\ S'il existe $\xi \in Wh(\pi)$ de carrŽ intŽgrable    non nulle, alors $\pi$ est de carrŽ intŽgrable. \eer
\subsection{Terme constant faible}
Soit $(\pi, V)$ une reprŽsentation lisse  admissible de $G$ et $\xi\in Wh(\pi)$ tempŽrŽe. Soit $P= MU$ un sous-groupe parabolique standard de $G$. On dŽcompose le module Jacquet $(\pi_P, V_P)$ en une somme directe: 
$$(\pi_P, V_P)= (\pi^w_P, V^w_P) \oplus (\pi^-_P, V^-_P) \oplus (\pi^r_P, V^r_P)$$
o
$$V^w_P= \oplus_{\chi \in {\cal E}xp(\pi_P), Re \chi=0} V_{P, \chi},
 V^-_P= \oplus_{\chi \in {\cal E} xp(\pi_P), Re \chi \not= 0, Re\chi \in {} ^-{{\overline \a}^G_P}'}, V_{P, \chi}, V^r_P= \oplus_{\chi \in {\cal E}xp(\pi_P) Re\chi \notin {} ^-{{\overline \a}^G_P}'} V_{P, \chi}.$$
 La Proposition  \ref{crit} montre que  la restriction de $\xi_P$ ˆ $ V^r_P$ est nulle. \\
  On appelle terme constant faible de $\xi$ la restriction de $\xi_P$ ˆ $ V_P^w$. On le note $\xi^w_P$. Gr\^ace ˆ la dŽcomposition ci-dessus on pourra aussi regarder $\xi^w_P$ comme une forme linŽaire sur $V_P$. De l'hŽrŽditŽ du terme constant (cf. [D3] , Proposition 3.16), on dŽduit que $\xi^w_P $ est tempŽrŽe.  \\
 Soit $P= MU$ un sous-groupe parabolique standard de $G$. Pour toute fonction $f: A_M\to \C$, on Žcrit
 $$lim_{a\opto P -\infty }f(a)=0$$ 
 si et seulement si pour tous $ \ep, \eta>0$, il existe $R>0$ tel que pour tout $a \in A_M$ vŽrifiant les conditions:
\\(i) $\langle \alpha, H_M(a) \rangle  <-R $ pour tout $\alpha \in \Sigma(P)$.\\
(ii) $\langle  \alpha, H_M(a)\rangle  < \eta \langle \beta, H_M(a)\rangle $ pour tout $\alpha, \beta \in \Sigma(P)$,
on ait $\vert f(a) \vert < \ep$. 
\begin{prop}
(i) Soit $f\in \AA$. Il existe un unique ŽlŽment $f^w_P\in \AAM$ tel que, pour tout $m \in M$, on ait:
$$ lim_{a\opto P - \infty } \delta_P^{1/2} f(ma)- f^w_P(ma) =0.$$
On appelle $f_P^w$ le terme constant faible de $f$ le long de $P$. 
\\(ii) Soit $\pi$ une reprŽsentation lisse admissible de $V$, soit $\xi$ un ŽlŽment tempŽrŽ de $Wh(\pi)$ et $v\in V$. Soit $f= c_{\xi, v}$.  Alors $f^w_P $ est Žgal ˆ $c_{\xi_P^w, v_P}$. \
\end{prop}
\dem  
La dŽmonstration est analogue ˆ celle du Lemme III.5.1 de [W]. 
\qed
Pour $f \in {\cal A}^w(U_0 \backslash G, \psi)$, on dŽfinit une application $f^{w,ind}_P: G\to \AAM$ par $ f^{w,ind}_P(g)= (\rho(g)f)^w_P$. L'application $f \mapsto f^{w,ind}_P$ est un entrelacement entre la reprŽsentation rŽgulire droite de $G$ dans ${\cal A}^w(U_0 \backslash G, \psi)$ et la reprŽsentation induite $i^G_P \AAM$, d'aprs (\ref{covct}). 
\subsection{Terme constant faible pour les fonctions de Whittaker et les coefficients lisses de reprŽsentations tempŽrŽes}
   \begin{lem}\label{ctct}Soit $(\pi, V)$ une reprŽsentation de longueur finie de $G$, $P$ un sous-groupe parabolique standard de $G$ et $\xi\in Wh(\pi)$. Soit $H$ un sous-groupe compact ouvert  de $G$ et $H'$ comme dans (\ref{H'}). On introduit le terme constant (resp. terme constant faible)  des coefficients lisses de reprŽsentations lisses admissibles  sur $G$ (resp. et tempŽrŽes) comme dans [W], Proposition I.6.2. Si. ${\check v} \in {\check V}$, on note $c_{{\check v}, v}$ le coefficient lisse  dŽfini par $c_{\check{v}, v}(g)= \langle \check{v}, \pi(g)v\rangle$. 
  \\Soit $v\in V^H$. On a $v':=e_{H'}\xi\in {\check V}$ et:  
   $$({c_{\xi, v}})_P(a_0a)= (c_{v',v})_P(a_0a), a_0 \in A_0^-, a\in A_M.$$
 Si de plus $(\pi,V)$ est tempŽrŽe, on a: 
  $$({c_{\xi, v}})^w_P(a_0a)= (c_{v',v})^w_P(a_0a), a_0 \in A_0^-, a \in A_M.$$
   \end{lem} 
   \dem
Prouvons la premire ŽgalitŽ. Pour $a_0\in A_0^-$ fixŽ,  les deux membres sont des fonctions $A_M$-finies qui sont Žgales sur $A_M\cap A_0^-(P, \ep)$
pour $\ep>0$ assez petit,  d'aprs les propriŽtŽs des deux  termes constants (cf. (\ref{ffp}) et [W] Proposition I.4.3) joint ˆ l'ŽgalitŽ (\ref{H'}). Elles sont donc Žgales. Cela prouve la premire ŽgalitŽ. \\ 
Supposons maintenant $(\pi, V)$ tempŽrŽe et notons $f$ la fonction sur $A_M$,  $a \mapsto c_{\xi,v}(a_0 a)$.
Comme $f$  est $A_M$ finie, il existe (cf.  [W] I.2) un ensemble fini ${\mathcal X}$, un entier $d$  tels que  pour tout $\chi \in {\mathcal X}$ , il existe un polynome $P_{\chi,f}$ sur $\a_M$, ˆ coefficients complexes et de degrŽ  infŽrieur ou Žgal ˆ $d$, de sorte que $$f(a) =
 \sum_{\chi \in{\mathcal X}}\chi(a) P_{\chi,f}(H_M(a)).$$ On peut prendre ${\mathcal X}= {\cal E}xp (\pi_P)$. Comme  $(\pi, V)$ est  tempŽrŽe, tenant compte de la premire ŽgalitŽ du Lemme,  les deux membres de la deuxime  sont Žgaux ˆ $$ \sum_{\chi \in {\cal E}xp (\pi_P), Re\chi=0}\chi(a) P_{\chi,f}(H_M(a)).$$ \qed
 \section{Fonctionnelles de Jacquet et terme constant faible}
 \setcounter{equation}{0}
\subsection{Fonctionnelles de Jacquet}
 Soit    $P=MU$ un sous-groupe parolique semi-standard.
$(\si, E)$ une  reprŽsentation lisse unitaire irrŽductible de  $M$. On note $\O$ (resp. $\O_\C$)
 l'ensemble des classes d'Žquivalences des reprŽsentations 
 $\si_\chi:= \si\otimes \chi$,  $ \chi \in X(M)$ ( resp. $X(M)_\C$) , qui est un tore compact (resp. complexe).  On appelle $\O_\C$ (resp. $\O$) l'orbite inertielle (resp. orbite inertielle unitaire) de $\si$. On utilise dans la suite la notion de fonction sur $\O_\C$ (resp. $\O$) ˆ valeurs dans un foncteur, qui dŽpendent des objets concrets $\si_\chi$, ou des reprŽsentations Žquivalentes (resp. unitares Žquivalentes)  ˆ l'une de celles-ci (cf. [D4], section 2.4). On appellera ces dernires objets de $\O_\C$ (resp. $\O$)avec des rgles de transformations pour tenir compte des Žquivalences ( resp. Žquivalences unitaires)   de reprŽsentations. On dispose de la notion de fonction polynomiale, rationnelle (cf. l.c.).\\On note que 
$Wh(\si _\chi ):= Wh(\si_ \chi, U_0\cap M)$ est indŽpendant de $\chi \in X(M) $, car $\chi$  est trivial sur  $ U_0\cap M$. Par restriction des fonctions ˆ $K$, les reprŽsentations $i^G_P  \si_\chi$ admettent une rŽalisation dans un espace indŽpendant de $\chi$, la rŽalisation compacte.  
\ber \label{ji}  Soit $P=MU$ un sous-groupe parabolique anti-standard de $G$,  $P^-=MU^-$ le sous-groupe parabolique opposŽ relativement ˆ $M$.  On note $(\si, E)$ une reprŽsentation lisse de longueur finie de $M$. \\
Il y a un unique isomorphisme entre  $Wh(\si) \to Wh (i^G_P  \si)$ notŽ  $\eta \mapsto \xi(P, \si, \eta)$ (Rodier [R], Casselman-Shalika [CS], Shahidi [Sh], Proposition 3.1) tel que, pour tout $v\in i^G_PE$ ˆ support contenu dans $PP^-$:
$$ \langle \xi(P,\si,  \eta), v \rangle   =  \int_{ U^- } \langle 
\eta,
v (u^-)\rangle   \psi(u^-)^{-1} du^-,  \eta \in Wh(\si).
$$ 
 De plus (cf. [D4], ThŽorme 1), pour tout  $v$ dans l'espace de la rŽalisation compacte et $\eta \in Wh(\si) $, $ \langle \xi (P, \si_\chi, \eta), v\rangle  $  est polynomiale en  $\chi \in X(M)$. \eer 
On dŽfinit les fonctionnelles de Jacquet  pour un sous-groupe parabolique semi-standard de $G$,  $P=MU$, par 
 transport de structure. Celles-ci ne seront qu'un outil pour les preuves et  n'apparaitront pas dans la formule de Plancherel. \ste
 On rappelle  que $W^G$ dŽsigne un ensemble de reprŽsentants dans $K$ du groupe de Weyl, ${\overline W}^G$,  de $G$ par rapport ˆ $M_0$. Si 
$M$ est un sous-groupe de LŽvi d'un sous-groupe parabolique semi-standard de $G$, $P$, on note $W^M= W^G\cap M$ qui est un ensemble de reprŽsentants dans 
$M$ de ${\overline W}^M$.  La longueur des ŽlŽments de $\overline{W}^G$ est dŽterminŽe par le choix de $P_0$. \\On suppose en outre ici $P$ anti-standard. 
 Il existe un  ensemble de reprŽsentants de $\ \overline{ W}^G / \overline{ W}^M $,  $\overline{ W}_M$, dans $\overline{W}^G$ tel que (cf. [War], Proposition 1.1.2.13):  \ber   \label{wM}Tout ŽlŽment $w$ de $\overline{W}^G$ s'Žcrive  sous la forme $w_Mw^M$, avec  $w_M \in \overline{ W}_M$, $w^M \in \overline{ W}^M$,   et tel que la longueur de $w$ soit Žgale ˆ la somme des longueurs de $w_M$ et $w^M$. \eer 
 \ber \label{wP} Soit $P=MU$ un sous-groupe parabolique semi-standard de $G$. On notera
$w_P$ ou parfois seulement $w$,  s'il n'y a pas d'ambiguitŽ,  l'ŽlŽment $w_P$ de $G$ 
tel que
$w_P^{-1}\in W^G$,  $P'= w_P.P$ soit anti-standard   et tel  que $w_P^{-1}$
reprŽsente l'ŽlŽment du groupe de Weyl  de longueur minimum dans
$w_P ^{-1}{\overline W}^{M'}= {\overline W}^{M}w_P^{-1}$. L'unicitŽ de $w_P$ rŽsulte du fait que deux sous-groupes paraboliques anti-standard de $G$ conjuguŽs sont Žgaux.\eer 
Soit $(\si, E)$ une reprŽsentation lisse de $M$.
On dispose de l'isomorphisme $\l(w): i^G_{P}E\mapsto
i^G_{w.P}w E$ entre les reprŽsentations $ i^G_{P}\si$ et $
i^G_{w.P}w\si$  qui ˆ 
$v$ associe $v_{w}$, o $v_w(g)= v(w^{-1}g)$ pour $g\in G$. 
Notons, que comme $w\in K$, pour $v \in i^K_{P\cap K}E$ et tout $\chi \in X(M)$, la restriction de 
$\l(w) v_\chi$ ˆ $K$ est Žgale ˆ $\l(w)v$. On voit aussi que si $\si$ est
unitaire,
$\l(w)$ est unitaire.  \ste
\ber  \label{whp} On dŽfinit: 
  $$Wh(P, \si):= Wh(w_P\si),$$
$$\xi(P, \si, \eta):= \xi(w_P .P, w_P\si, \eta)\circ \l (w_P),  \eta \in Wh(P, \si).$$
On a:
$$Wh(i^G_P \si)=\{\xi(P, \si, \eta)\vert \eta \in Wh(P, \si)\}$$
 \eer
On dŽfinit les intŽgrales de Jacquet par: $$ E^G_P(\si, \eta\otimes v) = c_{\xi, v}  \in  \cu , $$ o  $\xi= \xi(P, \si, \eta)$, $\eta\in Wh(P, \si),v \in i^G_P E$, ce qui, par bilinŽaritŽ permet de dŽfinir $ E^G_P(\si,\phi)$ pour $\phi \in Wh(P, \si\otimes i^G_PE) $.
On notera souvent plus simplement $ E^G_P(\phi)$ au lieu de  $E^G_P(\si, \phi)$.
De (\ref{whp} ), on dŽduit  l'ŽgalitŽ:  
 \beq\label {Ew} E^G_P(\si, \eta\otimes v)=E^G_{w_P.P}(w_P\si, \eta\otimes  \l (w_P) v ), \> v \in i^G_P E, \eta \in Wh(P, \si). \eeq 
 On suppose que $P$ est anti-standard. Soit  $(\si_1, E_1)$ une reprŽsentation de $M$  Žquivalente ˆ $(\si, E) $ et soit   $T$ est  un opŽrateur d'entrelacement bijectif entre $\si$ et $\si_1$.  On note $T^t$ l'application linŽaire de $Wh(\si_1)$ dans $Wh(\si)$ dŽduite de $T$ par transposition. 
 On note $ind T$ l'opŽrateur d'entrelacement induit de $T$. Alors, il rŽsulte du Lemme 4 (iii) de [D4]:  
\beq \label{ET} E^G_{P} (\si, T^t\eta_1\otimes  v) =E^G_{P} (\si_1, \eta_1\otimes (ind  T) v),  v\in i^K_{P\cap K}E,
\eta \in Wh(\si_1).\eeq 
 \subsection{Rappel sur les intŽgrales d'entrelacement et les matrices $B$ }
 On introduit les intŽgrales d'entrelacement, comme dans [W], Thorme  IV. 1.1, avec un changement de notation toutefois. On en rappelle ci-dessous certaines propriŽtŽs.
Soit $P=MU$, $P'=MU'$ deux sous-groupes paraboliques semi-standard de
$G$ de sous-groupe de LŽvi $M$. Soit $\O_\C$ l'orbite inertielle d'une reprŽsentation lisse irrŽductible de $M$.
\ber \label{intertw}
Il existe une fonction rationnelle dŽfinie sur $\O_\C$,  $A(P', P,.)$ ˆ valeurs dans $Hom_G(i^G_P., i^G_{P'}. )$
avec les propriŽtŽs suivantes:\\
 Pour tout  $(\si, E)$ objet de $\O_\C$,  il existe $R\in \R$
tel que, pour tout $\chi \in X(M)_\C$ vŽrifiant $\langle Re \chi, \alpha\rangle    > R$ pour tout  
$\alpha\in
\Sigma(P)\cap \Sigma({P'^-})$, on ait:$$\langle  (A(P', P, \si_\chi)v)(g), {\check e}\rangle   =
\int_{ U\cap U'\backslash U'} \langle v(u'g), {\check  e}\rangle   du',\> vÊ\in i^G_{P}V_{\chi}, {\check e}\in
{\check E}, $$ l'intŽgrale Žtant absolument convergente.
\\ Si $\si $ est tempŽrŽe, on peut prendre
$R=0$ (cf. [W], Proposition IV. 2.1) 
\eer 
La rationalitŽ s'entend dans le sens suivant:
\ber \label{ratintertw}
\ste Il existe une fonction
polyn\^ome sur $X(M)_\C$ non nulle, $b$, telle que pour
tout $v\in i^K_{K\cap P}V$, l'application qui ˆ $\chi\in X(M)_\C$ satisfaisant la condition ci-dessus associe la restriction ˆ $K$ de 
$b(\chi)A(P', P,\si_\chi)(v_\chi)$ est ˆ valeurs dans un
espace vectoriel de dimension finie de $i^K_{P\cap K}E$ et  se prolonge de faon polynomiale en $\chi\in X(M)_\C$ (cf. 
[W], ThŽorme IV.1.1) .\eer 
Il existe une application rationnelle sur $\O_\C$ ˆ valeurs dans $\C$, $j$,  telle que pour tout sous-groupe parabolique semi-standard de $G$, $P$,  de sous-groupe de LŽvi
$M$, on ait (cf.  [W], IV.3.(1)):
\ber   \label{hom} Pour $\si $ objet de $\O_\C$ tel que $A(P, P^-, \si) A(P^-,P, \si)$ soit dŽfini, cet opŽrateur est l'homothŽtie de rapport $j(\si)$.\eer 
On a (cf. [W] IV.3 (3)):
\ber \label{jw} Si $w\in W^G$, $j(w\si)=j(\si)$.Ê\eer
 D'autre part  on obtient facilement un analogue de [W] IV.1 (11)  pour les
adjoints des intŽgrales d'entrelacement.
 \\D'aprs [W] Lemme V.2.1 et  IV. 3 (6), on a:
\ber \label{reelunit} On suppose  que $\O$ est l'orbite inertielle unitaire d'une reprŽsentation  irrŽductible de carrŽ intŽgrable. Alors la fonction $j$ est positive ou nulle sur $\O$ et non identiquement nulle sur $\O$.  On notera $\mu$ la fonction rationnelle sur $\O$, $j^{-1}$. On notera parfois plus prŽcisŽment $\mu^G$ au lieu de $\mu$, $j^G$ au lieu de $j$. \eer 

Si $\alpha $ est ŽlŽment de l'ensemble $\Sigma_{red}(P)$ des racines rŽduites de
$\Sigma(P)$, on note $A_\aa$ la composante neutre du noyau de $\alpha$ dans $A_M$ et $M_\aa$ le centralisateur de $A_\aa$. On note  $j_\alpha(\si)$au lieu de $j^{M_\aa}$. \\
D'aprs [W] IV.3 (4), si $P, P', P''$ sont des
sous-groupes paraboliques semi-standard de
$G$ de sous-groupe de LŽvi $M$, on l'ŽgalitŽ de fonctions rationnelles sur  $\O$:
\ber  \label{aaa}$$A(P'', P', \si) A(P', P, \si)= j(P'', P', P, \si)
  A(P'',P, \si),$$
o $j(P'', P', P, \si)$ est le produit   
des $j_\alpha(\si)$ pour $\alpha \in  \Sigma_{red}(P) \cap
\Sigma_{red}(P'') \cap \Sigma_{red}(P'^-). $\eer
On a aussi:  
\ber \label{poleA} Pour $\aa \in \Sigma_{red}(P)$, les points o  l'application  rationnelle sur  $X(M)_\C$, $\chi \mapsto  j_\aa(\si_\chi)$,  a un p\^ole ou un zŽro sont de la forme $\chi= \chi_ \lambda$  avec $ \lambda$ ŽlŽment d'un nombre fini d'hyperplans de $Ê(\a'_M)_\C$ de la forme $\langle  \lambda, \check{\aa}\rangle   = c.$ \\ Les points o   l'application  rationnelle  sur $X(M)_\C$, $\chi \mapsto A(P', P, \si_{\chi})$  a un p\^ole ou bien o $ A(P', P, \si_{\chi})$ n'est pas inversible  sont de la forme $\chi= \chi_ \lambda$  avec $ \lambda$ ŽlŽment d'un nombre fini d'hyperplans de $Ê(\a'_M)_\C$ de la forme $\langle  \lambda, \check{\aa}\rangle  = c$,  
 avec $\aa \in \Sigma(P') \cap \Sigma(P^-)$(cf. [H],   p. 393).\eer
Par transport de structure, on a, pour $x\in G$, normalisant $M_0$:
\beq \label{xA} \lambda(x) A(P', P, \si) = A(x.P', x.P, x\si) \lambda(x). \eeq 
et:
\ber \label{jaawsi}  Si $\aa$ est une racine rŽduite de $\Sigma(P)$, $w \in W^G$, alors:
$$   j_\aa (\si)= j_{w\aa }(w \si).$$ \eer 
Les intŽgrales d'entrelacement transforment les fonctionnelles de Jacquet en des fonctionnelles de Jacquet, ce qui permet d'introduire les matrices $B$ (cf. [D4], section 5.2):
\ber\label{B} Soit  $\O$ l'orbite inertielle d'une reprŽsentation lisse  irrŽductible de $M$. Il existe une unique application rationnelle  dŽfinie sur $\O$, $
B(P,P',.)$  ˆ valeurs dans $Hom_\C(Wh(P',.), Wh(P,. ))$  telle que l'on ait
l'ŽgalitŽ de fonctions rationnelles sur $\O$: 
$$\xi(P', \si,  \eta) \circ  A(P', P,\si)= \xi(P, \si, B(P,P', \si)\eta), \> \eta \in Wh(P', \si).$$
La fonction  rationnelle sur $X(M)$, $\chi \mapsto B(P,P', \si_\chi)$ n'a de p\^oles qu'en des points ou l'application $\chi \mapsto A(P', P,  \si_\chi)$  a un  p\^ole. 
\eer 
 La rationalitŽ a ici le sens suivant:\\ 
Soit $\si$ un objet  de $\O$. Pour tout $\chi\in X(M)_\C$, $Wh(P,\si_\chi)= Wh(P, \si)$,  et,
avec les notations de (\ref{intertw}), pour tout $\eta\in Wh(P', \si)$, la fonction $\chi \mapsto b(\chi) B(P,P', \si_\chi)\eta$ est une fonction polynomiale sur $X(M)_\C$ ˆ valeurs dans $Wh(P, \si)$.\\
D'aprs (\ref{poleA}) et (\ref{B}), on a: 
\ber\label{poleB}
 Les points o   l'application  rationnelle  sur $X(M)_\C$, $\chi \mapsto B(P, P', \si_{\chi})$  a un p\^ole  sont de la forme $\chi= \chi_ \lambda$  avec $ \lambda$ ŽlŽment d'un nombre fini d'hyperplans de $Ê(\a'_M)_\C$ de la forme $\langle  \lambda, \check{\aa}\rangle  = c$,  
 avec $\aa \in \Sigma(P') \cap \Sigma(P^-)$. \eer 
 La relation qui dŽfinit $B(P',P, \si)$ comme fonction sur $\O$ est la suivante:  
\ber \label{foncB}Soit $(\si,E)$, $(\si_1, E_1)$ deux objets de $\O$ Žquivalents  et $T$ un entrelacement bijectif entre $\si$ et $\si_1$. La transposŽe de $T$, $T^t$, dŽtermine une bijection entre $Wh(P, \si_1)$ et $Wh(P, \si)$ d'une part,  $Wh(P', \si_1)$ et $Wh(P', \si)$ d'autre part et l'on a:
$$B(P,P', \si_1) = (T^t)^{-1} B(P, P', \si) T^t.$$
\eer
\subsection{\label{sssec}Terme constant faible des fonctionnelles et intŽgrales de Jacquet tempŽrŽes}
Soit $P= MU$ et $P'= M'U'$ deux sous-groupes paraboliques semi-standard de $G$. 
  On note $$\mathcal{W}(M'\vert G\vert M)= \{ s \in W^G \vert s.M\subset M'\}.$$
 On surligne pour indiquer l'image dans le groupe de Weyl.  Soit $$ {\overline  W} (M'\vert G\vert M):={\overline  W}^{M'} \backslash {\overline  {\cal W}} (M'\vert G\vert M). $$
On remarque que $$ {\overline  W} (M'\vert G\vert M):= {\overline  W}^{M'} \backslash {\overline  {\cal W}} (M'\vert G\vert M)/ {\overline  W}^{M}. $$
Ceci permet de choisir un ensemble de reprŽsentants $ W (M'\vert G\vert M)$ de  ${\overline  W} (M'\vert G\vert M)$   dans $W^{G}$
tel que: \ber \label{wm'gm} Pour $s\in W (M'\vert G\vert M)$, $s $ est de longueur minimum dans ${\overline W}^{M'}s$ qui contient $s\overline{W}^M$. \eer
De plus, en utilisant un sous-groupe parabolique standard conjuguŽ ˆ $P$, $x.P$, on voit gr\^ace  ˆ [War],  Proposition 1.2.1.10, que \ber \label{P'sP}
$W(M'\vert G\vert M)$ est un ensemble de reprŽsentants d'un sous-ensemble de $P' \backslash G/ P$. \eer  
\berÊ\label{wsP} Soient $P, P'$  deux sous-groupes paraboliques anti-standard de $G$ et soit  $s\in W (M'\vert G\vert M)$. Alors, avec les notations de (\ref{wP}), $w_{s.P} = s^{-1}$. \eer 
 \begin{defi} \label{reg} Soit $M$ le sous-groupe de LŽvi contenant $M_0$ d'un  sous-groupe parabolique semi-standard de $G$ et $(\si, E)$  une reprŽsentation lisse irrŽductible de carrŽ intŽgrable de $M$. On dit que $\si$ est $G$-rŽgulire si:
 \\1) Pour tout $s\in W(M \vert G\vert M)$  avec $s \not=1$, les reprŽsentations  $ s\si$ et $\si$ sont  non Žquivalentes.
 \\2) Si $P, P'$ sont des sous-groupes paraboliques de $G$ ayant $M$ pour sous groupe de LŽvi, les applications rationnelles sur $X(M)_\C$, $\chi \mapsto A(P, P', \si_\chi)$, $\chi \mapsto B(P, P', \si_\chi)$  n'ont pas pas de poles en $\chi=1$ et leurs valeurs en $1$ sont des opŽrateurs inversibles.\\ 3) La reprŽsentation  $i^G_P \si$ est irrŽductible. \end{defi}
 On remarque que: \ber \label{zar}Si $(\si, E)$  une reprŽsentation lisse irrŽductible de carrŽ intŽgrable de $M$, l'ensemble des $\chi \in X(M)_\C$ tel que $\si_\chi$ soit $G$-rŽgulire est un ouvert de Zariski non vide de $X(M)_\C$ et  si $\si$ est unitaire, l'ensemble des $\chi\in X(M)$  tels que $\si_\chi$ soit $G$-rŽgulire est   Zariski-dense dans $X(M)_\C$. \eer 
 
Soit $P=MU$, $P'=M'U'$  deux sous-groupes paraboliques semi-standard de $G$.
On introduit, pour $s \in W (M'\vert G\vert M)$,  les sous-groupes paraboliques de $G$:
\beq \label{ps} P_s= (M'\cap s.P)U', \> {\tilde P}_s= (M'\cap s.P) U'^-.\eeq
Soit $(\si,E)$,une reprŽsentation lisse irrŽductible de carrŽ intŽgrable de $M$, $G$-rŽgulire de $M$. Posons $(\pi, V)= (i^G_P \si,
i^G_PE)$. On dŽfinit une application $\alpha$:
$$ \alpha: V \to \oplus_{s \in W (M'\vert G\vert M)} i^{M'}_{M'\cap s.P} sE,  \> v \mapsto (v_s)_{s \in
W(M'\vert G\vert M)}.$$
par $$v_s(m')= \delta^{-1/2} _{P'}(m') (A(P_s, s.P, s\si) \l (s)v)(m'), \> m'\in M'.$$  
D'aprs [W], dŽbut  de la preuve de la Proposition V.1.1, on a:
\ber  \label{alpha} L'application $\alpha$ se factorise en un isomorphisme de $M'$-modules entre 
$V_{P'}^w$ et  l'espace d'arrivŽe, qui est une somme de $M'$-modules irrŽductibles non
Žquivalents, rŽduite ˆ zŽro si $W(M'\vert G\vert M)$ est vide. On identifie dans la suite $V^w_{P'} $ ˆ l'aide de $\alpha$ avec l'image cet isomorphisme.  \eer 
\ber \label{seta} Soit $\eta \in Wh(P, \si)$. Comme $\l(s)$ entrelace $i^G_P \si$ et $i^G_{s.P} s\si$,  d'aprs (\ref{whp}), il existe un unique $\s \eta \in Wh(s.P, s\si)$ tel que:
$$ \xi(s.P, s\si, \s \eta)= \xi(P, \si, \eta) \circ \l (s^{-1}).$$
\eer 
\begin{theo} \label{thct} Soit $P$ (resp.  $P'$)  un  sous-groupe parabolique semi-standard (resp. standard) de $G$. Soit $\si$ comme ci-dessus. \ste
(i) Pour $s \in W (M'\vert G\vert M)$, on a $Wh({\tilde  P}_s,s\si) = Wh (M'\cap s.P, s\si)$. \ste 
(ii)Soit $\eta\in Wh(P, \si)$.  On note $\xi= \xi(P, \si, \eta)$.  Alors,
dans l'isomorphisme ci-dessus,  
$\xi^w_{P'} $ qui est un ŽlŽment du dual de $V_{P'}$, est    nul si $W(M'\vert G\vert M)$ est vide et  sinon Žgal ˆ $(\xi_s)_{s \in W(M'\vert G\vert M)}$, avec: 
$$\xi_s = \xi(M'\cap s.P, s\si, B({\tilde P}_s, s.P, s\si)\s \eta), $$
l'expression Žtant bien dŽfinie gr\^ace ˆ (i).
\\(iii) Si $P$ est anti-standard, $\s \eta=\eta$. 
\end{theo}
\dem On procde comme  dans la preuve du ThŽorme 3  de [D4].
(i)   est entirement semblable ˆ celle  de (i) du l.c. 
\ste(ii)  Si  $W(M'\vert G\vert M)$ est vide, $V^w_{P'}$ est  rŽduit ˆ zŽro, donc $\xi_{P'}$ est nul. Ceci prouve la premire assertion de (ii). \ste On suppose maintenant que $W(M'\vert G\vert M)$ est non  vide.  On Žcrit: 
$$ \xi_{P'}=(\xi_s)_{s \in W(M'\vert G\vert M)}, $$
o $\xi_s$ est de la forme:
\ber \label{defetas}$\xi_s= \xi(M'\cap s.P, s\si, \eta_s)$,  pour un ŽlŽment $\eta_s$ de $Wh(M'\cap s.P, s \si)$.
\eer 
ll s'agit donc de montrer:
\ber \label{etas} $$\eta_s= B({\tilde P}_s, s.P, s\si) \s \eta .$$ \eer
a) Montrons d'abord:
\ber\label{eta1} Supposons $P$ semi-standard et  $P^-\subset P'$. Alors $\eta_{1_G}= \eta . $\eer
On procde comme dans la preuve du ThŽorme 3 (ii) de [D4].  Cela conduit ˆ l'ŽgalitŽ, pour $v\in V$ ˆ support dans $PP'= P'^-P'$: $$c_{\xi_{P'}, v_{P'}}(a)= \langle \xi(P\cap M', \si, \eta), \si^-(a) v_{1_G}\rangle, a \in A_{M'},$$ o $\si^-= i^{M'}_{P\cap M'} \si$. Comme $\si^-$ est unitaire, on en dŽduit que: $$c_{\xi^w_{PÔ},v^w_{PÔ}}(a)= \langle \xi(P\cap M', \si, \eta), \si^-(a) v_{1_G}\rangle, a \in A_{MÔ}.$$ Evaluant en $1$, on en dŽduit:
$$\langle  \xi_{P'}, v_{P'} \rangle = \langle \xi (P\cap M', \si, \eta), v_{1_G}\rangle.$$  
Tenant compte de [D4], (7.12) et (7.13), on en dŽduit a).\\
b) La fin de la dŽmonstration du ThŽorme est alors identique ˆ celle du ThŽorme 3 de [D4]. \qed 
Soit $P'=M'U'$ un sous-groupe parabolique semi-standard de $G$, $P=MU\subset M'$  un sous-groupe parabolique semi-standard de $M'$ et soit  $Q= PU'$. Soit $(\si, E)$ une reprŽsentation lisse irrŽductible  de $M$.  De l'application: 
$$E^{M'}_P:   Wh (P, \si)\otimes i^{M'}_P  E \to  C^{\infty}(U_0\cap M' \backslash M', \psi),$$ 
se dŽduit, par fonctorialitŽ de l'induction et l'identification canonique de $i^G_{P'}(i^{M'}_PE)$ avec  $i^G_Q E $, une application: 
 $$E_P^{P'}:Wh (P, \si)\otimes  i^G_Q E   \to  i^G_{P'} C^{\infty}(U_0\cap M' \backslash M', \psi).$$
 Si $\phi $ est un ŽlŽment de $  Wh (P, \si)   \otimes  i^G_Q E$, on le regarde comme
fonction sur $G$ ˆ valeurs dans $ Wh(P, \si) \otimes E$.  On a un isomorphisme de $G$-modules entre $i^G_{Q}E$ et 
$i^G_{P'}(i^{M'}_{P}\si ) $. 
L'Žvaluation en l'ŽlŽment neutre dans la deuxime rŽalisation de $i^G_Q E  $,  donne lieu ˆ
une application, notŽe $r_{M'}$,  de $i^G_Q E  $ dans $ i^{M'}_PE$. Soit $v \in i^G_Q E $.   On a (cf. [D4], (7.21)):
 \ber\label{rM'} $$ \rho(m')r_{M' }(v) = \delta_{P'}^{-1/2}(m')r_{M' }(\rho(m')v).$$ \eer On note encore
$r_{M'}$ l'application de $ Wh(P, \si) \otimes  i^G_Q ) $ dans $Wh(P, \si) \otimes  i^{M'}_{M'\cap P}E $ obtenue par tensorisation de   l'identitŽ de $Wh(P, \si)$ avec  $r_{M'}$.  On a:
\beq \label{epindep}[E^{P'}_P(\phi)(1)](m')= [E^{M'}_P(r_{M'}\phi)](m'),
m'\in M'. \eeq
De plus, si $P=M'$,  $i^{M'}_{M'}E$ s'identifie naturellement ˆ $E$. Avec cette identification on a:
\beq  \label{epindepbis} (E^{P'}_{M'} \phi )(g)= E^{M'}_{M'} (\phi (g))\eeq 

 \begin{prop} \label{cteis}  
 Soit $P=MU$(resp. $P'=M'U'$)   un sous-groupe parabolique anti-standard (resp. standard) de $G$, $\si$ une reprŽsentation lisse irrŽductible unitaire de carrŽ intŽgrable  $G$-rŽgulire de $M$.  \\
 Si  $s\in W(M'\vert G\vert M)$, on note $  C(s, P', P, \si) $ l'application linŽaire de $Wh(P, \si) \otimes i^G_PE  $  dans
$ Wh({\tilde P}_s, s \si)\otimes i^G_{P_s} sE  $ dŽfinie  par: $$ C(s, P', P, \si) =  B({\tilde P}_s, s.P, s\si)\otimes (A(P_s, s.P,  s\si) \l(s)). $$
 Alors, pour $ \phi \in  Wh(P, \si)  \otimes i^G_P E$,
$E^G_P(\phi)_{P'}^{ind} =0  $ si $W(M'\vert G\vert M)$ est vide et sinon, avec l'identification de $i^G_{P_s} s E$ avec
$i^G_{P'}(i^{M'}_{M'\cap s.P}s E )$ et celle de $Wh ({\tilde P}_s, s\si)$ avec $Wh(M'\cap s.P, s\si) $ (cf. ThŽorme \ref{thct},  (i)), on a: 
$$E^G_P(\phi)_{P'}^{w,ind} = \sum_{s \in W(M' \vert G\vert
M)}E^{P'}_{M'\cap s.P} (C(s, P', P, \si) \phi), $$
 $$E^G_P(\phi)^w_{P'}= \sum_{s \in W(M' \vert G\vert
M)}E^{M'}_{M'\cap s.P}(r_{M'} (C(s, P', P, \si) \phi)), \> \phi \in  Wh (P, \si)   \otimes  i^G_PE.$$
\end{prop}
 \dem 
 Les deux membres de la premire ŽgalitŽ sont des fonctions sur $G\times M'$. Par Žquivariance, on se rŽduit ˆ dŽmontrer l'ŽgalitŽ en $(1, m')$ pour tout $\phi$ et tout $m' \in M'$. Cette ŽgalitŽ se rŽduit, gr\^ace ˆ (\ref{epindep}),   ˆ la seconde.  Gr\^ ace  ˆ  (\ref{covct} ) et (\ref{rM'}), on se rŽduit ˆ prouver l'ŽgalitŽ ŽvaluŽe en $1$.  On peut se limiter ˆ dŽmontrer l'ŽgalitŽ pour $\phi= \eta \otimes v$, pour $ \eta \in Wh( \si)$, $v\in i^G_PE$.   
 Alors,  avec les notations du ThŽorme prŽcŽdent:
  $$E^G_{P' } (\phi)_{PÔ}(1) =\sum_{s \in W(M'\vert G \vert M)}\langle  \xi_s, v_s\rangle  .$$ En utilisant la dŽfinition de $\xi_s$ et $v_s$, ce ThŽorme montre la deuxime ŽgalitŽ ŽvaluŽe en $1$.   \qed
  
  \section{Transformation de Fourier-Whittaker}
  \setcounter{equation}{0}
\subsection{DŽfinition de la transformŽe de Fourier-Whittaker}
\begin{lem} \label{prodscal}
On suppose que $(\pi, V)$ est une reprŽsentation lisse de carrŽ intŽgrable   et irrŽductible de $G$. Alors:\\
(i) Il existe un unique produit scalaire hermitien sur $Wh(\pi)$ tel que:
$$\int_{A_G U_0\backslash G}c_{\xi,v}(g) \overline{ c_{\xi',v'}(g)} dg= (\xi, \xi') (v, v'), \> \xi, \xi' \in Wh(\pi), v, v' \in V.$$
o l'intŽgrale est absolument converente d'aprs la Proposition \ref{squaresquare}.\\
(ii) Si $\chi$ est un caractre unitaire  non ramifiŽ de $G$, le produit scalaire sur $Wh(\pi_ \chi)=Wh(\pi)$, ne dŽpend pas de
$\chi$.\\
(iii) Si $T$ est un opŽrateur d'entrelacement unitaire avec une autre reprŽsentation  lisse unitaire irrŽductible de carrŽ intŽgrable de $G$, $(\pi_1, V_1)$, l'opŽrateur transposŽ  $T^t$  dŽtermine un opŽrateur unitaire entre $Wh(\pi)$ et $Wh(\pi_1)$. 
\end{lem}
\dem
(i) Il s'agit d'une simple application du Lemme de Schur. \ste
(ii) rŽsulte immŽdiatement de la caractŽrisation du produit scalaire. \\
(iii) est immŽdiat. \qed
 \ste On appliquera  ces notations aux sous-groupes de LŽvi de $G$. 

\begin{prop} \label{fpo}
(i) Soit $P=MU$ un sous-groupe parabolique  semi-standard de $G$ et  $\O$ l'orbite inertielle unitaire d'une reprŽsentation lisse de carrŽ intŽgrable et irrŽductible  de
$M$. Soit $(\si, E)$ un objet de $\O$.  On munit  $Wh(P, \si)$ du produit scalaire dŽfini  gr\^ace au Lemme \ref{prodscal}  et ˆ la dŽfinition de $Wh(P, \si)$.  Par tensorisation avec le produit scalaire d'induite unitaire de $i^G_P\si$, on en dŽduit un produit scalaire sur $Wh(P, \si)\otimes i^G_PE $. \\ Soit $f \in \CCU$. Il existe  un unique ŽlŽment, $\hat{f}(P,\si)$ de $Wh(P, \si)\otimes i^G_PE $ tel que: 
\beq \label{fpsi}(\hat{f} (P, \si), \phi)= \int_{U_0 \backslash G} f(g) \overline{ E^G_P(\si, \eta\otimes  v)(g) } dg,
\>Ê\phi \in Wh(P, \si)\otimes   i^G_PE  .\eeq
l'intŽgrale du membre de droite Žtant convergente d'aprs la Proposition \ref{squaresquare} et (\ref{intff}). 
 (ii) Soit   $\si_1$  une reprŽsentation unitaire Žquivalente ˆ $\si$ de $M$ et $T$ un opŽrateur d'entrelacement unitaire entre $\si$ et $\si_1$. On note $indT $ l'opŽrateur d'entrelacement induit de $T$, qui est unitaire.  Alors on a: 
\ber \label{fpsi1}  $$\hat{f} (P, \si_1) =((T^t)^{-1} \otimes ind T)  \hat{f} (P, \si ), $$ \eer
ce qui dŽfinit $\hat{f}(P, \si)$ comme fonction sur $\O$ ˆ valeurs dans $Wh\otimes i^G_P$.\\
(iii) Soit $P=MU$,$P'=M'U'$ deux sous-groupes paraboliques antistandard de $G$ et $s\in W(M'\vert G\vert M)$.  Alors, pour $(\si,E)$ objet de $\O$:
$$ \hat{f}(s.P, s\si) = (Id \otimes \lambda(s)) \hat{f}(P, \si) . $$
 (iv)
 Soit $(\si, E)$ objet de $\O$.  
Soit $f\in \ccu$ et $g \in G$.  On a: $$(Id \otimes \rho (g))\hat{f} (P, \si)= (\rho(g) f )\hat{} (P, \si).$$ 
\end{prop}
\dem
(i) Pour $f$ fixŽ,  le second membre de  (\ref{fpsi}) dŽfinit une forme antilinŽaire sur $Wh (P, \si )\otimes i^G_PE$. Il est clair que si $ f$ est invariante ˆ droite  par un sous-groupe ouvert compact $H$ de $G$,  celle-ci est invariante par $H$. Par le ThŽorme de reprŽsentation de Riesz, en dimension finie, on en dŽduit (i).\\
(ii)  La relation 
(\ref{fpsi1})  rŽsulte  de la dŽfinition de $\hat{f}$ en (i) et de (\ref{ET}).\\
(iii) rŽsulte de (\ref{Ew}). \\(iv)  rŽsulte de l'ŽgalitŽ $\rho(g) E^G_P( \si, \eta\otimes v)=E^G_P( \si, \eta\otimes \rho(g)v)$, qui est une consŽquence immŽdiate de la dŽfinition  des intŽgrales de Jacquet. \qed 
On appellera $\hat{f}$ la transformŽe de Fourier-Whittaker de $f$, ou plus simplement sa transformŽe de Fourier. 

\begin{rem}

Nous avons dŽjˆ  dŽfini dans [D4] une transformation de Fourier-Whittaker pour les fonctions $C_c^{\infty}$ en utilisant les reprŽsentations cuspidales de $M$ au lieu des reprŽsentations de carrŽ intŽgrable. Ces deux transformations sont diffŽrentes. Nous nous excusons auprs du lecteur de cet abus de terminologie. \end{rem}
 \subsection{\label{majder} Majorations des dŽrivŽes d'intŽgrales de Jacquet } 
  On introduit l'espace $ C^{\infty}(\O,Wh \otimes  i_P^G)$ des fonctions $C^{\infty} $ sur $\O$ ˆ valeurs dans $Wh \otimes i^G_P$. 
Pour tout sous-groupe compact ouvert, $H$,  de $K$ on note $C^{\infty}(\O,Wh \otimes   i_P^G)^H$ l'espace  des fonctions $C^{\infty} $ $H$-invariantes  par l'action rŽgulire droite de $H$. 
Si $(\si,E) $ est un objet de $\O$, cet espace s'identifie ˆ un sous-espace fermŽ  $C^{\infty}(X(M), Wh (\si) \otimes( i_{K, K\cap P}E)^H$ muni de la topologie dŽduite de la topologie naturelle sur $C^{\infty}(X(M))$ . On vŽrifie aisŽment que cette topologie ne dŽpend pas du choix de $\si$. On note $ C^{\infty}(\O,Wh \otimes  i_P^G)$ la rŽunion des espaces $ C^{\infty}(\O,Wh \otimes  i_P^G)^H$ lorsque $H$ varie et on le   munit    de la topologie induite limite inductive. \\ 
On dŽfinit une reprŽsentation lisse de $G$ sur $ C^{\infty}(\O,Wh \otimes  i_P^G)$, $\rho_\bullet$, telle que pour tout $\phi \in C^{\infty}(\O,Wh \otimes  i_P^G)$ et $\si $ objet de $\O$:
$$(\rho_\bullet(g)\phi )(\si)= (Id \otimes \rho(g)) \phi (\si), g\in G. $$
Si $H$ est comme ci-dessus et $g\in G$, $\rho_\bullet (g)$ est un opŽrateur continu entre $\fhi^H$ et $\fhi^{H'}$, o $H'=g.H \cap  K$.\\
On note   $  Pol(\O,Wh \otimes  i_P^G)$ le sous-espace de $ C^{\infty}(\O,Wh \otimes  i_P^G)$ formŽ des applications polynomiales. 
 \begin{lem}  \label{HE} On fixe $H$ comme ci-dessus et on choisit $H' $ comme dans (\ref{H'}). 
 Si $(\sigma, E_\si) $ est un objet de $\O$, on  dŽfinit une application $p_{H'}$ de $Wh(\si)\otimes i^G_P E_\si$ dans $\check{ i^G_P  E_\si }\otimes i^G_P E_\si$, qui ˆ $\eta\otimes v$ associe $e_{H'} \xi(P, \si, \eta)$.  Utilisant (\ref{ji}), on voit que cela une application continue, notŽe encore $p_{H'}$ de $C^{\infty}(\O,Wh \otimes  i^G_P)^{H}$ dans $ C^{\infty} (\O,\check { i^G_P }\otimes i^G_P) ^{H'\times H} $, o on munit cet espace d'une topologie analogue ˆ celle de  $C^{\infty}(\O,Wh \otimes  i^G_P)^{H}$.   Pour $\phi'\in  C^{\infty} (\O, \check {i^G_P }\otimes i^G_P ) $,   soit 
   $E^G_P(\phi'(\si))$ la combinaison linŽaire de coefficients lisses de $i^G_P\si$ associŽe ˆ $\phi'(\si) \in \check{ i^G_PE_\si}  \otimes i^G_PE_\si$.
On a:
 \ber \label{restE} Pour $\si$ objet  de $\O$,  les restrictions ˆ $A_0^-$ de $E^G_P( \phi(\si))$ et $E^G_P ((p_{H'} \phi)(\si))$ sont  Žgales. \eer
 \end{lem}  
 \dem Cela rŽsulte de (\ref{H'}). \qed 
\begin{lem}
  Soit $\phi  \in C^{\infty}(\O,Wh \otimes  i_P^G) $. 
  Alors pour tout sous-groupe parabolique standard, $P'=M'U$ et $m'\in M'$ les applications   $\si\mapsto( E^G_P( \phi(\si) ))_{P' }(m')$ et $\si\mapsto( E^G_P( \phi(\si) ))^w_{P '}(m')$ sont $C^\infty $ sur $\O$. 
 \end{lem}
  \dem D'aprs (\ref{covct}) et une propriŽtŽ analogue pour le terme constant faible, il suffit de  prouver les assertions du Lemme   pour $m'= 1$.
 On suppose que $\phi$ est invariante ˆ droite par un sous-groupe compact ouvert de $G$, $H$. On utilise le Lemme prŽcŽdent.  On conclut g\^ace ˆ l'anlague de notre  
 Lemme pour les coefficients lisses (cf. [W] Lemme VI.2.1). \qed 
  \begin{lem}\label{foncw}
  Soit $w$ une fonction strictement positive sur $G$, invariant ˆ gauche par $U_0$. On suppose que pour tout $g_0 \in G$, il existe une constante $C_{g_0} >0$
telle que: $$w(gg_0^{-1} ) \leq C_{g_0}  w(g), \> g \in G.$$
(i) Soit $D$ un opŽrateur diffŽrentiel ˆ coefficients $C^{ \infty} $ sur $\O$. On suppose que pour tout $\phi \in \fhi$, il existe $C'>0$ telle que pour tout $a \in A_0^-$ et tout $ \si$ objet de $\O$:
$$ \vert w(a)D (E^G_P(\phi(\si) )(a)) \vert \leq C'.$$
Alors, pour tout $\phi \in \fhi$, il existe $C>0$ telle que pour tout $g\in G$, $ \si$ objet de $\O$:
$$ \vert w(g)D (E^G_P(\phi(\si) )(g)) \vert \leq C.$$
\\(ii) Soit $f\mapsto f_\phi$ un entrelacement entre $\fhi$ et $\cu$.  On suppose que  tout sous-groupe compact ouvert $H$ de $G$ il existe une semi-norme continue sur $\fhi^H$, $p'$,  tel que , $$ \vert w(a) f(a)\vert \leq  p'(\phi) , a \in A_0^-,  \phi\in \fhi ^H .$$ Alors il existe une 
semi-norme continue sur $\fhi ^H$, $p$,  tel que:  $$ \vert w(g)  f(g)\vert \leq  p(\phi) , g \in G,  \phi\in \fhi^H . $$
\end{lem} \dem
Prouvons (ii), la preuve de (i) Žtant similaire. \\
Pour $f \in \cu$, on note:  $$\nu(f)= Sup_{g\in G}\vert  w(g)  f(g)\vert , \nu'(f)= Sup_{a\in A_0} \vert w(a)  f(a)\vert , \nu'_- (f) = Sup_{a\in A^-_0} \vert w(a)  f(a)\vert .$$
 Soit $H$ un sous-groupe compact ouvert de $G$ contenu dans $K$. Soit $I$ un sous-ensemble fini de $K$  tel que $IH=K$ et soit $F_0$ un sous-ensemble fini de $M_0$ tel que $G= U_0A_0 F_0 K$. \\
 Il nous faut majorer $\nu (f_\phi)$. D'aprs l'hypothse sur $w$,  la dŽfinition de $I$ et $F_0$, et la propriŽtŽ  d'entrelacement de l'application $\phi\mapsto f_\phi $,  il existe $C>0$ tel que pour  tout ŽlŽment, $\phi $,  de $\fhi$ invariant par $H$ on ait:
 \beq \label{nudleq}  
\nu(f_\phi ) \leq C Sup_{f_0\in F_0, i\in I} \nu' (f_{ \rho_\bullet (f_0i)\phi)}) .\eeq
  D'aprs (\ref{trans}),  il existe $a_0 \in A$ tel que pour tout $\phi$, $f_0, i$ comme ci-dessus, la restriction  ˆ $A_0$ de $f_{ f(\rho(f_0ia_0) \phi}$ est ˆ support  dans $A_0^-$.
 \\Mais d'aprs les propriŽts de $w$, on a: 
 $$\nu'( f)\leq C_{a_0}\nu' ( \rho(a_0) f ), f \in \cu.$$
 Tenant compte de la propriŽtŽ de $a_0$, on a : 
 $$\nu(f_\phi ) \leq CC_{a_0}   Sup_{f_0\in F_0, i\in I} \nu'_- ( f_{ \rho_\bullet (f_0ia_0)  \phi }  ), \phi\in \fhi^H
 . $$
ll  existe un sous-groupe compact ouvert de $G$, $H_1$ tel que pour tout  $\phi \in \fhi$ invariant par $H$, tout $f_0\in F_0$, $i\in I$,  
 $ \rho(f_0ia_0)\phi $ est invariant par $H_1$. D'aps  l'hypothse, on en dŽduit  qu'il existe une semi-norme continue $p'$ sur $\fhi$ telle que:
$$ \nu(f) \leq p' (\rho(fia_0) \phi), \phi \in \fhi^H. $$
En tenant compte de la continuitŽ des opŽrateurs $\rho_\bullet (g)$, on en dŽduit le rŽsultat. 
\qed 
 \begin{lem}\label{deriv} 
 (i)  Soit   $D$ un opŽrateur diffŽrentiel sur $\O$, $\phi \in  C^{\infty} (\O, Wh\otimes i_{G,P})$.  Il  existe $d \in \N$ et $C>0$  tel que: 
  $$\vert ( D (E^G_P( \phi (\si))(g)))\vert \leq \delta_0^{1/2} (m_0(g) )( 1+ \Vert H_0(m_0(g))\Vert )^d, g\in G, \si \in \O $$
    
\end{lem} 
  \dem 
  On choisit un sous-groupe ouvert compact de $G$, $H$,  fixant $\phi$ et $H'$  comme dans (\ref{H'}). On utilise les notations du  Lemme \ref{HE}.   
  D'aprs [D1], Lemme 6, et [W] Lemme II.1.1, il existe $d\in \N$ et $C>0$  tel que:
  \beq \label{coef}  \vert (  1+ \Vert H_0 (a) \Vert)^{-d}  \delta_0 ^{-1/2}  (a) (E^G_P (p_{H'}  \phi (\si) ))  \vert  \leq C, a\in A_0 , \si\in \O\eeq 
  On introduit la fonction $w$ sur $G$ par $w(g)= \vert (  1+ \Vert H_0 (m_0(g) ) \Vert)^{-d}  \delta_0 ^{-1/2} (m_0(g))$. Montrons que : \ber   $w$ vŽrifie les hypothses du Lemme \ref{foncw}.  \eer   Soit $g_0\in G$.  Il suffit de voir qu'il existe un sous-ensemble compact de $M_0$, $\Omega$,  tel que pour tout $g \in G$, $m_0( gg_0) \in m_0(g) \Omega$. L'ensemble  $\{ m_0(kg_0)\vert  k \in K\}$  est  un sous-ensemble compact de $M_0$ qui  a les proprŽtŽs voulues. 
  Alors la combinaison de (\ref{coef}) et des  Lemmes \ref{HE} et  \ref{foncw} (i)  conduit au rŽsultat voulu.  \qed 
  \begin{rem} \label{remderiv}
  Soit $(\si,E)$ un objet de $\O$,  $\eta \in Wh(\si)$ et $v\in i^K_{K\cap P}E$.  Soit $D$ un opŽrateur diffŽrentiel ˆ coefficients $C^{\infty}$ sur $X(M)$.  Pour $\chi \in X(M)_u$, on note $v_\chi $ l'ŽlŽment de $i^G_P ( E_\chi)$ dont la restriction ˆ $K$ est Žgale ˆ $v$. 
  On dŽmontre,  de manire similaire au Lemme prŽcŽdent, qu'il existe $d\in \N$ et  $C>0$ telle que:  $$ \vert ( D (E^G_P (\eta\otimes v_\chi )(g))\vert \leq \delta_0^{1/2} (m_0(g)) ( 1+ \Vert H_0(m_0(g))\Vert )^d, g\in G, \chi \in X(M). $$
\end{rem} 
 Soit $f$ un ŽlŽment de ${\cal A}^w( U_0\backslash G, \psi)$. On note, pour $P'$ sous-groupe parabolique standard  de $G$:  \beq \label{f+} f^+_{P'}:= f_{P'}-f_{P'}^w.\eeq
\begin{lem}\label{E+a}
 Avec les notations ci-dessus, pour tout $\phi\in C^{\infty} (\O, Wh\otimes i^G_P)$,
 il existe $\varepsilon>0 $ et $C>0$ tels que:
$$(E^G_P( \phi(\si)) )^+_{P'} ( a)\leq C e^{- \varepsilon \Vert H_0(a)\Vert }, a \in A_{M'} \cap A_0^{-}.$$
\end{lem}
\dem 
On suppose que $\phi$ est invariante ˆ droite par un sous-groupe compact ouvert de $G$, $H$. On utilise le Lemme \ref{HE} et ses notations.
Alors d'aprs le Lemme \ref{ctct}, on a:
$$(E^G_P(\phi(\si) ))_{P'}  ^+(a)=(E^G_P (p_{H'}\phi (\si)) )_{P'} ^+(a), a \in A_0^-,$$
o dans le membre de droite on utilise les notations de [W] pour les coefficients lisses. Alors le lemme  Lemme VI.2. 3 de l.c. conduit au rŽsultat voulu car, avec les notations de celui-ci, $\Xi^{MÔ}(a)= 1$ pour $a\in A_{M'}$. \qed
\subsection{ Premires propriŽtŽs de la transformŽe de Fourier-Whittaker}
 Soit $P=MU$ un sous groupe parabolique anti-standard de $G$ et $\O$ l'orbite inertielle  d'une reprŽsentation lisse irrŽductible et de carrŽ intŽgrable  de $M$.   
\begin{prop}  \label{Ehatf} 
Soit $(\si, E)$ objet de $\O$.  Soit $f\in \CCU$, $g \in G$.  \\ 
(i) On a: $(\rho_\bullet (g)\hat{f}) (P, \si)= (\rho(g) f )\hat{} (P, \si)$. \ste 
(ii)  On rappelle que $Wh(\si) \otimes i^G_PE $ est muni du produit scalaire  obtenu par produit tensoriel du produit scalaire sur $Wh(\si)$ et sur $i^G_PE$. Si  $f,f' \in \cu$ et $f' \in \cu$,  sont tels que $ff'$ est intŽgrable sur $U_0\backslash G$, on pose  $(f,f')_G= \int_{U_0\backslash G}f(g) \overline{f' (g) }dg$. Alors: 
$$ (f, E^G_P(\phi))_G= ({\hat f} (P, \si), \phi), \phi \in Wh(\si) \otimes i^G_PE.$$
\end{prop} 
\dem 
(i) est une reformulation de la Proposition \ref{fpo} (iv) et
(ii) rŽsulte immŽdiatement de la dŽfinition de ${\hat f}$. \qed
\begin{prop} \label{fpoo}
On retient les notations prŽcŽdentes. Pour $f \in \CCU$ on note $\hat{f}_{P, \O}$, la restriction de $\hat{f} (P,.)$ aux objets de $\O$.\\
(i)Alors $\hat{f}_{P, \O} $ est un ŽlŽment de $ \fhi$. \\
(ii) De plus l'application $ f \mapsto \hat{f}_{P, \O}$ est continue de $\CCU$ dans $\fhi$. 
\end{prop}
\dem
(i) On utilise la dŽfinition de $\fhi$ (cf. section \ref{majder}). Avec les notations de la Remarque \ref{remderiv}, il suffit de voir que, pour $\si$ objet de $\O$, $v\in i^K_{K \cap P} \si$, $\eta\in Wh(\si)$ et $D$ opŽrateur diffŽrentiel ˆ coefficients $C^\infty$ sur $X(M)$, l'application $\chi \to (\hat{f}_{P,\O} (\si_\chi),  \eta \otimes v_\chi)$  est $C^\infty$. Mais d'aprs la Proposition prŽcŽdente:
$$  (\hat{f}_{P,\O} (\si_\chi),  \eta \otimes v_\chi)= \int_{U_0 \backslash G} f(g)\overline { E^G_P ( \eta \otimes v_\chi)(g)} dg.$$
La remarque \ref{remderiv} et la dŽfinition de $\CCU$ permettent d'utiliser le ThŽorme de dŽrivation sous le signe somme pour achever de prouver (i).\\
 (ii) Soit $H$ un sous-groupe compact ouvert.  Toujours en utilisant le ThŽorme de dŽrivation sous le signe somme,  on obtient une majoration  de  $\vert D(\hat{f}_{P,\O} (\si_\chi),  \eta \otimes v_\chi) \vert$  par une semi-norme continue de $\CCU^H$ ŽvaluŽe en $f$. Ceci achve de prouver la continuitŽ voulue. 
\section{Paquets d'ondes et leurs transformŽes  de Fourier-Whittaker  }
\setcounter{equation}{0}
\subsection{Paquets d'ondes dans l'espace de Schwartz}
\begin{prop}\label{waveS}
Soit $\O$  l'orbite inertielle unitaire d'une reprŽsentation lisse irrŽductible de carrŽ intŽgrable  de $M$. On  fixe sur $X(M)$ la mesure de Haar  de masse totale 1. On munit $\O$ d'une mesure $X(M)$-invariante et telle que pour tout objet de $\O$, $(\si, E)$, l'application de $X(M)$ dans $\O$, qui ˆ $\chi $ associe $[\si_\chi] $, prŽserve localement les mesures. Si  $\phi\in C^{\infty}(\O,Wh \otimes  i^G_P)$,  on dŽfinit $f_\phi \in \cu$ par:
  $$f_\phi(g):= \int_{\O}\mu(\si)  E^G_P(\phi(\si))(g) d\si, \> g \in G$$
 qu'on appellera paquet d'ondes de $\phi$.\ste  
 (i) On a, pour $g \in G$, $\rho(g) f_\phi= f_{\rho_\bullet (g)\phi}$. \ste 
 (ii) De plus  $f_\phi $ est ŽlŽment de $\CCU$.
 L'application $\phi\mapsto f_\phi $ de $ C^{\infty}(\O,Wh \otimes  i^G_P)$ dans $\CCU$  est continue. 
 \end{prop}
 \dem
 (i) est immŽdiat. 
\\Montrons (ii). 
On pose:
$$ \nu_{d, A_0^-}(f )= sup_{a \in A_0^-} \delta_{P_0} (a)^{-1} (1+ \Vert H_{M_0} (a)\vert )^d \Vert f(a)\vert, f\in \cu .$$
D'aprs le Lemme \ref{foncw} (ii), il suffit, pour tout sous-groupe compact ouvert, $H$, de $G$,  de majorer $ \nu_{d, A_0^-}(f_\phi )$ pour $\phi\in\fhi^H:$
\\On va montrer que la restriction de $f_\phi$ ˆ $A_0^-$ coincide avec la restriction d'un   paquet d'ondes sur le groupe ( cf [W] VI.13). On utilise les notations et le rŽsultat du Lemme \ref{HE}.  
  Dans [W], section VI.3, les  paquet d'ondes de coefficients lisses de reprŽsentations de  $G$ ont ŽtŽ dŽfinis.  D'aprs le Lemme \ref{HE}, on a:
 \ber  Pour $\si \in \O$,  les restrictions ˆ $A_0^-$ de $E^G_P( \phi(\si))$ et $E^G_P( (p_{H'} \phi) (\si))$ sont  Žgales. \eer
 D'o l'on dŽduit:
 \ber \label{rest} Les restrictions de $f_\phi$ et $f_{p_{H'}\phi}$ ˆ $A_0^{-} $ sont Žgales, \eer 
o dans le second membre, il s'agit des paquets d'ondes de [W].  
 D'aprs [W], Proposition VI.3.1, l'application $\phi' \mapsto f_{\phi'}$ est une application continue de $ C^{\infty} (\O, (\check {i^G_P }\otimes {i^G_P})^{H'\times H'} )$ dans l'espace de Schwartz de $G$, ce qui conduit notamment ˆ une majoration de $\nu_{d, A_0^-} (f_{\phi'})$.  Joint ˆ la continuitŽ de l'application $\phi\mapsto p_{H'} \phi$ et ˆ   (\ref{rest}), on conclut ˆ une  majoration de $\nu_d(f_\phi)$, o $\nu_d$ a ŽtŽ dŽinie  en (\ref{nud}).  Ceci achve de prouver la 
Proposition. \qed
\begin{rem} \label{refwaveS}
Soit $\si$ un objet de $\O$,  $\eta \in Wh(\si)$ et $v\in i^K_{K\cap P}\si$.   Pour $\chi \in X(M)$, on note $v_\chi $ l'ŽlŽment de $i^G_P ( \si\otimes \chi)$ dont la restriction ˆ $K$ est Žgale ˆ $v$. On dŽmontre de manire similaire ˆ la Proposition que,
 pour tout $\phi\in C^\infty(\O)$, $\Phi=\int_{X(M)} E^G_P ( \eta\otimes v_\chi) \phi(\chi) d\chi $ est un ŽlŽment de $\CCU$ et que l'application $\phi\to \Phi$ est continue de $C^\infty(\O)$ dans $\CCU$.
\end{rem}
 \subsection{Transformation unipotente et transformŽe de Fourier-Whittaker }
\begin{lem}\label{hsi} Soit $P=MU$ un sous-groupe parabolique anti-standard de $G$.\\ (i) Il    existe $C>0$ tel que:
\beq  (  (1+Ê\underline{\si} (u))  \leq C ( 1+ \Vert H_0(m_0(u) \Vert), u \in U .\eeq (ii) Il existe $d\in \N$  tel que l'intŽgrale: $$\int _U \delta^{1/2}_0(m_0(u) )(1+\Vert H_0(m_0(u)\Vert )^{-d}du  $$ soit absolument convergente.
\end{lem} 
\dem   (i) rŽsulte de  [W], Lemme II.3.4, utilisŽ pour $P_0$, de sorte que $U$ est contenu dans $U_0^-$,  et I.1(6). \\ (ii) rŽsulte alors du Lemme II.4.2 de [W]. 
\begin{lem}\label{jint} Soit $P=MU$ un sous-groupe parabolique anti-standard de $G$, $(\si, E)$ une reprŽsentation lisse irrŽductible et de carrŽ intŽgrable de $M$. 
\\Soit $\varepsilon> 0 $. On note $$X(M)^\varepsilon:\{ \chi=\chi'\delta_P^\varepsilon\vert \chi' \in X(M)\}.$$ Soit  $v\in i^K_{P\cap K}E$. \\(i)  Il existe $C>0$ tel que:
$$  \int_{U^-}\vert \langle \eta, v_\chi  (u^-) \rangle 
\vert du^-  \leq C, \chi\in X(M)^{\varepsilon} .$$ (ii) Pour $\chi\in X(M)^{\varepsilon}$, on a: $$ \langle \xi(P, \si_\chi , \eta),  v_\chi \rangle = \int_{U^-} \psi(u^{-})^{-1} \langle \eta, v_\chi  (u^-)  \rangle 
du^- $$  o l'intŽgrale est absolument convergente. 
\end{lem}
Commencons par dŽmontrer (i). 
On a: \beq \label {langeta}  \langle \eta, v_\chi   (u^-) \rangle= \chi (m_P({u^-}))\delta_P^{1/2} (m_P(u^-) )\langle \eta , \si(m_P(u^-)v(k(u^-))\rangle, u^-\in U^-. \eeq 
Soit $P'_0= (P_0\cap M)U$. C'est un sous-groupe parabolique minimal semi-standard de $G$ de sous-groupe de LŽvi $M_0$. D'aprs le Lemme \ref{squaresquare}, $\eta $ est tempŽrŽe. Tenant compte de (\ref{deftemp} ) et de l'ŽgalitŽ $m_{P_0\cap M}(m_P(u^-))= m_{P'_0}(u^-)$, on voit qu'il existe  $C'>0$ et $d' \in N$ tels que pour tout $k \in K$ et $m \in M$:\beq \label{esv}\langle\eta, \si(m_P(u^-)) v(k) \rangle  \leq C'\delta_{P_0 \cap M}^{1/2}( ( m_{P'_0}(u^-)) ( 1+\Vert  H_0 ( m_{P'_0}(u^-))\Vert )^{d'}Ê.\eeq  
Par ailleurs, par exponentiation de la premire inŽgalitŽ du   Lemme II.3.4 de [W], on voit qu'il existe $C_1, C'_1>0$ tels que: 
$$ \delta_P (m_P(u^-))^\varepsilon  \leq C_1 e^{-C_1' (1+ \underline{\si}(u^-)) }, \> u^-\in U^-.$$
Donc: \ber \label{Csec} Pour tout $d''\in  \N$,  il existe $C''>0$ tel que $$ \delta_P (m_P(u^-))^\varepsilon \leq  C''(1+ \underline{\si}(u^-)) ^{-d''}, u^-\in U^-.$$ \eer  Utilisant   (\ref{esv}), (\ref{Csec}) et le Lemme \ref{hsi} (i), en y remplaant $P_0$ par $P'_0$ et $U$ par $U^-$, on voit que 
pour tout $d\in \N$, il existe $C>0$ tel que pour tout $\chiÊ=  \chi' \delta_P^\varepsilon, \chi'\in X(M)$, $u\in U^-$:
$$ \vert  \langle \eta, v _\chi (u^-) \rangle \vert \leq C  \delta_{P_0\cap M}^{1/2} ( m_{P'_0}(u^-)) (1+ \Vert H_0(m_{P'_0}(u^-) )\Vert) ^{-d}.$$
On remarque que, d'une part :$$ \delta_{P'_0}(m_0)=  \delta_{P_0\cap M}(m_0)\delta_P(m_0), \> m_0 \in M_0, $$ et d'autre part: $$m_0 (m_P(u^-))= m_{P'_0} (u^-), \delta_P( m_P (u^-))= \delta_P( m_{P'_0} (u^-)).$$ Donc $$\delta_{P'_0}(m_{P'_0}(u^-)  )= \delta_{P_0\cap M}(m_0(m_P(u^-) ) \delta_P(m_P(u^-)).$$
De l'inŽgalitŽ prŽcŽdente,  de (\ref{langeta}) et du Lemme \ref{hsi} (ii), appliquŽ ˆ $P'_0$ au lieu de $P_0$ et $U^-$ au lieu de $U$, on dŽduit (i). \\
Alors le second membre de l'ŽgalitŽ de (ii)  dŽfinit un ŽlŽment de $Wh(i^G_P \si)$. Celui-ci est Žgal ˆ $\xi(P, \si, \eta)$ d'aprs (\ref{ji}). \qed  
 \begin{prop}\label{fpint}
(i) Si $f \in \CCU$  et $P=MU $ est un sous-groupe parabolique  anti-standard de $G$,  on dŽfinit: $$f^P(m)=
\delta_P^{1/2}(m)
\int_U f(mu) du, \> m \in M ,$$
l'intŽgrale Žtant absolument convergente. \\
 (ii) On appelle $f^P$ la transformŽe unipotente de $f$ relativement ˆ $P$. Soit $H$ un sous-groupe ouvert compact de $G$. Soit $d $ comme dans le Lemme \ref{hsi} (ii). Il existe une semi-norme continue, $p$, sur $\CCU^H$ telle que:
 $$ \Vert f^P(m)\Vert \leq  p(f) \delta^{1/2} _{P_0\cap M}(m_0(m)) (1+\Vert H_0(m_0(m)) \Vert )^d, m\in M .$$
\end{prop} 
\dem   Soit $f\in \CCU$. 
Soit $d$ comme dans le Lemme \ref{hsi} (ii). D'aprs (\ref{nud}), on a:
$$\vert f(mu)\vert \leq \nu_d(f)\delta_P^{1/2}(m)  \delta_0^{1/2} (m_0 (mu)) (1+ \Vert H_0(m_0(mu) \Vert)^{-d}, m\in M, u\in U.$$
Donc $$ f^P(m) \leq \nu_d(f) I(m), $$
o: $$I(m)=  \int_U \delta_P^{1/2}(m)\delta^{1/2}_0 (m_0 (mu)) (1+ \Vert H_0(m_0(mu) \Vert)^{-d}du.$$ 
On voit que  $I(m)$ est invariante ˆ droite par $K\cap M$, en changeant $u$ en $kuk^{-1}$ pour $ k \in  K\cap M$. Pour majorer $I(m)$, on  peut se ramener ˆ  $m= u_1 m_1$ avec $m_1 \in M_0$, $u_1 \in U_0$. On voit immŽdiatement qu'alors $m_0( mu)= m_0(m) m_0(u)$. Donc \beq \label{ddd} \delta_P (m) \delta_0(m_0 (mu))=\delta _{P_0\cap M} (m_0(m))  \delta_0 (m_0(u)), u\in U. \eeq  On pose $X= H_0(m_0( m))$,
$Y= H_{0}(m_0(u))$,
$Z= H_0(m_0(mu)) $. On a $Y= Z -X$. L'inŽgalite triangulaire implique $(1+\Vert Z\Vert ) (1+\Vert X\Vert ) \geq 1+\Vert Y\Vert  $.  On en dŽduit:
 \beq \label{HHH} ( 1+  \Vert H_0(m_0(mu))\Vert )^{-d} \leq  (1 + \Vert  H_0(m_0(u)) \Vert )^{-d}  ( 1 + \Vert H_0(m_0(m)) \Vert )^d.\eeq 
 Tenant compte de (\ref{HHH}),  et de (\ref{ddd}), on en dŽduit:
$$I(m) \leq  \delta^{1/2} _{P_0\cap M}(m_0(m))(1+ \Vert H_0(m_0(m)) \Vert ) ^d  \int_U \delta_{0} ^{1/2} (m_0(u))(1+ \Vert H_0(m_0(u))\Vert )^{-d} du .$$
L'intŽgrale, $I$, du membre de droite  Žtant convergente d'aprs le Lemme \ref{hsi}, on a:
$$\vert f^{P} (m)\vert \leq  I \nu_d (f)  \delta^{1/2} _{P_0\cap M}(m_0(m)) (1+ \Vert H_0(m_0(m)) \Vert ) ^d, m\in M.$$
Ceci achve de prouver le Lemme.\qed\\
On dŽfinit 
$f^{P, ind}$  par:
\beq \label{fPind}f^{P, ind} (g, m)=(\rho (g) f)^{P}(m), \> g\in G, m\in M. \eeq
Alors $f^{P, ind}\in i_P^G\cmu$.
   \begin{prop}\label{findfhat} Soit $P=MU$ un sous-groupe parabolique anti-standard de $G$, $(\si, E)$ une reprŽsentation lisse irrŽductible de carrŽ intŽgrable de $M$.  Soit $f \in \CCU$ tel que pour tout $g\in G$,  $f^{P,ind}(g) $ est ŽlŽment de  $\CCMU$. 
  Alors: 
 $$ {\hat f }(P,\si)(g)= (f^{P,ind}(g) {\hat )} (M, \si), g \in G,$$ 
 ŽgalitŽ qu'on rŽŽcrit:
  $${\hat f }(P,\si)=(f^{P,ind} {\hat )} (M, \si).$$
  \end{prop}
 \dem 
 \\Soit $ v \in i^K_{K\cap P} E$,    $\eta \in Wh(\si)$ et  $ \phi\in C^\infty(X(M))$.  Pour $\chi \in X(M)$, on  note
$v_\chi$ l'ŽlŽment de  $i^G_{P}E_\chi$ dont la restriction ˆ $K$ est Žgale ˆ $v$.  Soit 
$$I :=\int_{X(M)}  ({\hat f}(P, \si_\chi),  \eta \otimes v_ \chi)\phi(\chi) d\chi,$$
l'intŽgrale Žtant convergente d'aprs la Proposition \ref{fpoo} (i). On utilise  la Proposition \ref{Ehatf} (ii) et la Remarque \ref{remderiv} pour utiliser  le ThŽorme de Fubini et en dŽduire:
\beq \label{defI} I=(f,  \int_{X(M)} E^G_P(\eta \otimes  v_\chi  )\phi( \chi) d\chi )_G, \eeq 
c'est ˆ dire: 
  $$I= \int_{U_0\backslash G} f(g)\int_ {X(M)}  {\overline {\langle  \xi(P, \si_\chi , \eta), i^G_P\si (g)v_\chi \rangle  \phi(\chi) } }d\chi dg.$$
  On suppose maintenant que $\phi $ ŽlŽment de l'espace $Pol( \O)$ des fonctions polynomiales sur $\O$ et que $f$ est ˆ support compact modulo $U_0$. Soit $\varepsilon>0$.  On peut remplacer l'intŽgrale sur  $X(M)$ par l'intŽgrale sur  $X(M)^\varepsilon $, car $\hat{f}(P, \si_\chi)$ est polynomiale (cf. la preuve de la Proposition 5 (i) de [D4] qui s'adapte ici).
 D'aprs [D4], (8.19), on a: 
  \ber  \label{tauu}
  Soit $\Omega$ un sous-ensemble compact modulo l'action ˆ gauche de $U_0$.
  Il existe une fonction continue ˆ support compact sur $G$
  telle pour tout 
  $g\in \Omega$, $\int_{U_0} \tau (ug)dg = 1$.
\eer   
   \\ Appliquant ceci au support de $f$ et reprenant le calcul de $I$, on trouve:
    $$I= \int_{ G} f(g) \tau (g)\int_{X(M)^\varepsilon}   \overline{\langle  \xi(P, \si_\chi, \eta), i^G_P \si(g) v_\chi  \rangle \phi(\chi) }d\chi dg.$$
 D'aprs le Lemme \ref{jint}, pour $\chi \in X(M)^\varepsilon$, on a:
\beq \label{xii} \langle \xi(P, \si_\chi, \eta), i^G_P\si_\chi (g) v_\chi \rangle = \int_{U^-} \psi(u^{-})^{-1}\langle \eta, v_\chi  (u^-g)  \rangle  
du^- , \eeq o l'intŽgrale est absolument convergente et  pour tout $g\in G$: \beq\label{xiii}  \int_{X(M)^\varepsilon} \int_{U^-}\vert  \psi(u^{-})^{-1}\langle \eta, v_\chi (u^-g)  \rangle \vert du^- <+\infty .\eeq 
Donc:  $$ I= \int_G \tau(g) f(g) \int_{X(M) ^\varepsilon} \int_{U^-}\overline { \ \psi(u^{-})^{-1}\langle \eta, v_\chi (u^-g)  \rangle \phi(\chi) }
du^- d\chi dg.$$ 
L'application  $\tau f$  est une application continue sur $G$, ˆ support  compact. Ce support est donc contenu  dans un nombre  fini de $H$-classes a droite, o $H$ est un sous-groupe compact ouvert de $G$ fixant $v$. Il rŽsulte de (\ref{xiii}) que      le ThŽorme de
Fubini s'applique et l'intŽgrale sur $G\times X(M)^\varepsilon\times U^-$ est absolument convergente.  Il en ira de m\^eme des intŽgrales suivantes car elles sont obtenues gr\^ace au ThŽorme de Fubini ou par un changement de variable.\\
 En utilisant l'ŽgalitŽ $f(u^-g) =\psi(u^-) f(g) $ et le fait que $\psi$ est unitaire,  on a:  
  $$I= \int_{X(M)^\varepsilon} \int_{U^-}   \int _G \tau(g) f(u^-g)\overline{ \langle \eta, v_\chi (u^-g) \rangle \phi(\chi) }dg  du^- d\chi .$$
 On change $g $ en $u^-g$ dans l'intŽgrale intŽrieure:  
  $$I=\int_{X(M)^\varepsilon}  \int_{U^-}   \int _G \tau((u^-)^{-1} g) f(g) \overline{ \langle \eta, v_\chi(g) \rangle \phi(\chi)} dg   du^- d\chi .$$
Tenant compte de l'invariance  ˆ droite par $U_0\cap M$ de $f(g) \overline{\langle \eta, v_\chi (g) \rangle 
}$, on obtient: 
$$I= \int_{X(M)^\varepsilon } \int_{U^-} \int _ {U_0\cap M\backslash G}\int_{U_0 \cap M } \tau((u^-)^{-1} u_0g)
f(g)\overline{ \langle
\eta, v_\chi (g) \rangle  \phi(\chi)}du_0 dg  dh d\chi .$$
Transformant la succession des intŽgrales sur $U^-$ et $U_0\cap M$ en une intŽgrale sur
$U_0$ et en utilisant les propriŽtŽs de
$\tau$ (cf. 
(\ref{tauu})), on en dŽduit: 
 $$I=\int_{X(M)^\varepsilon} \int _{U_0\cap M \backslash G}  f(g) \overline{\langle \eta, v_\chi (g) \rangle \phi(\chi)} dg d\chi.  $$
 Utilisant la formule intŽgrale (\ref{iwas}) et tenant compte du fait que $v_\chi $ est invariante
ˆ gauche par
$U$, il s'ensuit que $I$ est Žgal ˆ:
 $$\int_{X(M)^\varepsilon}  \int_{U\times (U_0\cap M \backslash M) \times K} f(umk) \overline{\langle
\eta, v_\chi(mk) \rangle \phi(\chi)}Ê\delta_P^{-1} (m) du dm dk  d\chi.$$
Mais on a:  
 $$\langle \eta, v_\chi  (mk) \rangle =  \delta_P^{1/2}(m) E^M_M( \si_\chi, \eta \otimes v_\chi(k) )(m) $$
 et $$\int_U f(umk)du =  \delta_P ^{1/2} (m)   [f^{P, ind} (k)] ( m). $$
 Donc les fonctions modules disparaissent et l'on a, gr\^ace ˆ une utilisation du ThŽorme de Fubini, rendue possible par la convergence absolue des intŽgrales dŽjˆ mentionnŽe:
 $$I=   \int_{(U_0\cap M \backslash M)
\times K } [f^{P, ind} (k)] ( m)
\int_{X(M)^\varepsilon} \overline { E^M_M(  \si_\chi, v_\chi (k) \otimes \eta)(m) \phi(\chi)} d\chi dm dk.$$
Par holomorphie, on a: $$\int_{X(M)^\varepsilon } E^M_M( \si_\chi, v(k) \otimes \eta)(m) \phi(\chi) d\chi= \int_{X(M)} E^M_M( \si_\chi, v(k) \otimes \eta)(m) \phi(\chi) d\chi .$$
 D'aprs la Proposition \ref{waveS}, appliquŽe ˆ $M$ au lieu de $G$ et $P$,  le premier membre de cette ŽgalitŽ dŽfinit un ŽlŽment de  $\CCMU$ qu'on note   $\Phi_k$. On a donc:
\beq \label{Iegal} I= \int_{( U_0\cap M \backslash M)
\times K } [f^{P, ind} (k)] ( m)
\overline {\Phi_k(m)}  dm dk .\eeq 
D'abord, $f $ Žtant fixŽe,  d'aprs (\ref{defI}) et la Remarque \ref{remderiv}, l'application $\phi\to I $ est une forme linŽaire continue sur $C^{\infty} (\O)$. Par ailleurs,  le second membre, $J$, de  (\ref{Iegal}) est bien dŽfini pour $ \phi \in C^\infty(\O)$,   car  pour tout $k \in K$, $f^{P,ind}(k) \in  {\cal A}^w(U_0\cap M \backslash M, \psi)$    et $\Phi_k \in \CCMU$. De plus l'application $\phi \mapsto \Phi_k$ est continue de $ C^{\infty} (\O) $ 
dans $\CCMU$ d'aprs la Remarque \ref{refwaveS}.  Donc l'application  $\phi\mapsto J$  est une forme linŽaire continue sur $C^{\infty}(\O)$. Par densitŽ, l'ŽgalitŽ (\ref{Iegal}) est  vraie pour tout $\phi\in C^\infty(\O)$. On fixe $\phi \in C^{\infty}(\O)$. 
D'aprs (\ref{defI}) et la Remarque \ref{remderiv}, l'application $f\mapsto I$ est continue sur $\CCU$ . De m\^eme l'application $f\mapsto J$ se prolonge en une forme linŽaire continue   sur $\CCU$: cela rŽsulte du fait que pour tout $k\in K$, $\Phi_k\in \CCMU$ et de la  Proposition \ref{fpint} (ii). Donc, par densitŽ (\ref{Iegal}) est vrai pour tout $\phi \in C^{\infty}(\O)$ et pour tout $f\in \CCU$.\\ 
Supposons que $ f$ vŽrifie l'hypothse du ThŽorme. On peut alors, gr\^ace ˆ  la Remarque  \ref{remderiv}, utilisŽe pour $M$ au lieu de $G$,  appliquer le ThŽorme de Fubini pour obtenir, pour tout $ \phi \in C^{\infty}(X(M))$:
 $$ \int_{XM)}(\hat{f}(P, \si_\chi), \eta\otimes v_\chi)d\chi   = \int_{X(M)   \times U_0\cap M \backslash M\times K } [f^{P, ind} (k)] ( m)
\overline { E^M_M(  v_\chi(k) \otimes \eta)(m)  \phi(\si) }dm dkd\chi. $$
Comme ceci est vrai pour tout $\phi\in C^{\infty} (X(M))$, on en dŽduit en particulier l'ŽgalitŽ:
$$(\hat{f}(P, \si ), \eta\otimes v)=  \int_{K\times U_0\cap M\backslash  M} ( f^{P,ind}(k) (m)\overline{ E^M_M (\si_\chi,\eta \otimes v(k)) (m) } dk dm,$$
o on note encore $v$ l'ŽlŽment de $i^G_P E$ dont la restriction ˆ $K$ est Žgale ˆ $v$.  
D'aprs la dŽfinition de la transformŽe de Fourier-Whittaker  pour $M$, on en dŽduit:
$$(\hat{f}(P, \si ), \eta\otimes v)= (( f^{P,ind}(k))\hat{}, \eta \otimes v).$$
Comme ceci est vrai pour tout $v\in i^G_PE$,  on en dŽduit:
 $${\hat f }(P,\si)=(f^{P,ind} {\hat )} (M, \si).$$ Ceci achve la preuve de la Proposition.\qed 
  \subsection{Assertion A et relations de Maass-Selberg}
  Soit $P=MU$ un sous-groupe parabolique anti-standard de $G$ et $P' =MU'$ un sous-groupe parabolique semi-standard de $G$ de sous-groupe de LŽvi  $M$. Soit $\O$ l'orbite inertielle unitaire  d'une reprŽsentation lisse irrŽductible de carrŽ intŽgrable de $M$.  
 On sait (cf.  [W],  V.1 (11)),  que l'on a l'ŽgalitŽ de fonctions rationnelles sur $\O$:
\beq\label{adjA} A(P', P,  \si)^*= A(P, P',  \si)  .\eeq
{\bf Nous dŽmontrerons plus loin (cf. ThŽorme \ref{adjB}),  l'assertion suivante dite Assertion A:  }
\ber \label{B*} On a l'ŽgalitŽ de fonctions rationnelles sur $\O$,
$$ B(P, P', \si )^*= B(P', P,  \si) .$$,\eer
 \begin{prop}\label{MS}  Supposons l'assertion A vraie. Soit $P= MU$ (resp. $P'=M'U'$ ) un sous-groupe parabolique antistandard (resp. standard)   de $G$ .  Soit $\O$ l'orbite inertielle unitaire  d'une reprŽsentation lisse irrŽductible de carrŽ intŽgrable de $M$. Soit $(\si,E)$ un ŽlŽment $G$-rŽgulier de $\O$. Soit $\phi\in Wh(\si)\otimes i^G_P E$ et $s \in W(M'\vert G\vert M)$. Alors on a:
 $$ \mu^G(\si) \Vert C(s, P',P, \si) \phi\Vert^2= \mu^{M'} (s\si) \Vert \phi\Vert ^2.$$
 \end{prop}
 \dem 
 On peut se limiter ˆ prouver l'assertion pour $\phi= \eta\otimes v$, avec $\eta\in Wh(\si)$, $v\in i^G_P E$. 
 Alors d'aprs la dŽfinition des fonctions $C$ (cf. Proposition \ref{cteis}):
  $$ \Vert C(s, P',P, \si) \eta\otimes v \Vert^2 = \Vert B({\tilde P}_s, s. P, s\si)\eta \Vert ^2 \Vert A(P_s,s.P, s\si)v\Vert^2
  .$$
  D'aprs la formule d'adjonction pour les intŽgrales d'entrelacement (cf (\ref{adjA})) et l'assertion A qu'on suppose vraie, on a:
  $$\Vert  B({\tilde P}_s, s. P, s\si)\eta \Vert ^{2}=x \Vert \eta\Vert ^2, \>\>  \Vert A(P_s,s.P, s\si)v\Vert^2=y\Vert v\Vert^2,  $$
  o $x$ et $y$ sont les scalaires tels que:
 $$B(s.P, {\tilde P}_s, s\si)  B({\tilde P}_s, s.P,s\si)=xId, \>\> A(s.P, P_s, s\si)A(P_s, s.P, s\si)= yId.$$
 On a donc:
 $$\Vert C(s, P',P, \si) \eta\otimes v \Vert^2 = xy \Vert \eta\otimes v\Vert^2. $$
 On calcule $xy$ ˆ l'aide de (\ref{aaa}) et de  la dŽfinition des matrices $B$.  On trouve, par un calcul analogue ˆ celui de [W], fin de la preuve du Lemme V.2.2:
 $$xy= \mu^{M'}(s\si) \mu^G(s\si)^{-1} $$ 
 puis,  en utilisant (\ref{jw}):  $$xy= \mu^{M'}(s\si) \mu^{G}(\si)^{-1}.$$
 Le Lemme en rŽsulte. \qed
 Nous aurons besoin du Lemme suivant.
 \begin{lem}\label{holbor} Soit $E$ un espace vectoriel complexe de dimension finie.
 Soit $f$ une fonction holomorphe sur le complŽmentaire, $U$,  d'un nombre fini d'hyperplans  d'un voisinage de $0$ dans  $E$ et ˆ valeurs dans un espace vectoriel complexe de dimension finie. On la suppose  en outre bornŽe sur $U$.  On suppose qu'il existe une famille finie de  fonctions affines sur $E$, $(l_i)_{i\in I}$, telle que le produit de $f$ par le produit des $l_i$ s'Žtende en une fonction holomorphe, $g$, au voisinage de $0$ dans $\C^n$. Alors $f$ s'Žtend en une fonction holomorphe au voisinage de $0$ dans $E$. 
 \end{lem}
 \dem 
 On peut supposer que $I$ a ŽtŽ choisi minimal et il faut montrer qu'il est alors vide. Raisonnons par l'absurde et supposons $I$ non vide. Soit $i_0\in I$.  On note $f_0:= (\prod_{i\in I\setminus \{i_0\} } l_i ) f$. Alors $f_0$ est holomorphe sur $U$ et $l_{i_0}f_0 $ s'Žtend en une fonction holomorphe  au voisinage de $0$, ˆ savoir $g$.  D'aprs la minimalitŽ de $I$, $l_{i_0}$ s'annule en 0. Choisissons des coordonnŽes affines $(z_1, \dots, z_n)$ telles que  $l_{i_0}$ soit Žgale  ˆ la fonction $z_1$. On a:
 $$f_0(z_1, \dots  , z_n)= z_1^{-1} g(z_1,\dots, z_n).$$
 Comme $f_0$  est bornŽe au voisinage de 0, on en dŽduit que $g$ s'annule dans un voisinage de 0 sur l'hyperplan d'Žquation $z_1=0$ et  $f_0$ s'Žtend en une fonction holomorphe au voisinage de 0. Ceci contredit le fait que $I$ est minimal. Ceci achve de prouver que $I$ est vide.\qed 
 \begin{prop}\label{mumuC}Supposons l'assertion A vraie. On note encore $C(s,P',P,.)$ le prolongement rationnel de  $C(s,P',P,.)$ (cf. Proposition \ref{cteis}) ˆ $\O_\C$. 
 Avec les notations de la Proposition prŽcŽdente, l'application $\si\mapsto \mu^{G} (\si) \mu^{M'} (s\si)^{-1}C(s, P', P, \si)$ est holomorphe au voisinage de  $\O$.
 \end{prop}
 \dem Soit $(\si, E)$ objet de $ \O$ et $\phi\in \fhi$. Gr\^ace ˆ (\ref{poleA}), (\ref{aaa}) et (\ref{poleB}),  on peut appliquer le Lemme prŽcŽdent ˆ la fonction dŽfinie sur le complŽmentaire d'un nombre fini d'hyperplans affines de $(\a_M ')_\C$: $$\lambda \mapsto 
 \mu^{G} (\si\otimes \chi_{\lambda} ) \mu^{M'} (s(\si\otimes \chi_\lambda)) ^{-1}C(s, P', P, \si\otimes \chi_\lambda)\phi (\si\otimes \chi_\lambda)$$ dont les valeurs sont de norme constante,  d'aprs la Proposition \ref{MS}. La Proposition en rŽsulte. \qed 
  \begin{prop} \label{Chol} Soient $P=MU$, $P'=M'U'$ des sous-groupes paraboliques anti-standard de $G$ tels que $M$ et $M'$ soient conjuguŽs,  $s\in W(M'\vert G\vert M) $et soit $\O_\C$ l'orbite inertielle d'une reprŽsentation lisse irrŽductible de carrŽ intŽgrable de $M$. On pose, pour $(\si,E)$ objet de $\O_\C$ et lorsque cela ˆ un sens: 
 $$C^0(s, P', P,  \si ) : = B(s.P, P', s\si)^{-1} \otimes  A(P', s.P, s\si) \lambda (s). $$
  Tenant compte de l'ŽgalitŽ $Wh(  s. P, s\si) = Wh( \si) $ (cf. (\ref{whp} )), on voit que lorsqu'il est dŽfini,
   $C^0(s, P', P,  \si )$ est ŽlŽment de $Hom (Wh(\si) \otimes i^G_PE, Wh(s\si) \otimes i^G_{P'} sE)$.
   Cela dŽfinit une fonction  rationnelle sur $\O_\C$. 
   De plus  l'application $\si \mapsto C^0(s, P' , P, \si)$ est holomorphe au voisinage de $\O$ et,   pour tout $\si$ objet de $\O$, $C^0(s,P', P, \si)$ est unitaire.
 \end{prop}\dem D'aprs (\ref{hom}) et (\ref{B}),  on a:
 $$B(s.P, P', s\si)^{-1}=  j(s\si)^{-1} B(P',s.P,s\si).$$
   De (\ref{adjA}) et de  l'assertion A, on dŽduit que $C^0(s,P', P, \si)^* C^0(s,P', P,\si)$ est Žgal ˆ
  $$
 [B(s.P, P', s\si)B(P',s.P, s\si)]^{-1}\otimes [ \lambda(s^{-1} ) A(s.P, P', s\si)  A(P', s.P, \si) \lambda (s)] .$$
 Utilisant (\ref{aaa}) et son analogue  pour $B$, on en dŽduit.  
 $$C^0(s,P', P, \si)^* C^0(s,P', P, \si)= Id.$$ De faon similaire,  on voit que:
 $$C^0(s,P', P,\si)C^0(s,P', P, \si )^*= Id$$ ce qui prouve que,  lorsqu'il est dŽfini, $C^0(s,P', P)$ est unitaire.  Pour tout $\phi\in \fhi$,  d'aprs  (\ref{poleA}), (\ref{aaa}) et (\ref{poleB}),  on peut  appliquer le Lemme \ref{holbor} ˆ la fonction dŽfinie sur le complŽmentaire d'un nombre fini d'hyperplans affines de  $(a_M ')_\C$,  $\lambda \mapsto C^0(s,P', P,\si_\chi)\phi(\si_\chi)$. Le Lemme en rŽsulte. \qed 
 \subsection{Transformation unipotente de paquets d'ondes}
La preuve du ThŽorme suivant suit au plus prs la preuve de la Proposition VI.4.1 de [W]. Soit $P=MU$ un sous-groupe parabolique anti-standard de $G$ et soit  $\O$ l'orbite inertielle unitaire d'une reprŽsentation lisse irrŽductible et de carrŽ intŽgrable de $M$. 
Soit $\phi\in Pol(\O, Wh \otimes i^G_P) $. On dit que  que $\phi$ est trs rŽgulire si $\phi$ vŽrifie les propriŽtŽs suivantes:
\ber\label{trc} Pour tout sous-groupe parabolique standard, $P'=M'U'$, de $G$ et  $s\in W(M'\vert G \vert M) $, la fonction rationnelle sur  $\O_\C$, $\si\mapsto \mu^{G} (\si) \mu^{M'} (s\si)^{-1}C(s, P', P, \si)\phi(\si)$,  est  holomorphe sur un voisinage de $\O$. \eer
\ber \label{tgc0}Pour tout  sous-groupe parabolique antistandard $P=MU$, $P'=M'U'$, de $G$ et  $s\in W(M'\vert G \vert M) $, la fonction rationnelle sur $\O_\C$,  $C^0(s,P',P) \phi(\si)$  est  holomorphe sur un voisinage de $\O$.   \eer 
D'aprs les propriŽtŽs de rationalitŽ  (\ref{ratintertw}), (\ref{B}):
\ber\label{pO}
 Il existe $p_\O \in Po(\O)$ non identiquement nul tel que pour tout $\phi\in  Pol(\O, Wh \otimes i^G_P)$,  $p_\O\phi$ soit trs rŽgulire. 
\eer 
On note que si $g\in G$ et si $\phi$ est trs rŽgulire, $\rho_\bullet(g)\phi$ est trs rŽgulire d'aprs les propritŽs d'entrelacement des fonctions $C$ et $C^0$. 
On remarque que, d'aprs les Propositions  \ref{mumuC} et  et \ref{Chol}:
\ber \label{Atr} Si l'Assertion A est vraie, tout $\phi\in Pol(\O, Wh \otimes i^G_P) $ est trs rŽgulire. \eer 
  \begin{theo} \label{fphipind} Soit $P=MU$,   $P'=M'U'$   deux sous-groupes paraboliques anti-standard de $G$, soit $\O$ l'orbite inertielle d'une  reprŽsentation, lisse irrŽductible  et de carrŽ intŽgrable de $M$, et $\phi\in Pol(\O, Wh \otimes i^G_P) $.   Pour $\phi$   trs rŽgulire, on a: \\ 
(i) Si le rang semi-simple de $M'$ est infŽrieur ou Žgal ˆ celui de $M$  et si
$M'$ n'est pas conjuguŽ ˆ $M$, $f_\phi^{P', ind}$ est nul.  
 \\ (ii)  Supposons $M$ et $M'$ conjuguŽs. 
Pour tout $m'\in M'$ et $g\in G$,  on a: $$(f_\phi^{P',ind}(g)) (m') = \sum_{s\in W(M'\vert G\vert M)} 
\int _{\O _u}    E^{M'}_{M'}  [(C^{0}  (s, P', P,\si) (\phi(\si) ) (g)](m')d\si. $$
\\ (iii) De plus $f^{P, ind} $ est ˆ valeurs dans $\CCU$. 
\\(iv) Si l'assertion A est vraie, les mmes conclusions sont valables pour tout $\phi \in C^{\infty} (\O, Wh \otimes i^G_P) $. 
 \end{theo}
\dem
Prouvons (i) et (ii). 
   Soit $\phi\in Pol(\O, Wh \otimes i^G_P) $, $g\in G$ et $m'\in M'$. Ona:
 $${\rho(g)}f_\phi= f_{\rho_\bullet (g)  \phi },$$ 
$$f_\phi^{P'}(m')= (\rho(m') f_\phi^{P'})(1)= \delta_{P'}^{-1/2} (m') (\rho(m') f_\phi)^{P'}(1)$$
et pour $ v' $ ŽlŽment de l'espace de  $i^G_{P'} s\si$: 
$$((i^G_{P'} s\si )(m') )v' ) (1)=   \delta_{P'}^{1/2}(m') (s\si(m') )(v'(1)).$$  On est  donc  ramenŽ ˆ calculer $f_\phi ^{P'}(1)$ pour tout $\phi$ trs rŽgulire.
 On fixe un sous-groupe compact ouvert
$H$  comme dans (\ref{iwa}) tel que 
 $\phi$ est invariante par $\rho_\bullet (H)$. 
\ber 
Il existe une constante $C_H>0$ telle que pour tout $\varphi\in
C^{\infty}(H)$: $$\int _H
\varphi(g)dh=c_H \int _{H\cap U'^-\times H\cap M' \times  H\cap U'} \varphi(u'^-m'_1 u'
)du'^- dm' du', $$ 
o les mesures sont celles induites par les mesures de Haar sur $G, U'^-, M',U'$.
\eer 
On fixe $a \in A_{M'}$ tel que $\vert \alpha (a)\vert _F <1$ pour tout $\alpha
\in \Sigma(P')$. Pour $n \in \N$, posons $ U'_n= a^{-n} (H\cap U')a^n$.  Alors on a:
 $$f_\phi^{P'} (1)=\int_{U'} f_\phi(u') du'.$$
 Donc \beq \label{lim} f_\phi^{P'} (1) =  lim_{n \to \infty} \int_{U'_n}f_\phi (u') du' = lim_{n \to \infty} \delta_{P'}(a)^{-n}
\int _{H\cap U' } f_\phi (a^{-n} u'a^{n})du'.  \eeq 
 Pour tout $n \in \N$:  
 $$\int _{H\cap U' } f_\phi (a^{-n} u'a^{n}) du'  = \int _{H\cap U'}\int _{\O} 
E^G_P(\phi (\si) ) (a^{-n} u' a^n ) d\si du'.  $$ 
 On remarque que: 
 $$E^G_P(\phi (\si)) (a^{-n} u' a^n ) = E^G_P((\rho_\bullet (u'a^{n}) \phi  )
(\si) ) (a^{-n}).$$ 
 On pose alors: 
 $$\phi_n= \int _{ H \cap U'} \rho_\bullet (u'a^{n}) \phi du' \in Pol(\O, 
 Wh \otimes  i_P^G),$$ l'intŽgrale se rŽduisant ˆ une somme finie, puisque $\phi$ est $H$-invariante.
 Alors:
  \beq \label{fp'}  f_\phi^{P'} (1) =lim_{n \to \infty} \delta_{P'}(a)^{-n})
\int_{\O} E^G_P(
\phi_n (\si)) (a^{-n}) d \si.\eeq
   Comme $a^{-n} (H \cap M') (H\cap U'^- )a^n \subset H$, on a l'ŽgalitŽ:  
  \beq \label{phin} \phi_n = c  \int _{H\cap U' \times  H\cap M' \times H\cap  U'^-} \rho_\bullet  (
u' m' u'^-  a^n ) \phi du'^- dm' du' \eeq 
 $$ =  cc_H^{-1} \int _H \rho_\bullet (ha^n) \phi dh. $$
 o $c= vol (H\cap M')^{-1} vol (H \cap U'^-) ^{-1} $. 
 Donc $\phi_n \in Pol (\O,   Wh \otimes i_P^G) ^H$. \\
 D'aprs les propriŽtŽs du terme constant (cf.  (\ref{ffp})),  et en tenant compte de l'ŽgalitŽ $\delta_{P'^-}= \delta_{P'}^{-1}$,  on
peut choisir
$N$ assez grand pour que,  pour tout
$n\geq N$, on ait:
 $$ \label{ep'}E^G_P(\phi_n (\si ) ) (a^{-n}) = \delta_{P'} (a)^{n/2}  E^G_P(
\phi_n(\si) ) _{P'^-} (a^{-n}) .$$ 
et avec la notation (\ref{f+}):
\beq  \label{ep'} E^G_P(\phi_n (\si ) ) (a^{-n}) = \delta_{P'} (a)^{n/2} [ E^G_P(
\phi_n(\si) )^w _{P'^-} (a^{-n}) +E^G_P(
\phi_n(\si) )^+ _{P'^-} (a^{-n})]. \eeq
En reportant dans les ŽgalitŽs (\ref{fp'}) puis (\ref{lim}), on obtient: 
\beq f_\phi^{P'}(1)= X^w+ X^+, \eeq
o
\ber \label{Xw+}$$X^w = lim_{n\to \infty}  \delta_{P'} (a)^{n/2}  \int_{\O} \mu(\si) [ E^G_P(
\phi_n(\si) )^w _{P'^-} (a^{-n})d\si,$$
$$X^+ = lim_{n\to \infty}  \delta_{P'} (a)^{n/2}  \int_{\O} \mu(\si)  [ E^G_P(
\phi_n(\si) )^+ _{P'^-} (a^{-n})d\si,$$\eer 
pour peu que ces limites existent. \\
Montrons: 
\ber \label{ineg} Il existe $\varepsilon>0$, $C>0$ et $d\in \N$ tels que pour tout $n\in \N$ et tout $\si\in \O$ on ait:
$$  [ E^G_P(
\phi_n(\si) )^+ _{P'^-} (a^{-n})\leq C \delta_{P'} (a)^{n/2} (1+ \Vert H_0(a^n)\Vert )^d e^{-\varepsilonÊ\Vert H_0 (a^n)\Vert}. $$
\eer 
Par compacitŽ de $\O$, il suffit de prouver que pour $\si_0 $ objet de $\O$, on a une telle majoration pour $\si$  dans un voisinage convenable de $\si_0$. Fixons $\si_0\in \O$ et un ensemble fini ${\cal F}\subset C^{\infty} (\O, Wh\otimes i^G_P)^H$    tel que $\{\phi'(\si_0) \vert \phi' \in {\cal F}\}$ forme une base de $Wh(\si_0) \otimes (i^G_P E_{\si_0})^H$. Il existe un voisinage  $\cal U$ de $\si_0$ dans $\O$ tel que pour tout $\si \in {\cal U}$,  $\{\phi'(\si) \vert \phi' \in {\cal F}\}$ forme une base de $Wh(\si) \otimes (i^G_P E_{\si})^H$. Pour tout  $\phi''\in C^{\infty} (\O, Wh\otimes i^G_P)^H$ et tout $\si \in {\cal U}$, on peut donc Žcrire de faon unique:
$$\phi'' (\si)= \sum_{\phi'\in {\cal F} } C(\phi', \phi'', \si)\phi'(\si).$$
Quitte ˆ restreindre $\cal U$, on peut supposer qu'il existe $C_1>0$ tel que pour tout $\phi'' \in C^{\infty} (\O, Wh\otimes i^G_P)^H$ et tout $\si \in {\cal U}$, on ait: 
$$\sum_{\phi' \in {\cal F} }\vert C( \phi',\phi'', \si)\vert \leq C_1 Sup\{ ( \phi'(\si), \phi''(\si))\vert \phi'\in {\cal F}\} 
.$$
D'aprs le Lemme \ref{E+a}, il existe $C_2>0$ et  $\varepsilon >0$ tels que pour tout $\si$ objet de $\O$ et $ \phi'\in {\cal F}$:
$$ \vert (E^G_P (\phi' (\si)) )^+_{P'} (a^{-n})\vert  \leq C_2e^{-\varepsilon n\Vert H_0(a)\Vert }.$$
Donc \beq \label{egpa}\vert (E^G_P \phi_n ( \si) )^+_{PÔ} (a^{-n})\vert  \leq C_1 C_2Sup\{ ( \phi'(\si), \phi(\si))\vert \phi'\in {\cal F}\} e^{-\varepsilon n\Vert H_0(a)\Vert }.\eeq
Fixons $\phi'\in {\cal F}$. Pour $(\si,E) \in {\cal U}$, on a:
$$( \phi'(\si), \phi_n(\si)) = \int_{H\cap U} ( \phi'(\si),\rho (ua^n)\phi(\si)) du. $$
Comme $\phi'$ est invariante ˆ droite par $H$ et le produit scalaire sur $Wh( \si ) \otimes  i^G_PE$ Žtant $G$-invariant, on a donc :
$$( \phi'(\si), \phi_n(\si)) = vol( H\cap U) ( \phi'(\si),\rho (a^n)\phi(\si)).   $$
De la dŽfinition de ce  produit scalaire rŽsulte l'existence de  $\phi'_1 \in C^{\infty} (\O, \check{i^G_P } \otimes i^G_P)$ telle que, pour tout $\si$ objet de $\O$, l'application $g\mapsto  ( \phi'(\si),\rho (g)\phi(\si)) $ soit de la forme $g \mapsto E^G_P ( \phi_1'(\si))(g)$, o $E^G_P$ est l'application coefficient de [W] section V.1.
Donc, d'aprs [W], Lemmes VI.2.2, II.1.1 et  Equation I. 1(6), il existe $C_3>0$ tel que pour tout $\si\in \O$:
$$\vert (\phi'(\si), \phi_n(\si)\vert \leq C_3 \delta_0(a^{n/2})( 1 + n\Vert H_0(a)\Vert)^d.$$
Joint ˆ (\ref{egpa}), cela prouve (\ref{ineg}).
Joint ˆ (\ref{Xw+}) et au choix  de  $a$, on en dŽduit que $X^+=0$.
\\ En utilisant successivement  (\ref{covct}),  l'ŽgalitŽ  $\delta_{P'^-}=
\delta_{P'}^{-1}$ et  la dŽfinition de
$U'_n$ pour effectuer un changement de variable, on Žtablit les ŽgalitŽs:  
$$ E^G_P(\phi_n (\si) )_{P'^-} (a^{-n}) = \delta_{P'}^{1/2} (a^{-n} )   \int _{H\cap
U'}  [E^G_P(\phi(\si)) _{P'^-} ^{ind} (a^{-n} u'a^n )] (1) du'   $$
 $$ =  \delta^{1/2}_{P'} (a^{n} ) \int _{U'_n}[E^G_P(
\phi(\si)) _{P'^-} ^{ind} (u')] (1)du' . $$
 En tenant compte de (\ref{fp'}) et (\ref{ep'}), on en dŽduit:

$$X^w =lim_{n\to \infty} \int_{\O}  \int _{U'_n}\mu(\si) [E^G_P(
\phi(\si)) _{P'^-} ^{ind} (u')] (1)du' d\si. $$
Si le rang semi-simple de $M'$ est infŽrieur ou Žgal ˆ celui de $M'$ et si $M$ et $M'$ ne sont pas conjuguŽs, on voit que $X^w=0$. Ceci achve de prouver (i).\\
Prouvons (ii).
\\ Suposons $M$ et $M'$ conjuguŽs. 
On va utiliser la Proposition \ref{cteis} pour donner une expression de $E^G_P(
\phi(\si)) _{P'^-} ^{ind}$. Comme $M$ et $M'$ sont 
conjuguŽs, pour tout $s \in W(M'\vert G\vert M) $, $M'\cap s.P= M'$. Dans notre utilisation du ThŽorme \ref{cteis},  $P'$ est  ici remplacŽ par $P'^-$ et 
$P_s$ est ici Žgal ˆ 
$P'^-$, 
${\tilde P}_s$ est ici Žgal ˆ $ P'$. On a alors, gr\^ace ˆ la Proposition \ref{cteis} et (\ref{epindepbis}):
$$X^w= lim_{n\to \infty} X^w_n$$ o $$X^w_n=  \sum_{s\in W(M'\vert G\vert M)}\int _{\O}   \int _{U'_n} \mu(\si)
[E^{M'}_{M'}(C(s, P'^-, P, \si) \phi(\si))(u')](1) d\si du' $$
 et on veut passer ˆ la limite sur $n$. Mais un dŽplacement de contour d'intŽgration est nŽcessaire. Si $\chi$ est un caractre non ramifiŽ de
$M$, on note $\O\chi$ l'ensemble des classes d'Žquivalence des reprŽsentations
$\si _ \chi$  lorsque $\si$ dŽcrit les objets de $\O$.  On munit $\O\chi$
de la mesure obtenue par transport de structure de la mesure  sur $\O$. Posons:   \beq\label{phisc} \phi_s(\si):=\mu(\si)  C(s, P'^-, P, \si) \phi(\si) \in Wh(P', s\si)\otimes  i^G_{P'^-}(sE). 
 \eeq  Comme $\phi
$ est trs rŽgulire, il existe $ \varepsilon_1>0$ tel que, pour tout $s\in W(M'\vert G\vert M)$ $\phi_s$,  est holomorphe  sur $\{\si \otimes \chi\vert \Vert Re \> \chi  \Vert < \varepsilon_1 \Vert Re \>\delta_P \Vert \}$. On choisit $\varepsilon>0$ strictement plus petit que $\varepsilon_1$. On pose \beq \label{lse} \Lambda_{s, \varepsilon}:= s^{-1} \delta_{P'}^{-\varepsilon}. \eeq   Pour des
raisons  d'holomorphie,, on a: 
 \beq \label{fpoun} X^w_n=  \sum_{s\in W(M'\vert G\vert M)}\int _{\O
\Lambda_{s, \varepsilon}}  
\int _{U'_n} [ E^{M'}_{M'}(\phi_s(\si)(u') )](1) d\si du'.\eeq
 On veut passer ˆ la limite sur $n$ dans cette expression.   
 Il faut  majorer $[ E^{M'}_{M'}(\phi_s(\si) )(u') ](1) $. \\
 Soit $(\si, E)$ un objet de $\O$. 
 Soit $v_s \in i^G_{P'^-} (sE)_{s \Lambda_{s, \varepsilon }}$,     $\eta_s\in Wh (P', s\si)$.  On remarque que  $Wh(P', s\si)$ est Žgal ˆ $Wh(s \si)$, puisque $P'$ est anti-standard.  
Alors, il rŽsulte de la dŽfinition que:   
$$[E^{M'}_{M'}(v_s\otimes
\eta_s )(g) ](1)=\langle  \eta_s, v_s    (g)\rangle  .$$ On Žcrit:   $$u'= u^{'-} (u) m'
(u') k(u'), \> u'^-(u')\in U'^-, \> m'(u')\in M', k(u') \in K,$$
de sorte que $m'(u')=m_{P'^-} (u')$, etc.. 
 Tenant compte des propriŽtŽs de covariance de $v_s$, on voit que: 
 $$\langle  \eta_s,v_s    (u')\rangle    = (s\Lambda_{s, \varepsilon })(m'(u')) \delta_{P'}^{-1/2}(m'(u'))\langle  \eta_s, s\si (m'(u')) 
v_s (k(u')) \rangle       $$ 
 et,  tenant compte de la dŽfinition de  $\Lambda_{s, \varepsilon }$, on obtient: 
$$\langle \eta_s, v_s(u)\rangle=  \delta_{P'}^{-1/2- \varepsilon}(m'(u'))\langle  \eta_s, s\si (m'(u'))  v_s
(k(u')) \rangle  . $$
Tenant compte du fait que la restriction de $v_s$ ˆ $K$ ne prend qu'un nombre fini
de valeurs, on dŽduit de (\ref{Csec}), du Lemme \ref{hsi} (i) et  de la tempŽrance de $\eta_s$ (cf. Proposition \ref{squaresquare} ), que pour tout $d\in N$, il existe $C>0$ tel que :
$$ \vert \langle \eta_s, v_s(u)\rangle \vert \leq C  \delta_{P'}^{-1/2}(m'(u))\delta_{P_0\cap M} ^{1/2} ( m_0(m'(u'))(1+ \Vert H_0(m_0(m'(u'))\Vert )^{-d}$$
 ou $C$ est une constante qui ne dŽpend que de la restriction de $v_s$ ˆ $K$ et de
$\eta_s$ mais pas de $\si\in \O$.  Comme $P'$ est anti-standard, on a l'ŽgalitŽ:
$$ \delta_{P'} ^{-1} (m'(u'))\delta_{P_0\cap M}  ( m_0(m'(u'))= \delta_{P_0} (m_0( m'(u')). $$ Donc, on a:
$$ \vert \langle \eta_s, v_s(u)\rangle \vert \leq C  \delta_{P_0}^{1/2}  (m_0(m'(u'))) (1+ \Vert H_0(m_0(m'(u'))\Vert )^{-d}.$$
\\ On en dŽduit que , pour tout $d\in \N$, il existe $C'>0$ tel que  pour tout $s \in W(M'\vert G\vert M)$ et $\si$ objet de $ \O
\LL_{s, \varepsilon}$: 
$$\vert E^{M'}_{M'}(\phi_s(\si))(u') ](1)\vert \leq C' \delta_{P_0} (m'(u')) (1+ \Vert H_0(m_0(m'(u'))\Vert )^{-d}  , u'\in U'.$$
Comme le second membre de cette inŽgalitŽ est une fonction intŽgrable sur $\O
\LL_{s, \varepsilon}  \times U'$, pour $d $ bien choisi, d'aprs le Lemme \ref{hsi}. On peut appliquer le ThŽorme de convergence dominŽe et 
dŽduire de  (\ref{fpoun}): 
$$X^w=  \sum_{s\in W(M'\vert G\vert M)}\int _{\O \Lambda_{s, \varepsilon } }  \int
_{U'}  \mu(\si) [ E^{M'}_{M'}(\phi_s(\si)(u')) ](1) d\si du',$$
la fonction sous le signe intŽgrale Žtant intŽgrable pour la mesure produit. On
peut appliquer le thŽorme de Fubini et commencer par calculer l'intŽgrale sur $U'$. Soit $\si$ objet de $\O
\LL_{s, \varepsilon}$.  
 Si $v $ est invariante  par un sous-groupe compact ouvert de $G$, $H$, on a: 
 $$\langle \eta, (A(P', P'^-, s \si) v)(1) \rangle = \langle e_{M'\cap H} \eta, (A(P', P'^-, s \si) v)(1) \rangle  .    $$
 Donc,  d'aprs la dŽfinition des intŽgrales d'entrelacement, notamment la dernire assertion de (\ref{intertw}), on a:
 $$\langle \eta, (A(P', P'^-, s \si) v)(1)> = \int_{U'} \langle e_{M'\cap H} \eta, v(u')>du'.$$ 
Comme $v(u')$ est $M'\cap H$-invariante,  on obtient finalement:
$$\langle \eta, (A(P', P'^-, s \si) v)(1)> = \int_{U'} \langle \eta, v(u') \rangle du'. $$ De cette dernire ŽgalitŽ et de (\ref{phisc}),   on en dŽduit:  
 $$\int_{U'}  [ E^{M'}_{M'}(\phi_s(\si))(u') ](1)du'=  
 E^{M'}_{M'}[((Id_{Wh(P', s\si) } \otimes A(P', P'^-, s\si))C(s, P'^-, P, \si)
\phi(\si))(1)](1).$$  
  Tenant compte de la dŽfinition des fonctions $C$  (cf. Proposition \ref{cteis}) et de (\ref{aaa}), on voit que  
 $f_\phi^{P'}(1)$ est Žgal ˆ la somme sur ${s \in W(M'\vert G\vert M)}$ de: $$ \int _{O \Lambda_{s, \varepsilon }}\mu(\si)     j (P', P'^-, s.P, s\si)  E^{M'}_{M'} 
 [ ((B(P', s.P, s\si, s\chi) \otimes  (A(P', s.P,
s\si)\l(s))
\phi(\si) )(1)](1) d\si .$$ 
Montrons, avec les notations de la Proposition \ref{Chol}:   \beq \label{mujc}\mu(\si)     j (P', P'^-, s.P, s\si)   
 [ ((B(P', s.P, s\si, s\chi) \otimes  (A(P', s.P,
s\si)\l(s))= C^0(s, P', P,  \si ) .\eeq 
Comme $C^0(s, P', P,  \si ) = B(s.P, P's\si)^{-1} \otimes  A(P', s.P, s\si) \lambda (s)$, il faut donc  montrer $$\mu(\si)     j (P', P'^-, s.P, s\si)   
 B(P', s.P, s\si, s\chi)= B(s.P, P', s\si)^{-1}.$$
 Mais $$B(P', s.P, s\si) B(s.P, P', s\si) = j(P', s.P, P', s\si). $$
 Donc $$\mu(\si)     j (P', P'^-, s.P, s\si)   
 B(P', s.P, s\si, s\chi)$$ 
 est Žgal ˆ: 
 $$    \mu(\si)     j (P', P'^-, s.P, s\si)  j(P', s.P, P') B(s.P, P', s\si)^{-1}. $$
 On dŽduit de (\ref{aaa}),  $  j (P', P'^-, s.P, s\si)  j(P', s.P, P')= j(s\si)$.
 Mais, cf. [W], IV. 3(3)  $\mu(\si)= \mu(s\si)= j(s\si)^{-1}$.
 Cela prouve (\ref{mujc}). \\
 En choissant $\varepsilon$ suffisamment petit, comme $\phi$ est trs rŽgulire,   on peut, pour des raisons d'holomorphie,  remplacer l'intŽgrale sur $\O \LL_{s, \varepsilon}$ par l'intŽgrale sur $\O$. Ceci achve la preuve de  (ii).   \\Alors (iii) rŽsulte de (i) et (ii) et de la Proposition \ref{waveS} appliquŽe ˆ $M$ au lieu de $G$ et $P$.\\ (iv) Si l'Assertion A est vraie,  le rŽsultat est aussi valable pour   $\phi$ polynomiale, d'aprs (\ref{Atr}). 
 D'aprs les Propositions 
\ref{waveS} et \ref{fpint} (ii), l'application ,  $\phi\to f_\phi^P (1)$ est une forme linŽaire  continue sur $C^\infty(\O, Wh \otimes i^G_P)$.  Son Žvaluation en $\phi$ est Žgale au premier membre de l'ŽgalitŽ de (ii)  pour $g=m'=1$.  Mais, si l'Assertion A est vraie,  le deuxime membre, pour $g=m'=1$, est Žgalement une forme linŽaire continue sur ce m\^eme espace membre est Žgalement  d'aprs la Proposition \ref{waveS}. La densitŽ de $Pol(\O, Wh \otimes i^G_P)$ dans $C^\infty(\O, Wh \otimes i^G_P)$ et (i) et (ii) impliquent l'analogue de (i) et (ii) pour $\phi \in \fhi$. 
\\ Utilisant la Proposition \ref{waveS} appliquŽe ˆ $M$ au lieu de $G$ et $P$ , on dŽduit l'analogue  de (iii).
\qed 

\subsection{ TransformŽe de Fourier  de paquets d'ondes}
 \begin{theo} \label{fphihat} Soient $P=MU, P=M'U'$  des sous-groupes paraboliques anti-standard  de $G$ . Soit  $\O$ l'orbite inertielle unitaire  d'une reprŽsentation lisse irrŽductible   et de carrŽ intŽgrable de $M$.  Soit $\phi \in Pol (\O, Wh \otimes i^G_P)$  trs rŽgulire ou, si l'Asserion A est vraie,  $\phi \in C^\infty (\O, Wh \otimes i^G_P)$  . Soit  $(\si_1, E_1)$ une reprŽsentation lisse irrŽductible   et de carrŽ intŽgrable de  $M'$. Si $M$ et $M'$ sont  conjuguŽs, on a:
 $$\hat{f}_ \phi(g) (P', \si_1)= \sum_{ s \in W(M'\vert G\vert M), s^{-1} \si_1\in
{\O}}C^0 (P', P, s, s^{-1}\si_1) \phi (s^{-1} \si_1).$$
Si l'Assertion A est vraie, le m\^eme rŽsultat est vrai pour $\phi \in C^\infty(\O, Wh \otimes i^G_P)$.
 \end{theo}
  \dem La preuve est entirement analogue ˆ celle du Thorme 5 de [D4] , en y remplaant $\tilde (C(s,P',P,s\si )\phi)(s\si)$ par $C^0(s,P',P, \si) \phi(\si)$.   
\subsection{Produit scalaire de paquets d'ondes}
\begin{prop} 
 \label{prodscalpaq}
Soit $P=MU, P'=MU'$ deux sous-groupes paraboliques anti-standard de $G$. Soit $\O$ (resp. $\O_1$) l'orbite inertielle d'une reprŽsentation lisse, irrŽductible et  cuspidale de $M$ (resp. $M'$). Soit $\phi \in Pol (\O, Wh i^G_P )$, $\phi_1 \in Pol(\O_{1}, Wh \otimes i^G_{P'})$ trs rŽgulires (resp. $\phi \in C^\infty (\O,  i^G_P), \phi_1\in C^\infty (\O_{1}, Wh \otimes i^G_{P'})$ si l'Assertion A est vraie ).  \\
(i) Si $M$ et $M'$  ne sont pas conjuguŽs
dans $G$, $(f_\phi, f_{\phi_1})_G$ est nul.   \ste 
(ii) Si $M$ et $ M'$  sont conjuguŽs dans $G$, $(f_\phi, f_{\phi_1})_G$ est
Žgal ˆ:
 $$ \int_ {\O_{1} }\sum_{ s \in W(M'\vert G\vert M),  s{\cal O}= {\cal
O}_1 } \mu(\si_1)(C^0 (s, P', P,s^{-1}  \si_1) \phi(s^{-1}\si_1), \phi_1(\si_1))d\si_1.$$
\end{prop}
 \dem La  dŽmonstration est semblable ˆ celle de [W], Proposition VII.2.2. Nous la donnons pour la commoditŽ du lecteur.\ste 
On a l'ŽgalitŽ $$(f_\phi, f_{\phi_1})_G=  \int_{U_0\backslash G}
 f_\phi(g) \int_{{\O_1}} \overline{\mu( \si_1)(E^G_{P'}(\phi_1(\si_1))(g)}d\si_1
dg. $$
 Comme ${\O_1}$ est compact et  comme $f_\phi\in \CCU$, on dŽduit du Lemme \ref{deriv}  et (\ref{intff})  que 
 l'intŽgrale double est absolument convergente et l'on obtient:   
 $$(f_\phi, f_{\phi_1})_G=  
  \int_{{\O_1}} \mu(\si_1) (f_\phi, E^G_{P'}(\phi_1(\si_1)))_Gd\si_1 $$ 
et d'aprs la dŽfinition de $\hat{f}$:
$$(f_\phi, f_{\phi_1})_G=  
  \int_{{\O_1}}\mu(\si_1) (\hat{f_\phi}(P', \si_1), \phi_1(\si_1))d\si_1 . $$
Supposons le rang semi-simple de $M'$ infŽrieur ou Žgal ˆ celui de $M$.
En utilisant la Proposition  \ref{findfhat} et le ThŽorme \ref{fphipind} (i),
on voit que  $(f_\phi, f_{\phi_1})_G=0$ si $M$ n'est pas conjuguŽ ˆ $M'$.
Si le rang semi-simple de $M'$ est strictement plus grand que celui de
$M$, il suffit d'appliquer la relation $(f_\phi, f_{\phi_1})_G=
\overline{(f_{\phi_1}, f_\phi)_G}$, pour achever la preuve de (i).
\ste Supposons maintenant $M$ et $M'$ conjuguŽs dans $G$.
Le ThŽorme  \ref{fphihat} calcule  $\hat{f_\phi}(P', \si_1)$, ce qui
conduit ˆ (ii).\qed
 \begin{cor}
Avec les notations de la Proposition prŽcŽdente, en supposant l'assertion A vraie, si $M$ et $M'$ ne sont pas conjuguŽs, $\hat{f}(P',)$ est nul.\end{cor}
\dem Pour tout $\phi_1$ comme ci-dessus , on a $(f_\phi, f_{\phi_1})_G=0$. Soit encore, d'aprs la Proposition \ref{Ehatf}:
$$\int_{\O_1 } (\hat{f}(P', \si_1), \phi_1(\si_1) \mu(\si_1)d\si_1=0 .$$
Comme cela est vrai pour tout $\phi_1\in C^\infty (\O_{1}, Wh \otimes i^G_{P'})$, on en dŽduit la nullitŽ voulue. \qed
\subsection{ Adjoint de la matrice: preuve de l'assertion A}
\begin{theo} \label{adjB}
Soit $P$, $Q$ des sous-groupes paraboliques  semi-standard de $G$ possŽdant le mme
sous-groupe de LŽvi semi-standard,  $M$. On suppose $P$ anti-standard. Soit $\O$ l'orbite
inertielle d'une reprŽsentation lisse, irrŽductible et  de carrŽ intŽgrable de $M$. On a l'ŽgalitŽ de
fonctions rationnelles sur $\O$:
$$B(P, Q, \si)^*= B(Q, P, \si)$$
\end{theo}
\dem
Soit  $P'$ le sous-groupe parabolique anti-standard de $G$ auquel $Q$
est conjuguŽ et $\O_1$ une orbite inertielle de $M'$ conjuguŽe de $\O$ par un ŽlŽment de $W(M'\vert G\vert M)$. 
Soit $\phi\in Pol (\O, Wh \otimes  i^G_P)$, 
$\phi_1 \in Pol (\O_{1}, Wh \otimes i^G_{P'}) $ trs rŽgulires.   Alors,  d'aprs la Proposition prŽcŽdente: 
$$(f_\phi, f_{\phi'})_G =     \int_ {{\cal O}_{1}}
 \sum_{ s \in
W(M'\vert G\vert M), s{\cal O}= {\cal O}_1 } \mu(\si_1)(C^0 (s, P', P, s^{-1}\si_1) \phi (s^{-1} \si_1), \phi_1(\si_1))d\si_1.
$$
Puis en utilisant $(f_\phi, f_{\phi'})_G = \overline{(f_{\phi'}, f_\phi)_G}$, on a: 
\beq \label{ff} (f_\phi, f_{\phi'})_G=\int_ {{\cal O}} \sum_{ t \in W(M\vert G\vert M'),  t^{-1}  \O=  \O_1}\mu(\si) (\phi(\si), C^0 (t, P, P',t^{-1} \si) \phi_1 (t^{-1}\si)) d \si. \eeq 

On pose 
$\si= s^{-1}_1 \si$ dans (\ref{ff}). Comme $ \mu(\si)= \mu(\si_1)$ d'aprs [W]  IV.3 (3), on a:
\beq \label{barff} (f_\phi, f_{\phi'})_G= \int_ {{\cal O}} 
\sum_{ s \in W(M'\vert G\vert M),  s{\cal O}= {\cal O}_1 } \mu(\si)(C^0 (s, P', P, \si) \phi (\si), \phi_1(s\si))d\si.
\eeq
  A $s \in  W(M'\vert G\vert M) $ correspond un unique ŽlŽment $t$ de $W(M\vert G\vert M')$ tel que $ts=m\in
M_0 \cap K$.  
Si $\phi$ est trs rŽgulire, il en va de mme de
$p\phi$ pour tout $p\in Pol(\O)$. Par ailleurs si $F\in
C^\infty(\O)$ est tel que:  $$\int_{\O}p(\si) F(\si) d\si=0,
p\in Pol(\O),  $$ alors  $F=0$.\\
 Donc, pour tout $\si$ objet de $ \O$, les
expressions sous le signe intŽgrale dans les membres de droite des
ŽgalitŽs (\ref{ff}) et (\ref{barff}) sont Žgales. \\ Soit $\O'$ l'ensemble
des 
$\siÊ\in
\O$ tel que, avec les notations de (\ref{pO}), 
$p_{\O}(\si)$ soit non nul et   tels que si $w\in W(M, \O)$, $w\si$ ne
soit Žquivalente ˆ
$\si$ que si $w=1$. On remarque que $\O '$ est dense dans $\O$. \\Soit $(\si,E) $ un objet de $\O'$, 
$ s\in W(M'\vert G\vert M)$ et $t$ comme ci-dessus. D'aprs (\ref{pO}) on
peut choisir
$\phi_1$ trs rŽgulire tel que $\phi_1(s \si )$ soit non nul et arbitraire dans
$i^G_P(tE_1)$ et tel que
$\phi_1(s' \si)= 0$ si $ s'\in W(M'\vert G\vert M)$ est distinct de $s$. De mme $\phi(\si)$ peut tre choisi arbitrairement.  Alors, tous les termes sous le signe intŽgral de (\ref{ff}) (resp.(\ref{barff})) sont nuls exceptŽ celui correspondant ˆ $s$ (resp. $t$).  On en dŽduit: 
\beq \label{Cphiphi}(C^0 (s, P', P, \si) \phi (\si), \phi_1(s\si))= (\phi(\si), C^0 (t, P, P', t^{-1} \si) \phi_1 (t^{-1} \si)) . \eeq 
Les deux membres de cette ŽgalitŽ sont des fonctions scalaires sur $\O$.
Comme $m\si $ est Žquivalente ˆ $\si$ , on ne change pas le deuxime membre en remplacant $\si$ par $m\si$. Comme on a   $s=t^{-1} m \si$, on a donc :

$$( \phi(\si), C^0 (t, P, P', t^{-1} \si) \phi_1 (t^{-1} \si))= (\phi(m\si), C^0 (t, P, P', s \si) \phi_1 (s \si)). $$Mais $\phi$ est une fonction sur $\O$ ˆ valeurs dans $Wh \otimes i^G_{P}$. L'opŽrateur $\si(m)$ entrelace $\si_1$ et $m\si_1$ et le transposŽ de $\si(m)$, $\si'(m^{-1})$, induit une bijection de $Wh( \si)$ sur $Wh(m\si ) $, 
Donc $$ \phi_1(m\si)=T \phi(\si),$$
o $T=  \si'( m^{-1} )\otimes \l(m)$, dont l'adjoint est Žgal ˆ $ \si'(m) \otimes \l(m^{-1})$.  Donc:
$$ (\phi(\si), C^0 (t, P, P', t^{-1} \si) \phi_1 ( t^{-1}\si)) =(   \phi(\si),(\si'(m)  \otimes \l(m^{-1} )) C^0 (t, P, P', s \si) \phi_1( s\si)),$$
o $*$ dŽsigne l'adjoint.
Joint ˆ l'ŽgalitŽ  (\ref{Cphiphi}), on en  dŽduit: 
\beq\label{C00} C^0 (s, P', P, \si)^*=  (\si'(m)  \otimes \l(m^{-1} ))C^0 (t, P, P', s \si). \eeq  
D'aprs la dŽfinition des fonctions $C^0$ (cf. Proposition \ref{Chol}) et l'ŽgalitŽ $ts=m$, on a: 
$$C^0(s, P', P, \si)= B(s.P, P', s\si )^{-1} \otimes A(P', s.P, s\si) \lambda (s) $$ $$ C^0(t. P, P', s\si ) = B(t. P', P, m\si )^{-1} \otimes A(P,t. P', m\si ) \lambda(t).$$
On a, par conjugaison par $m$ et en tenant compte de l'ŽgalitŽ $m^{-1}t=s^{-1}$ 
$$\l(m^{-1})  A(P,t. P', m\si ) \lambda(t)= A(P, t.P', \si) \lambda(s^{-1}).$$
Par ailleurs (cf. (\ref{adjA}): $$((A(P',s.P, s\si) \l(s))^*= \lambda (s^{-1}) A(s.P, P', s\si)$$
puis par conjugaison par $s^{-1}$: 
$$(A(P',s.P, s\si) \l(s))^*= A(P, t.P, \si) \lambda(s^{-1}).$$
Alors les ŽgalitŽs prŽcŽdentes jointes ˆ  l'ŽgalitŽ (\ref{C00}),  montrent, aprs  que l'on ai pris les inverses:
 \beq \label{BmB}   B(s. P, P', s\si )^*=  B(t.  P', P, m\si) \si'(m^{-1}) \eeq
 Par conjugaison par $m$, on a: 
 $$B(t.  P', P, m\si)  \si'(m^{-1})=  \si'(m^{-1}) B(t.P, P, \si) $$
   Mais  d'aprs l'analogue de [D4],  (7.20),  dont la preuve, par transport de structure,  est identique et o $\s$ et $\t$ sont triviaux car $P$ et $P'$ sont standard, on a:
  $$  \si'(m^{-1}) B(t.P, P, \si) = B(P',s. P, s\si)$$
  Joignant les deux ŽgalitŽs prŽcŽdentes ˆ (\ref{BmB}), on obtient:
  $$B(s.P, P', s\si))^*= B(P',s.P, s\si)$$
   Echangeant le role de $P$ et $P'$, on en dŽduit le ThŽorme. 
\qed
\section{ Formule de Plancherel}
\setcounter{equation}{0}
\subsection{Equation fonctionnelle pour les intŽgrales de Jacquet}
Soit $\Theta$ l'ensemble des couples $(P, \O)$ o $P=MU$ est un sous-groupe parabolique anti-standard de $G$ et $\O$ l'orbite inertielle unitaire d'une reprŽsentation lisse irrŽductible et  de carrŽ intŽgrable de $M$.
\begin{prop}\label{Efonc}
 Soit  $P= MU$, $P'= M'U'$ deux sous-groupes paraboliques anti-standard de $G$  tels que $M$ et $M'$ soient conjuguŽs et $s\in W(M'\vert G\vert M)$. Soit $(P,\O)\in \Theta$ et  $\phi \in \fhi $.  
 On a l'ŽgalitŽ de fonctions sur $\O$:
 $$E^G_{P'} (C^0(s, P',P, \si) \phi(\si) )= E^G_P (\phi(\si) ). $$
\end{prop}
 \dem 
 Soit $\eta' \in Wh(s\si)$. D'aprs la dŽfinition des matrices $B$, on  a l'ŽgalitŽ de fonctions sur $\O$:
 $$ \xi(s.P, s\si, B(s.P, P', s\si, \eta')  = \xi( P', s\si,\eta') \circ A(P',s.P, s\si).$$ 
 Mais d'aprs (\ref{whp}), on a   $Wh(s.P, s\si)=Wh(\si)$ et l'ŽgalitŽ.:
 $$\xi(s.P, s\si, B(s.P, P', s\si)\eta')= \xi(P, \si, B(s.P, P', s\si) \eta') \lambda(s^{-1}).$$
 Donc on a:
 \beq \label{xisiB} \xi(P, \si, B(s.P, P', s\si) \eta') = \xi( P', s\si,\eta') \circ A(P',s.P, s\si)\circ \lambda(s). \eeq
 Posant $ \eta= B(s.P, P', s\si) \eta'\in Wh(s.P, s\si)= Wh(\si) $, on a, lorsque $B(s.P, P', s\si)^{-1}  $ est dŽfini:
 $$ \xi(P, \si, \eta)= \xi( P', s\si,B(s.P, P', s\si)^{-1}\eta)  A(P',s.P, s\si)\circ \lambda (s).$$
 Utilisant la dŽfinition des intŽgrales de Jacquet et la dŽfinition de $C^0(s,P',P,\si)$, on en dŽduit la Proposition.\qed 
 \begin{rem}
 Notons $A(s, P', P, \si):=A(P',s.P, s\si) \lambda(s) $ l'opŽrateur d'entrelacement entre $i^G_P\si$ et $i^G_{PÔ} s\si$, lorqu'il est dŽfini.  Comme dans la dŽfinition des matrices $B$, on voit qu'il existe une unique application linŽaire, notŽe $B(s^{-1} , P, P', \si)$, entre $Wh( s\si)$ et $Wh(\si)$  telle que:
 $$\xi( P,\si,  B(s^{-1}, P, P',s \si)\eta')= \xi(P', s\si, \eta')\circ A(s, P',P, \si) , \eta' \in Wh(s\si). $$ 
On dŽduit de (\ref{xisiB}), l'ŽgalitŽ:
$$B(s.P, P', s\si) =B(s^{-1}, P, P',s \si) .$$
Cela permet une dŽfinition des fonctions $C^0$ sans avoir recours, pour les intŽgrales de Jacquet,  ˆ des sous-groupes paraboliques autres qu'anti-standard. 
\end{rem}
 \subsection{Equation fonctionnelle pour la TransformŽe de Fourier-Whittaker}
\begin{prop}\label{eqfoncF}
Soit $f \in \CCU$.\\
(i)   Il existe un ensemble fini de couples $(P, \O)\in \Theta$ tels que  que $ \hat{f}_{P, \O}$ soit non nul.\\ (ii) Soit  $P= MU$, $P'= M'U'$ deux sous-groupes paraboliques antistandard de $G$ et $s\in W(M'\vert G \vert M)$
tels que $M$ et $M'$ soient conjuguŽs. Soit $(P, \O) \in \Theta$.
On a l'ŽgalitŽ, pour tout $(\si,E)$ objet de $\O$:
$$ \hat{f}(P', s\si)   = C^0(s, P', P, \si)  \hat{f}(P, \si).$$
\end{prop}
\dem
(i)  Soit $H$ un sous-groupe compact ouvert de $G$ fixant $f$. (i) rŽsulte  du fait  qu'il n'y a qu'un nombre fini de couples $(P, \O)\in \Theta$  tels que pour un objet de $\O$, $(\si,E)$, l'espace $i^G_P E$ contienne un vecteur non nul invariant par $H$ ( cf [W], ThŽorme VIII.1.2).\\
 (ii) Soit $ \phi \in Wh(\si) \otimes i^G_PE$. On a par dŽfinition de $\hat{f}$:
 $$(\hat{f}(P, \si),\phi)= \int_{U_0\backslash G} f(g) \overline{  E^G_P(\phi)(g)}dg.$$
 Tenant compte de la Propostion \ref{Efonc}, on en dŽduit:
 $$(\hat{f}(P, \si),\phi)= \int_{U_0\backslash G} f(g) \overline{  E^G_{P'} (C^0(s,P',P, \si)\phi)(g)}dg.$$
 De la dŽfinition de $\hat{f}$, on en dŽduit:
  $$(\hat{f}(P, \si),\phi)= (\hat{f}(P', s\si),  (C^0(s,P',P, \si)\phi).$$
  Mais, l'Assertion A tant maintenant dŽmontrŽe, l'adjoint de $(C^0(s,P',P, \si)$ est Žgal ˆ $(C^0(s,P',P, \si)^{-1}$, daprs la Proposition \ref{Chol}. Donc:
$$ (\hat{f}(P,\si),\phi)=  ((C^0(s,P',P, \si)^{-1} \hat{f}(P', s\si), \phi).$$
Comme cela est vrai pour tout $\phi$, on en dŽduit l'ŽgalitŽ voulue. \qed
\subsection{Espace ${\cal S}^{inv}$ et applications $\cal{F}$ et $\cal W$}
 On dŽfinit:
 $${ \cal S} = \oplus_{(P, \O)\in \Theta} C^{\infty}(\O, Wh \otimes i^G_P).$$
Si $P=MU$ est un sous groupe parabolique anti-standard de $G$, on note $p(P)$ le nombre de sous-groupes paraboliques anti-standard dont le LŽvi standard est conjuguŽ ˆ $M$. Si $I$ est un ensemble fini, on  note $\vert I\vert$ son cardinal. On dŽfinit:
\beq \label{cP} c(P):= p(P) \vert W(M\vert G\vert M)\vert .\eeq 
 On munit $\cal{S}$  du produit scalaire qui fait de cette somme une somme directe orthogonale et qui,  pour tout $(P, \O)\in \Theta$, $\phi, \phi' \in C^{\infty}(\O, Wh \otimes i^G_P)$ vŽrifie:
 \beq\label{prodS} (\phi, \phi')=  c(P)^{-1}\int_\O \mu(\si) (\phi(\si), \phi(\si')) d\si.\eeq
 On note ${\cal S}^{inv}$ l'ensemble des ŽlŽments $(\phi_{P,\O})_{(P, \O) \in \Theta}$ tels que, si $(P=MU, \O)\in \Theta$, $ P'=M'U'$ est  un  sous-groupes parabolique anti-standard de $G$  telsque $M$ et $M'$ soient conjuguŽs et $s \in W(M'\vert G\vert M)$, on a:
 $$\phi_{P', s\O} ( s\si) = C^{0}(s, P',P,\si)\phi_{P,\O} (\si),$$
 pour tout $\si$ objet de $\O$.\\
 Pour $f\in \CCU$, on note ${\cal F}f$, la famille $(\hat{f}_{P,\O} )_{(P, \O) \in \Theta}$. D'aprs les Propositions  \ref{fpoo} et \ref{eqfoncF}, c'est un ŽlŽment de ${\cal S}^{inv} $.\\
Tenant compte de la Proposition \ref{waveS}, on  note $\cal{W}$, l'application de $\cal S$ dans $\CCU$ qui ˆ 
$(\phi_{P,\O})_{(P, \O)\in \Theta}$ associe 
$Ê\sum_{(P, \O) \in \Theta} c(P)^{-1} f_{\phi_{P,\O}}$, la somme ne comportant qu'un nombre fini de termes non nuls. 
\subsection{ InjectivitŽ de $\cal{F}$} Soit $(\pi, V)$ une reprŽsentation lisse de $G$,  $\xi \in Wh(\pi)$. On suppose que $\pi$ est tempŽrŽe de sorte que $\xi$ est Žgalement tempŽrŽe.  Soit  $ f \in \CCU$.  On note  $f^*$ la fonction sur $G$ dŽfinie par  $ f^*(g) = \overline{f(g^{-1})}$ pour $ g \in G$. On dŽfinit $\pi' (f^{*}) \xi  \in \check{V}$, par:
\beq \label{pifxi}  \langle \pi'(f^*) \xi, v\rangle  =  \int_{G/U_0} f^*(g)\langle  \pi'(g) \xi, v\rangle  dg. \eeq
l'intŽgrale Žtant convergente puisque $\xi$ est tempŽrŽe. 
La reprŽsentation rŽgulire droite de $G$ dans $\lu$ se dŽcompose en une intŽgrale hilbertienne    de reprŽsentations unitaires irrŽductibles de $G$,  $(\int_Z^\oplus\pi_z d\mu(z), \int_Z^\oplus H_z d\mu(z))$. On note, pour $z\in Z $, $(\pi_z^\infty, H_z^{\infty})$  la reprŽsentation lisse de $G$, $\pi_z^\infty$, dans l'espace $H_z^{\infty} $,  des vecteurs lisses de $H_z$,  i.e. fixŽs par un sous-groupe compact ouvert de $G$. 
\begin{prop}  \label{bernstein} Pour $\mu$-presque tout $z\in  Z$, il existe un morphisme de $G$-modules, $\beta_z$, entre $\CCU$ et  $H_z^\infty$ et il existe  $\xi_z \in Wh( \pi^\infty _z)$ tels que:
\\(i)  Pour tout $ f \in \CCU$, $f= \int_Z^\oplus\beta_z( f)  d\mu(z).$\\
(ii) Pour  $\mu$-presque tout $z\in Z$, $\pi_z$ est tempŽrŽe.
\\(ii)  Pour $\mu$-presque tout $z\in Z$, on a $$(\beta_z ( f), v)= 
 \overline{\langle \pi'_z( f^{*} )\xi_z,  v \rangle }, f \in \CCU, v\in H^{\infty}_z.$$
\end{prop}
\dem  D'aprs [B1], section 3.6 et section 1.4, pour montrer (i) il suffit de dŽmontrer que l'injection de $\CCU$ est ''fine'' dans le sens de l.c., section 1.4.
On introduit pour cela la fonction $w$ sur $U_0\backslash G$:
$$w(g) = (1 + \Vert H_0(m_0(g))\Vert)^d $$
o $d=2 dim\>  A_0$. Soit $S_w:= \{ f\in L^2( U_0\backslash G, \chi)\vert \int_{U_0 \backslash G} \vert f(g)\vert ^2w(g)dg < \infty \}$.  
Alors $\CCU$ est un sous espace de $S_w$ et l'injection est continue d'aprs (\ref{intcv}). 
Il suffit donc de montrer que l'injection de $S_w$ est ''fine'' (cf. l.c section 1.4) . Pour cela utilisant la section 3.6  et le ThŽorme 3.4 de l.c., il suffit de voir que $w$ est un poids sommable dans le sens de l.c.,  section 3.2. 
Pour cela on introduit $\Lambda= \Lambda(A_0)F_0$, o $F_0$ est  un sous-ensemble fini de $M_0$ tel que $U_0 A_0F_0K$. Alors  c'est un ''net'' au sens de l.c. section 3.2.
On vŽrifie facilement que:  $$\sum_{g\in \Lambda} w(g)^{-1}< \infty.$$
Donc $w$ est bien sommable. Comme on l'a dit, (i) en rŽsulte. \\
(ii) a ŽtŽ dŽmontrŽ dans [D4], ThŽorme 7. Ceci implique que pour  $\mu$-presque tout $z\in Z$, $\xi_z$ est tempŽrŽe.\\
On dŽduit (iii) de (i) comme dans la Proposition 15 de [D4] . \qed 
\begin{cor}Soit $f\in \CCU$, si pour toute reprŽsentation unitaire irrŽductiblle tempŽrŽe lisse de $G$ et $\xi \in Wh(\pi)$, $ 
 \pi'(f^*) \xi=0$, alors $f$ est nulle. 
 \end{cor}
 \begin{prop}
$\cal{F}$ est une application linŽaire injective de $\CCU$ dans ${\cal S}^{inv}$. 
\end{prop}
\dem
Soit $f\in \CCU$ tel que ${\cal F} f=0$.  Comme dans la preuve de la Proposition 11 de [D4], on voit que,  pour tout $(P, \O) \in \Theta$, $(\si,E) $ objet de $\O$ et $ \xi\in Wh(\si)$:
$$(i^G_P\si )'(f^*) \xi(P,\si, \eta)=0$$
Comme toute reprŽsentation lisse unitaire irrŽductible et  tempŽrŽe de $G$, $(\pi,V)$, est un facteur direct d'une reprŽsentation $i^G_P \si$ et  que l'espace $Wh(i^G_P \si)$ est Žgal ˆ $\{\xi (P, \si, \eta)\vert \eta \in Wh(\si)\}$, on en dŽduit:
$$ \pi'(f^{*}) \xi =0, \xi \in Wh(\si).$$
La proposition rŽsulte alors du Corollaire prŽcŽdent.  \qed 
\subsection{Formule de Plancherel} 
\begin{theo}
L'application $\cal{F}$ est une bijection de $\CCU$ sur ${\cal S}^{inv}$, dont l'inverse est donnŽ par la restriction de $\cal{W} $ ˆ ${\cal S}^{inv}$. 
En particulier: 
$$f= \sum_{(P=MU,\O)\in \Theta} c(P)^{-1} \int_\O \mu(\si)E^G_P ( \hat{f}_{P, O} (\si) )d\si  , f \in \CCU .$$
 De plus $ \cal{F}$  est une isomŽtrie, lorsque $\CCU$ est muni du produit scalaire $(.,.)_G$ et ${\cal S}$ du produit scalaire dŽfini en (\ref{prodS}).
\end{theo} 
\dem
On a dŽjˆ montrŽ que $\cal{F}$ est injective. 
Montrons que la restriction de $\cal{FW}$ ˆ ${\cal S}^{ inv}$ est l'identitŽ. 
Soit $\Phi= (\phi_{P,\O})_{(P, \O)\in \Theta} \in {\cal S}^{inv}$. Calculons ${\cal FW}\Phi$. 
Soit $(P=MU, \O)$, $(P=M4U4)', \O')$ deux ŽlŽments de $\Theta$.D'aprs le Corollaire  D'aprs le Corollaire de la Proposition \ref{prodscalpaq}, $\hat{f}_{P',\O'}$ est nul si $M$ et $M'$ ne sont pas conjuguŽs. Si $M$ et $M'$ sont conjuguŽs, d'aprs le Thorme \ref{fphihat},  et le fait que $\Phi \in {\cal S}^{inv}$, on a pour $\sigma_1$ objet de $\O'$:  
 $$ (f_{\phi_{P, \O}}) \hat{}_{ P', \O'} (\si_1) =
 \sum_{s \in W(M'\vert G\vert M), s^{-1} \si_1 \in \O } C^{0}(s, P',P, s^{-1} \si_1)\phi_{P, \O} (s^{-1} \si_1).$$
De la dŽfinition de ${\cal S}^{inv}$, on en dŽduit , que, pour $\sigma_1$ objet de $\O'$ on 
 $$  (f_{\phi_{P, \O}}) \hat{}_{ P', \O'} (\si_1) = \vert \{s \in W(M'\vert G\vert M), s^{-1} \si_1 \in \O\} \vert \phi_{P', \O'}(\si_1).$$
 Appliquons cela pour calculer $({\cal W}\Phi) \hat{}_{P, \O}$.
 D'aprs ce qui prŽcde et en Žchangeant le r\^ole de $P$ et $P'$, $c(P)({\cal W}\Phi)\hat{ } _{P, \O}$
 est la somme sur $(P'=M'U', \O')\in \Theta$, avec $M$ et $M'$ conjuguŽs,    de $\vert \{s \in W(M'\vert G\vert M), s^{-1} \O' =\O\}   \vert \phi_{P, \O}$.
 Mais l'ensemble des couples $( s, \O')$ tels  $s^{-1} \O' =\O)$ est Žgal ˆ $ \{ (s, s\O) \vert s \in W(M'\vert G\vert M)\}$, dont le nombre d'ŽlŽment est Žgal ˆ $\vert W(M'\vert G \vert  M)\vert$. \\
 Montrons que si $M'$ et $M$ sont conjuguŽs  par un ŽlŽment de $\overline{W}^G$, on a:  \beq \label{WW}  \vert W(M'\vert G \vert  M)\vert=  W(M\vert G \vert  M)\vert .\eeq 
En effet si $M'= x.M$ avec $x$ ŽlŽment de $\overline{W}^G$,  on a $ \overline{W}^{M'} = x. \overline{W}^M$  ainsi que, avec les notations de la section \ref{sssec},  
$ \overline{\cal{W}}(M'\vert G\vert M)= x\overline{\cal{W}}(M\vert G\vert M) $.  Alors (\ref{WW}) en rŽsulte immŽdiatement. En sommant sur les $(P',\O')$ et tenant compte de la dŽfinition de ${\cal W}$, on voit que
$$ ({\cal W}\Phi)\hat{} _{P, \O}=p(P) \vert W(M \vert G \vert  M)\vert c(P)^{-1} \phi_{P, \O} .$$ 
 Finalement :
$$({\cal W}\Phi)\hat{} _{P, \O}=  \phi_{P, \O}. $$
Donc $${\cal FW} \Phi= \Phi, \Phi \in {\cal S}^{inv}. $$
Ceci montre que $\cal F$ est surjective sur ${\cal S}^{inv}$ . Comme on  a vu que $\cal F$ Žtait aussi injective, c'est bien une bijection sur  entre $\CCU $ et ${\cal S}^{inv}$ . Ce qu'on vient de montrer prouve Žgalement que l'inverse de $\cal F$ est Žgal ˆ la restriction ˆ ${\cal S}^{inv}$ de $\cal W$.\\
Il ne reste plus qu'ˆ montrer que $\cal W$ restreint ˆ ${\cal S}^{inv}$   est isomŽtrique. \\Soit $\Phi= (\phi_{P,\O})_{(P, \O)\in \Theta} , \Phi' = (\phi'_{P,\O})_{(P, \O)\in \Theta} \in {\cal S}^{inv}$. 
Alors: \beq \label{calww}({\cal W} \Phi, {\cal W} \Phi')=  \sum_{(P,\O), (P', \O') \in \Theta } c(P)^{-1}c(P')^{-1} ( f_{\phi_{P,\O}}, f_{\phi'_{P', \O'}})_G. \eeq
D'aprs la Proposition \ref{prodscalpaq}, $(f_{\phi_{P,\O}}, f_{\phi'_{P', \O'}})_G$  est nul si $M$ et $M'$ ne sont pas conjuguŽs par $W^G$ et sinon:
$$(f_{\phi_{P,\O}}, f_{\phi'_{P', \O'}})_G=  \sum_{s\in W(M'\vert G\vert M), s\O= \O'} \int_{\O'} \mu(\si_1)(C^{0}(s, P',P, s^{-1} \si_1 )\phi_{P, \O}( s^{-1}\si_1), \phi'_{P', \O'}(\si_1)) d \si_1$$
Comme $ \phi \in   {\cal S}^{inv}$, on a :
$$(f_{\phi_{P,\O}}, f_{\phi'_{P', \O'}})_G= \vert \{s\in W(M'\vert G\vert M), s\O= \O'\} \vert (\phi_{P', \O'}, \phi'_{P', \O'}). $$
Si $P$ et $P'$ sont associŽs, on a $c(P)=c(P')$.
Reportant dans (\ref{calww}) et regroupant les termes correspondants ˆ $\phi'_{P', \O'}$, on trouve:
$$({\cal W} \Phi, {\cal W} \Phi') =\sum_{(P', \O')\in \Theta}p(P') \vert W(M'\vert G\vert M') c(P')^{-2}  (\phi_{P', \O'}, \phi'_{P', \O'}).$$
 Soit encore:$$({\cal W} \Phi, {\cal W} \Phi') =\sum_{(P', \O')\in \Theta} c(P')^{-1}  (\phi_{P', \O'}, \phi'_{P', \O'}).$$
 Finalement 
 $$({\cal W} \Phi, {\cal W} \Phi')= ( \Phi, \Phi').$$
 Donc l'inverse de $\cal F$ est bien isomŽtrique.Il en va de m\^eme de $\cal F$. Ceci achve la preuve du ThŽorme. \qed
 \begin{rem} 
 La Proposition  \ref{fphipind} et La surjectivitŽ de ${\cal W}$ montre que pour tout $f \in\CCU$ et tout sous-groupe parabolique anti-standard de $G$, $P$, l'application $f^P$ est un ŽlŽment de $\CCMU$.
 \end{rem}

\section{RŽfŽrences}

 \noindent [A]  Arthur J., A local trace formula. Publ. Math.  Inst. Hautes \' Etudes Sci.   No. 73  (1991), 5--96.

 \noindent[B1] Bernstein, J.,  On the support of the Plancherel measure,  J. Geom. Phys.  5  (1988),   663--710 (1989). 
 
 \noindent[B2]  Bernstein J., Representations of p-adic groups. Lectures given at Harvard
University, Fall 1992, Notes by K. E. Rummelhart.
 
 \noindent [B3] Bernstein J.,  Second adjointness theorem for representations of reductive p-adic groups, unpublished
manuscript. 

\noindent [BZ] Bernstein J., Zelevinsky V., Induced representations of reductive $p$-adic groups I,  Annales de l'Ecole Normale SupŽrieure 10 (1977) 441-472.

\noindent [BlD] Blanc P., Delorme P.,  Vecteurs distributions H-invariants de
reprŽsentations induites pour un espace symŽtrique rŽductif $p$-adique
G/H, Ann. Inst. Fourier, 58 (2008)  213--261.

\noindent[BoT] Borel A., Tits J., Groupes rŽductifs, Publ. Math.  Inst. Hautes \' Etudes Sci. 27 (1965) 55--150. 

\noindent [BrD] Brylinski, J-L.,  Delorme, P.
Vecteurs distributions $H$-invariants pour
 les sŽries principales gŽnŽralisŽes d'espaces symŽtriques rŽductifs et prolongement mŽromorphe,
d'intŽgrales d'Eisenstein. Invent. Math.
109 (1992) 619--664. 

\noindent [Bu] Bushnell, C., Representations of reductive $p$-adic groups: localization of Hecke algebras and applications.  J. London Math. Soc. 63  (2001),  no. 2, 364--386. 
 
 \noindent [BuHen] Bushnell, C.,  Henniart  G., 
Generalized Whittaker models and the Bernstein center.
Amer. J. Math. 125 (2003),  513--547.

\noindent [C] Casselman W.,  Introduction  to the theory of admissible representations of $p$-adic reductive groups,
http://www.math.ubc.ca/$\tilde{}$ cass/research.html. 

\noindent[CS]  Casselman, W., Shalika, J.,  The unramified principal series of $p$-adic groups. II. The Whittaker function. Compositio Math. 41 (1980), 207--231.

\noindent [DeliB] Deligne P.,  Le ``centre'' de Bernstein rŽdigŽ par Pierre Deligne. Travaux en Cours,  Representations of reductive groups over a local field,  1--32, Hermann, Paris, 1984.

\noindent [D1] Delorme  P., Espace des coefficients de reprŽsentations admissibles d'un groupe
rŽductif $p$-adique,  131--176, Progr. Math., 220, Birkhauser Boston, Boston, MA, 2004.

\noindent [D2] Delorme P., The Plancherel formula on reductive symmetric spaces from the point of view of the Schwartz space.  Lie theory,  135--175, Progr. Math., 230, Birkhauser Boston, Boston, MA, 2005. 

\noindent[D3] Delorme P., Constant term of smooth $H_\psi$-invariant functions,  Trans. Amer. Math. Soc.  362  (2010),  933--955. 

 \noindent[D4] Delorme P., ThŽorme de Paley-Wiener pour les fonctions de Whittaker sur un groupe rŽductif $p$-adique (deuxime version,uniquement sur ma page web, celle-ci Žtant boquŽe sur le site HAL, et n'apparait pas encore sur arXiv).
 
  \noindent [H] Heiermann, V.,  Une formule de Plancherel pour l'algbre de Hecke d'un groupe rŽductif $p$-adique. Comment. Math. Helv. 76 (2001),  388--415.
 
 \noindent [Hu] Humphreys J. E, Linear algebraic groups, Graduate Text In Math. 21, Springer, 1981.
 
 \noindentÊ[K] Knapp A., Representation theory of semisimple groups. An overview based on examples. Reprint of the 1986 original. Princeton Landmarks in Mathematics. Princeton University Press, Princeton, NJ, 2001
 
 \noindent [Ma] , Matringe  N., Derivatives and asymptotics of Whittaker functions, arXiv:1004.1315. 
     
\noindent [M]  Michael E.,  Selected selections theorems, Amer. Math. Monthly 63 (1956) 223-283.

\noindent [R] Rodier, F., Modles de Whittaker des reprŽsentations admisssibles des groupes rŽductifs p-adiques quasi-dŽployŽs, manuscript  non publiŽ.

\noindent [Sh] Shahidi F.,  On certain $L$-functions.  Amer. J. Math.  103  (1981)  297--355.

\noindent [T] Tits J., ReprŽsentations linŽaires irrŽductibles  d'un groupe rŽductif sur un corps quelconque, J.  Reine Angew.   247 (1971),  196--220. 

  \noindent [W] Waldspurger J.-L.,  La formule de Plancherel pour les groupes
$p$-adiques (d'aprs Harish-Chandra), J. Inst. Math. Jussieu  2  (2003),  235--333.

 \noindent [Wall] Wallach, N., Real reductive groups. II. Pure and Applied Mathematics, 132-II. Academic Press, Inc., Boston, MA, 1992.
 
  \noindent [War] Warner, G., Harmonic analysis on semi-simple Lie groups. I. Die Grundlehren der mathematischen Wissenschaften, Band 188. Springer-Verlag, New York-Heidelberg, 1972.
 \\\\
 Patrick Delorme\\Institut de Math\'ematiques de Luminy, UMR 6206 CNRS,
Universit\'e  de la M\'editerran\'ee, 163 Avenue de Luminy,
  13288 Marseillle Cedex 09, France
\\delorme@iml.univ-mrs.fr
\\
L'auteur a bŽnŽficiŽ du soutien du programme ANR-BLAN08-2-326851 pendant l'Žlaboration de ce travail.

 \end{document}